\documentclass[11pt,a4paper,fleqn,reqno]{amsart}
\usepackage[T1]{fontenc}
\usepackage[utf8]{inputenc}
\usepackage[english]{babel}
\usepackage{graphicx}
\usepackage[font=small]{quoting}
\usepackage{caption}
\usepackage{enumitem}

\usepackage{float}
\usepackage{mathrsfs}
\usepackage[top=1.0in,bottom=1.0in,left=1.in,right=1.in]{geometry}
\usepackage{csquotes}
\usepackage[url=false,doi=false,isbn=false,giveninits=true,
  style=numeric, 
  maxbibnames=99]{biblatex}
\usepackage{mathtools}
\usepackage{verbatim}
\usepackage{color}
\usepackage{microtype}
\usepackage{amsmath,amsthm,amsfonts,amssymb}
\usepackage{comment}
\usepackage{etoolbox}
\usepackage{xurl}
\usepackage{hyperref}
\hypersetup{colorlinks,linkcolor={blue},citecolor={blue},urlcolor={blue},draft
}
\newtheorem{thm}{Theorem}[section]
\newtheorem{lem}[thm]{Lemma}
\newtheorem{cor}{Corollary}[section]
\newtheorem{prop}{Proposition}[section]
\newtheorem{defn}{Definition}[section]
\theoremstyle{definition}
\newtheorem{rem}{Remark}[section]
\newtheorem*{rem*}{Remark}
\theoremstyle{remark}

\theoremstyle{definition}

\def\n{\nabla}

\def\tr{\mathrm{tr}}

\def\<{\langle}
\def\>{\rangle}

\def\RR{\mathbb{R}}

\def\SS{\mathbb{S}}

\def\HH{\mathbb{H}}

\def\wh{\widehat}

\def\wt{\widetilde}

\def\R{\mathbb{R}}

\DeclareMathOperator{\dive}{div}

\DeclareMathOperator{\dist}{dist}

\patchcmd{\abstract}{\scshape\abstractname}{\textbf{\abstractname}}{}{}
\makeatletter
\def\@makefnmark{}
\makeatother

\numberwithin{equation}{section}
\numberwithin{defn}{section}

\title[Neumann problem in hyperbolic space]{Neumann problem for Hessian equation and geometric inequalities in hyperbolic space}



\begin{document}

\author[L. Weng]{Liangjun Weng}
\address[L. Weng]{Centro di Ricerca Matematica Ennio De Giorgi, Scuola Normale Superiore, Pisa, 56126, Italy}
\email{\href{mailto:liangjun.weng@sns.it}  {liangjun.weng@sns.it}}

\begin{abstract}
We establish the existence of solutions to the Neumann problem for the Hessian equation on a smooth, bounded strictly static convex domain in hyperbolic space. The boundary condition determines a unique compatibility constant, and admissible solutions are unique up to a multiple of the static potential. As an application, we give a direct analytic proof of some weighted geometric inequalities for the static convex domain in hyperbolic space. 
\end{abstract} 
 
\maketitle

 
\section{Introduction}

The use of elliptic boundary value problems in convex geometry goes back, in particular, to Reilly's integral formula and his interpretation of the Minkowski inequality through a Neumann problem \cite{Reilly, ReillyNeumann}. In Euclidean space, elliptic methods include the use of
Hessian equations in quermassintegral inequalities
\cite{trudinger-1994} and the Poisson equation approach to the isoperimetric
inequality described by Cabr\'e \cite{Cabre2008}, see also \cite{brendle-2021-jams,xia-zhang-2017} etc.  
Lions--Trudinger--Urbas \cite{LTU} developed the Neumann
theory for Monge--Amp\`ere equations, including solvability
with a unique additive boundary constant and uniqueness
of solutions modulo constants. For Hessian equations with Neumann boundary,
Trudinger treated balls \cite[Theorem~2]{Trudinger1987},
and Ma--Qiu established solvability on uniformly convex
domain in 
\cite[Theorem~1]{MQ}. Qiu--Xia extended the formulation
with an unknown boundary constant to general Hessian
equations \cite[Theorem~1]{QX} and applied the resulting
solutions to a family of Euclidean Alexandrov--Fenchel
inequalities \cite[Theorem~2]{QX}.
These results complement the classical mixed-volume convex body theory \cite{Schneider} and hypersurface curvature flow proofs of quermassintegral inequalities, including the extension to star-shaped $k$-convex domains by Guan--Li in \cite{GuanLi2009}; among many others.

The hyperbolic counterpart is governed by a static potential. Let $(\HH^n,g)$ have constant sectional curvature $-1$, and  $\Omega\Subset\HH^n$ be a bounded smooth domain. Fix a point $p$, and define $V(x)\coloneqq \cosh d(p,x)$. For $u\in C^2(\bar\Omega)$,  define 
\begin{eqnarray}\label{eq:static-pair-intro}
A[u]\coloneqq \nabla^2u-ug,\qquad
B_Vu\coloneqq u_\nu-\frac{V_\nu}{V}u.
\end{eqnarray}
Both expressions are invariant under $u\mapsto u+sV$, for all $s\in\RR$. These quantities occur naturally in the weighted Reilly formulas of Qiu--Xia \cite{QXReilly} and Li--Xia \cite{LiXia}. At the linear level, Xia \cite{Xia} used the solution to the   Neumann problem for linear PDEs and proved the sharp weighted Minkowski inequality
\begin{eqnarray}\label{eq:xia-intro}
\left(\int_\Sigma V d\sigma_g \right)^2
\ge n\left(\int_\Omega V d\mu_g \right)
\left(\int_\Sigma VH d\sigma_g \right)
\end{eqnarray}
under the static convexity assumption, namely: $h\ge(\frac{V_\nu}{V})g_\Sigma$. Here $\Sigma=\partial\Omega$, $h$ is the second fundamental form with respect to the outward normal and $H$ is the mean curvature of $\Sigma$, see Eq. \eqref{eqn:H-defn}. Equality holds for every geodesic ball, without requiring its center to be $p$. This linear theory suggests asking whether
$\sigma_k(A[u])=f$ admits the same boundary condition,
with its constant determined as part of the solution,
and what further information these nonlinear Neumann potentials
provide about geometric inequalities.

There is now a substantial geometric context and progress along this direction. The Heintze--Karcher comparison principle \cite{HeintzeKarcher}, its weighted warped-product form due to Brendle \cite{Brendle}, and the Minkowski-type inequality of Brendle, Hung and Wang \cite{BrendleHungWang} show the importance of static weights in negatively curved geometry. Wang and Xia established sharp comparisons between hyperbolic quermassintegrals on horospherically convex domains \cite{WangXia}; Ge, Wang and Wu proved further hyperbolic quermassintegral inequalities \cite{GeWangWu}. Quermassintegral-preserving curvature flows were studied by
Andrews--Wei \cite{AndrewsWei}, and volume-preserving flows
by Andrews--Chen--Wei \cite{AndrewsChenWei}.
Locally constrained flows and their geometric applications
were developed in \cite{ScheuerXia,HuLiWei}, among many others; see also the survey by Guan--Li \cite{Guan-Li2021-JMS}. 
 
 For inequalities involving the static weight, Hu--Li proved weighted inequalities under static convexity \cite{HuLi}, and Wu established weighted comparisons for a class of
nondecreasing convex weights satisfying an additional
differential inequality \cite[Theorem~1.1]{Wu}. Hu, Li, Wan and Xu used the relation between static convexity and Euclidean convexity under spatial projection to obtain affine isoperimetric inequalities \cite{HuLiWanXu}. Li--Xu \cite{LiXu} removed convexity assumptions from the weighted isoperimetric comparison of Scheuer--Xia \cite[Theorem 1.6]{ScheuerXia}. For more related results, see the references therein. The Neumann applications below concern the weighted quantities in \eqref{eq:xia-intro}.

Nonlinear Neumann problems on Riemannian manifolds were studied
by Cherrier \cite{Cherrier-1984}. 
For $\det(\nabla^2u-ug)$, Hanani
\cite[Theorems~1 and~3]{Hanani} proved solvability on strictly
convex Riemannian manifolds under supersolution hypotheses, both for prescribed Neumann data and for boundary data strictly increasing in $u$ with respect to his inward normal convention. 
 Jin--Li--Li \cite{Jin-Li-Liyanyan2007} studied the fully
nonlinear Yamabe problem with Neumann-type boundary condition on a manifold with umbilic boundary.
The static pair in
\eqref{eq:static-pair-intro} presents two linked difficulties: its
boundary zeroth-order coefficient has the critical sign associated
with the kernel generated by $V$, and negative curvature contributes
an unfavorable term to the directional Hessian estimate. The
estimates below address both features on nonumbilic domains.

The main theorem assumes \textit{strict static convexity} relative to $V$
and places the pole $p$ in $\Omega$, namely, 
\begin{eqnarray}\label{eqn:static-convexity}
 h-\frac{V_\nu}{V}g_\Sigma>0.
\end{eqnarray}For the non-strict positive support formulation of
static convexity, see
\cite[Definition~3]{brendle-wang-mu-tao2014}. 
The projection correspondence of
\cite[Corollary~2.1]{HuLiWanXu} implies, when $p\in\Omega$,
that $V_\nu>0$ and $h>0$; see
Lemma~\ref{lem:static-convex-domain}.
Thus neither starshapedness nor geodesic convexity needs
to be imposed separately. In this setting,
\eqref{eqn:static-convexity} is equivalent to
$h>(V_\nu/V)g_\Sigma>0$.
By contrast, Xia's weighted Minkowski inequality
\cite[Theorem~1.1]{Xia} assumes the non-strict inequality
$h\ge(V_\nu/V)g_\Sigma$ and does not require $p\in\Omega$.
For an interior pole, horospherical convexity
$h\ge g_\Sigma$ implies \eqref{eqn:static-convexity},
since $|\nabla V|<V$.
The arbitrary-pole extension below assumes $h>0$ and
\eqref{eqn:static-convexity} separately.
 A function is $k$-admissible if the
eigenvalues of $\nabla^2u-ug$ belong to the G\aa rding cone $\Gamma_k$
on $\bar\Omega$, cf. Definition \ref{def-k-admissible}.

\begin{thm} \label{thm:existence}
Let $\Omega\Subset\HH^n$, $n\ge2$, be a connected bounded domain with
$C^{4,\alpha_0}$ boundary, $0<\alpha_0<1$. Fix $p\in\Omega$, put
$V=\cosh d(p,\cdot)$, and assume that $\Omega$ is strictly static 
convex relative to $V$. 
Let $1\le k\le n$, and 
$f\in C^2(\bar\Omega)$ be strictly positive and
$\varphi\in C^3(\Sigma)$. Then there exists a unique constant $\lambda \in\RR$ for which
\begin{eqnarray}\label{eq:main-problem}
\begin{cases}
\sigma_k(\nabla^2u-ug)=f(x) \qquad &\text{ in }\Omega,\\
u_\nu=\dfrac{V_\nu}{V}u+\lambda +\varphi(x) \qquad &\text{ on }\Sigma
\end{cases}
\end{eqnarray}
has a $k$-admissible solution $u\in C^{3,\alpha}(\bar\Omega)$ for
some $0<\alpha\le\alpha_0$. Any two solutions differ by a multiple
of $V$.
\end{thm}

\begin{rem}\ 
\begin{enumerate}[label=(\arabic*), leftmargin=*, labelindent=0pt]
    \item 

Static convexity in \eqref{eqn:static-convexity} is relative to the
specified pole $p$. Since $p\in\Omega$,
Lemma~\ref{lem:static-convex-domain} gives $V_\nu>0$, so static
convexity implies strict geodesic convexity.  However, the static convexity assumption in
Theorem~\ref{thm:existence} cannot, in general, be weakened to strict
geodesic convexity. Indeed,
Proposition~\ref{prop:strict-convexity-counterexample} shows that,
for every $n\ge3$ and $k=n$, there are a smooth strictly
geodesically convex domain containing the prescribed pole and
smooth data $f>0$ and $\varphi$ for which no
$C^2(\bar\Omega)$ admissible solution of Eq. \eqref{eq:main-problem} exists for any
$\lambda\in\RR$.

\item For $k=1$, neither convexity nor pole containment is needed for \eqref{eq:main-problem}. On every
bounded connected smooth domain, putting $w=u/V$ reduces the problem \eqref{eq:main-problem} to
\begin{eqnarray*}
\begin{cases}
\dive(V^2\nabla w)=Vf \qquad &\text{ in }\Omega,\\
V^2w_\nu= (\lambda+\varphi) V \qquad &\text{ on }\Sigma , 
\end{cases}
\end{eqnarray*}
where the compatibility condition gives $
 \lambda=\frac{\int_\Omega Vf-\int_\Sigma V\varphi}{\int_\Sigma V}.$ 
This is the classical weighted linear PDE with a Neumann boundary problem that has been used by Xia
\cite{Xia}.  
\item The solvability and uniqueness conclusions of
Theorem~\ref{thm:existence} also hold for an arbitrary pole
$p\in\HH^n$, with the same boundary regularity and data
assumptions, provided that
\begin{eqnarray*}
h>0,\qquad h-\frac{V_\nu}{V}g_\Sigma>0
\quad\text{on }\Sigma.
\end{eqnarray*}
The a priori estimates and the continuity argument use
these two inequalities separately, as explained in the
proof of Theorem~\ref{thm:existence}.
\end{enumerate} 
\end{rem}

For $1\le j\le n$, let $\lambda_j $ be the boundary constant for the
normalized $j$-Hessian equation $\sigma_j(A[u_j])=\binom nj$ in $\Omega$ with
$B_Vu_j=\lambda_j $ on $\Sigma$. The existence result of Theorem \ref{thm:existence} defines all these constants and the $j$-admissible solutions $u_j$ on
the same domain. We show that their comparison gives a natural geometric consequence and rigidity.  

\begin{thm}
\label{thm:hierarchy}
Let $\Omega, V$ satisfy the domain and pole hypotheses of
Theorem~\ref{thm:existence}. For $1\le k\le n$, let $\lambda_k $ be
the unique constant associated with Eq. 
\begin{eqnarray}\label{eq:main-problem-lambda-j}
\begin{cases}
 \sigma_k(\nabla^2u_k-u_kg)=\binom nk \qquad &\text{ in }\Omega,\\
(u_k)_\nu=\frac{V_\nu}{V}u_k+\lambda_k  \qquad &\text{ on }\Sigma. 
\end{cases}
\end{eqnarray}
Then
\begin{eqnarray}\label{eq:constant-hierarchy}
 0<\frac{n\int_\Omega V d\mu_g }{\int_\Sigma V d\sigma_g }
 =\lambda_1\le \lambda_2\le\cdots\le \lambda_n \leq \left(\frac {n\int_\Omega V d\mu_g} { \int_\Sigma VH_{n-1} d\sigma_g} \right)^{1/n} 
 \le\max_\Sigma \left( \frac{V_\nu}{V} \right) <1.
\end{eqnarray}
For any $1\le k<\ell\le n$, equality $\lambda_k=\lambda_\ell$ holds
if and only if $\Omega$ is a geodesic ball. 
\end{thm}
The formula for $\lambda_1$ in \eqref{eq:constant-hierarchy} is the standard 
weighted compatibility condition used in \cite{Xia}. The assertion
here compares the constants of different fully nonlinear Hessian equations with Neumann boundary conditions on the
same hyperbolic domain and identifies their common rigidity case.

 \ 

Let $(u,\lambda) $ be the normalized
solution supplied by Theorem~\ref{thm:existence} with $k=2$,
$f=\binom n2$, and $\varphi=0$. A further application of Theorem \ref{thm:existence} is that we use the normalized nonlinear Neumann potential to give the exact representations of the geometric deficits.
 Denote 
\begin{eqnarray}\label{eqn:defn-X}
  X\coloneqq \nabla u-u\nabla\log V.
\end{eqnarray}
\begin{thm}\label{thm:deficits}
Let $\Omega, V$ satisfy the domain and pole hypotheses of
Theorem~\ref{thm:existence}, 
and let $(u,\lambda) $ be the $2$-admissible 
solution to
\begin{eqnarray}\label{eqn:sigma-2-neumann}
\begin{cases}
\sigma_2(\nabla^2u-ug)=\binom{n}{2} \qquad &\text{ in }\Omega,\\
u_\nu=\dfrac{V_\nu}{V}u+\lambda   \qquad &\text{ on }\Sigma.
\end{cases}
\end{eqnarray} 
Then  
\begin{eqnarray}\label{eq:main-minkowski}
\begin{aligned}
&\lambda ^2\biggl[\left(\int_\Sigma V d\sigma_g \right)^2
 -n\left(\int_\Omega V d\mu_g \right)
       \left(\int_\Sigma VH d\sigma_g \right)\biggr]\\
&\quad=\left(\int_\Omega V(\tr A[u]-n) d\mu_g \right)
       \left(\int_\Omega V(\tr A[u]+n) d\mu_g \right)\\
&\qquad+\frac{n}{n-1}\left(\int_\Omega V d\mu_g \right)
       \int_\Sigma V\left(h-\frac{V_\nu}{V}g_\Sigma\right)
                  (X^T,X^T) d\sigma_g , 
\end{aligned}
\end{eqnarray}
and
\begin{eqnarray}\label{eq:main-hk}
\begin{aligned}
&\int_\Sigma\frac VH d\sigma_g -n\int_\Omega V d\mu_g \\
&=\int_\Sigma VH\left(\frac1H-
 \frac{\int_\Sigma V d\sigma_g }{\int_\Sigma VH d\sigma_g }\right)^2 d\sigma_g  +\frac{\left(\int_\Omega V(\tr A[u]-n) d\mu_g \right)
             \left(\int_\Omega V(\tr A[u]+n) d\mu_g \right)}
            {\lambda^2 \int_\Sigma VH d\sigma_g }\\
&\quad+\frac{n\int_\Omega V d\mu_g }
            {(n-1)\lambda^2 \int_\Sigma VH d\sigma_g }
       \int_\Sigma V\left(h-\frac{V_\nu}{V}g_\Sigma\right)
                  (X^T,X^T) d\sigma_g .
\end{aligned}
\end{eqnarray}
\end{thm}

The contribution of \eqref{eq:main-minkowski} and
\eqref{eq:main-hk} is their representation through the
normalized nonlinear Neumann potential
\eqref{eqn:sigma-2-neumann}.
A different nonnegative linear remainder already follows from the weighted
Reilly argument in \cite{LiXia}; see \eqref{eq:linear-deficit} for comparison.
Explicit Heintze--Karcher deficit identities in the substatic
setting were established by Fogagnolo--Pinamonti
\cite[Theorem~3.6]{FogagnoloPinamonti2022}, using a linear
Dirichlet problem. For the fully nonlinear
Neumann potential used here, the
normalization of $\sigma_2(A[u])$ gives
\begin{eqnarray*}
\int_\Omega V(\tr A[u]-n) d\mu_g
=\int_\Omega
 \frac{V|A[u]-g|^2}{\tr A[u]+n-2} d\mu_g.
\end{eqnarray*}

 Theorems \ref{thm:hierarchy} and \ref{thm:deficits}  give
an alternative proof of the weighted Minkowski inequality
\cite[Theorem~1.1]{Xia} and the Heintze--Karcher inequality
\cite[Eq.~(1.6)]{LiXia}. For results in more general ambient settings, see
\cite[Section~3]{Brendle},
\cite[Theorem~1.3]{QXReilly}, and
\cite[Sections~4 and~5]{LiXia}.
\begin{cor}
    Let $\Omega, V$ satisfy the domain and pole hypotheses of Theorem~\ref{thm:existence}, then 
\begin{eqnarray*}
\left(\int_\Sigma V d\sigma_g \right)^2
\ge n\left(\int_\Omega V d\mu_g \right)
\left(\int_\Sigma VH d\sigma_g \right), 
\end{eqnarray*} 
and 
\begin{eqnarray*}
\int_\Sigma\frac VH d\sigma_g \ge n\int_\Omega V d\mu_g . 
\end{eqnarray*}
Equality in either inequality holds if and only if $\Omega$ is a geodesic ball. 
\end{cor}

We finally give an analytic proof of the following weighted
Alexandrov--Fenchel inequalities using 
Theorem~\ref{thm:existence}.
\begin{thm} 
\label{thm:af}
 Let $\Omega, V$ satisfy the domain and pole hypotheses of Theorem~\ref{thm:existence}, 
then 
\begin{eqnarray}\label{eq:af-two-hessian-bounds}
 \frac{n\displaystyle\int_\Omega V d\mu_g}
      {\displaystyle\int_\Sigma V d\sigma_g}
 \le
 \left(
 \frac{n\displaystyle\int_\Omega V d\mu_g}
      {\displaystyle\int_\Sigma VH_j d\sigma_g}
 \right)^{1/(j+1)},\qquad 1\leq j\leq n-1.
\end{eqnarray}
Moreover,
\begin{eqnarray}\label{eq:af-adjacent-analytic}
 \left(\int_\Sigma VH_j d\sigma_g\right)^2
 \ge
 \left(\int_\Sigma VH_{j-1} d\sigma_g\right)
 \left(\int_\Sigma VH_{j+1} d\sigma_g\right),
 \qquad 1\le j\le n-2,
\end{eqnarray}
and
\begin{eqnarray}\label{eq:global-weighted-AF}
 \left(\int_\Sigma V d\sigma_g\right)^{j+1}
 \ge
 \left(n\int_\Omega V d\mu_g\right)^j
 \int_\Sigma VH_j d\sigma_g,
 \qquad 1\le j\le n-1.
\end{eqnarray}
Equivalently,
\begin{eqnarray}\label{eq:global-weighted-AF-unnormalized}
     \left(\int_\Sigma V d\sigma_g\right)^{j+1}
 \ge
 \frac{n^j}{\binom{n-1}{j}}
 \left(\int_\Omega V d\mu_g\right)^j
 \int_\Sigma V\sigma_j(h) d\sigma_g.
\end{eqnarray}
Equality in \eqref{eq:af-adjacent-analytic} or
\eqref{eq:global-weighted-AF} holds if and only if $\Omega$ is a geodesic ball.
\end{thm}

Here $H_j$ is defined by \eqref{eqn:H-j-defn}, and $\lambda$ denotes the boundary constant
in \eqref{eqn:sigma-2-neumann}. The case $j=1$ of \eqref{eq:global-weighted-AF} is
Xia's weighted Minkowski inequality \cite[Theorem~1.1]{Xia}.
For $n\ge3$, the case $j=1$ of
\eqref{eq:af-adjacent-analytic} is the weighted
Alexandrov--Fenchel inequality of Li--Xia
\cite[Theorem~1.9]{LiXia}; together, these two inequalities
give the case $j=2$ of \eqref{eq:global-weighted-AF}. 
Here \eqref{eq:af-adjacent-analytic} compares the weighted
curvature integrals of three consecutive orders, and
\eqref{eq:af-adjacent-deficit} gives an exact remainder. The results of Hu--Li \cite{HuLi} and Hu--Li--Wei \cite{HuLiWei}
give lower bounds for weighted curvature integrals
in terms of weighted volume or ordinary quermassintegrals.
Our comparison \eqref{eq:global-weighted-AF} gives
complementary upper bounds in terms of weighted volume
and weighted boundary area.

\begin{rem}
After completing this manuscript, we became aware of the work of
Xu \cite[Theorem~1 and Remark~2(ii)]{Xu2026WeightedAF},
who established \eqref{eq:af-adjacent-analytic} and
\eqref{eq:global-weighted-AF} through Euclidean mixed-volume
theory under non-strict static convexity and $h\ge0$, as well as in more general settings.
Here we provide an analytic proof of these inequalities based on
the nonlinear Neumann problem  \eqref{eqn:sigma-2-neumann} and an auxiliary linear equation \eqref{eq:af-boundary-equation} on $\Sigma$.
\end{rem}

We conclude this section with an outline of the proofs.
The main difficulties in Theorem~\ref{thm:existence}
come from the kernel generated by $V$ and the negative
ambient curvature. We first bound the boundary constant
by comparison principle with constant functions. Inspired by the idea used in \cite{Jin-Li-Liyanyan2007,MX-2016,MaQiuXu2016} etc., a boundary
correction gives homogeneous Neumann data and allows
a gradient estimate which is independent of a positive
lower bound for the rescaled right hand side.
Applying this estimate to the normalized oscillation
of $u/V$, and combining the weighted divergence identity
with an energy estimate and the Poincar\'e inequality,
gives the oscillation bound. Under the normalization
$\int_\Omega u/V d\mu_g=0$, this yields the required
$C^1$ estimate.

For the second derivative estimate, we follow the
reduction to a boundary double normal estimate used
in \cite{LTU,MQ}, with modifications adapted to the
static boundary condition. Tangential differentiation
of that condition determines a bounded correction to
$A[u]$ which cancels its mixed normal--tangential
components on $\Sigma$. After subtracting a suitable
multiple of $V$, we estimate the largest eigenvalue
of the corrected tensor by a logarithmic test function
containing $\gamma\log(-u)$, with $\gamma>1$.
This term compensates for the unfavorable trace term
in the curvature commutator. Since the corrected tensor
need not commute with $A[u]$, the interior estimate
uses a spectral comparison together with the positive
terms in the second derivative formula for its largest
eigenvalue. At a tangential boundary maximum, strict
static convexity supplies the required sign.
Upper and lower barriers for the boundary residual,
constructed in a fixed thin neighborhood of $\Sigma$ as in
\cite{MQ}, then control the double normal derivative
and close the global $C^2$ estimate.

We apply the continuity method with $\lambda$ included
among the unknowns. Division by $V$ removes the
zeroth-order term from the linearized equation, and
the normalization removes its constant kernel.
The maximum principle and the Fredholm alternative
give invertibility of the augmented linearized operator,
while the a priori estimates and oblique regularity
give compactness along the continuity path.
The same maximum principle argument proves uniqueness
of $\lambda$ and uniqueness of $u$ modulo multiples
of $V$.  

The geometric applications use these solutions with
constant right hand sides and constant boundary data.
For Theorem~\ref{thm:hierarchy}, the Newton--Maclaurin
inequalities and the comparison principle order the
boundary constants of the normalized Hessian equations.
The weighted compatibility identity determines the
linear endpoint. Equality between two distinct orders
forces equality in the Newton--Maclaurin inequality
and hence $A[u]=g$.
For Theorem~\ref{thm:deficits}, we apply the weighted Reilly identities to the normalized
$2$-Hessian solution and retain the remainder terms.
The equation expresses the trace excess through
$|A[u]-g|^2$, giving the exact Minkowski deficit.
Combining this identity with the remainder in the
weighted Cauchy--Schwarz inequality gives the
Heintze--Karcher deficit.  

To prove Theorem~\ref{thm:af}, we first obtain lower
and upper estimates for the boundary constant of
the same normalized $2$-Hessian solution from the Reilly identities. For each higher order,
we solve an auxiliary linear equation on $\Sigma$
whose right hand side is chosen so that a corrected
tangential Hessian has zero contraction with the
corresponding Newton tensor. The Codazzi identity
and integration by parts give the mixed boundary
formula, and the concavity of the elementary
symmetric polynomial gives the remaining nonnegative
term. This proves the adjacent-order inequalities
and their exact deficit identity.
The resulting monotonicity of consecutive integral
ratios transfers the initial estimate for the
boundary constant to every curvature order.
Eliminating that constant proves the global
weighted comparisons.
Equality in these comparisons reduces to $A[u]=g$,
while equality in an adjacent-order comparison
forces $\Sigma$ to be totally umbilical.
These conditions give the geodesic ball rigidity.

\bigskip
\noindent\textit{The paper is organized as follows.}
Section~\ref{sec:preliminaries} records the static
identities, the algebraic inequalities, and the
weighted Reilly formula. In Section~\ref{sec:existence}, we establish the a priori
estimates for Eq. \eqref{eq:main-problem}, prove Theorem~\ref{thm:existence}, and give
a counterexample to solvability under strict geodesic convexity alone.
Section~\ref{sec:geometry} proves the comparison of
boundary constants, the geometric deficit identities,
and the weighted Alexandrov--Fenchel inequalities.

\section{Preliminaries}\label{sec:preliminaries}

\subsection{Elementary symmetric polynomial}

Let $A$ be a real symmetric $n\times n$ matrix, and let
$\lambda(A)=(\lambda_1(A),\ldots,\lambda_n(A))$ be its eigenvalues.
For $1\le k\le n$, define
\begin{eqnarray*}
\sigma_k(A)\coloneqq
\sum_{1\le i_1<i_2<\cdots<i_k\le n}
\lambda_{i_1}(A)\lambda_{i_2}(A)\cdots\lambda_{i_k}(A).
\end{eqnarray*}
We use $\sigma_0=1$ and $\sigma_j=0$ when $j$ exceeds
the dimension of the matrix under consideration. 
The notation $\sigma_k(\lambda|i)$ means that the
$i$-th component of $\lambda$ is replaced by zero.
The notation $H_j$ is reserved for the normalized
boundary curvatures in \eqref{eqn:H-j-defn} below.

Let
\begin{eqnarray*}
\Gamma_k= \left\{\lambda\in \mathbb{R}^n:\sigma_j(\lambda)>0\text{ for }1\le j\le k \right\}.
\end{eqnarray*}

\begin{defn}\label{def-k-admissible}
Let $1\le k\le n$ and $u\in C^2(\Omega)$.
We say that $u$ is $k$-admissible in $\Omega$ if
\begin{eqnarray*}
 \lambda(A[u])\in\Gamma_k\qquad\text{in }\Omega.
\end{eqnarray*}
For $u\in C^2(\bar\Omega)$, admissibility on
$\bar\Omega$ means that the same condition holds at
every point of $\bar\Omega$.
\end{defn}

 It is well-known that 
\begin{eqnarray*}
F(A)\coloneqq \sigma_k(A)^{1/k}
\end{eqnarray*}
is smooth, concave, homogeneous of degree one, and elliptic on the cone $\Gamma_k$.  
Denote  
\begin{eqnarray*}
F^{ij}\coloneqq \frac{\partial F}{\partial A_{ij}},\qquad
\mathcal{F} \coloneqq \tr(F^{ij}),  \qquad \mathscr L\coloneqq F^{ij}\nabla_{ij}. 
\end{eqnarray*}
In a frame diagonalizing $A$, denote $$F_i\coloneqq F^{ii},\qquad \lambda_i\coloneqq A_{ii}.$$  

We collect the properties of $F_i$ that will be used later.

\begin{lem}\label{lem:spectral}
Let $1\le k\le n$ and let $A$ be a symmetric matrix with $\lambda(A)\in\Gamma_k$. Then the following hold. 
\begin{enumerate}[label=(\arabic*), leftmargin=*, labelindent=0pt]
\item $\sum_{i=1}^nF_i\lambda_i=F(A).$

\item $\mathcal F\ge F(g)=\binom nk^{1/k}.$
\item $\frac{\sigma_1(A)}n
\ge\left(\frac{\sigma_k(A)}{\binom nk}\right)^{1/k}.$

\item If $\lambda_i<0$, then
  $F_i\ge\mathcal F/(n-k+1)$.
\item Given sufficiently small $\epsilon,\delta>0$, there exists $c_*>0$, depending only on $n,k,\epsilon,\delta$, such that
\begin{eqnarray*}
\lambda_i>0,\quad \lambda_i\ge\delta\lambda_{\max},\quad
-\lambda_{\min}\ge\epsilon\lambda_i
\quad\Longrightarrow\quad F_i\ge c_*\mathcal{F} .
\end{eqnarray*}
\end{enumerate}
\end{lem}

The assertions (1)-(3) in Lemma \ref{lem:spectral}
follow from homogeneity and the Newton--Maclaurin inequalities;
see \cite{CNS,Lin-Trudinger1994}.
For $2\le k<n$, the assertions (4)-(5) of Lemma~\ref{lem:spectral}
follow from \cite[Lemmas~2.5 and~2.7]{CZ}, respectively,
by taking the denominator to be $\sigma_0=1$ and canceling
the common positive factor in the derivatives of
$\sigma_k^{1/k}$. The case $k=1$ is immediate, while
the hypotheses involving a negative eigenvalue cannot
occur when $k=n$.

\subsection{The static potential}

Throughout the paper, $\HH^n$ has constant sectional curvature $-1$, and $\nabla$, $\Delta$, and $\langle \cdot, \cdot\rangle$ refer to its metric $g$. The outward unit normal to $\Sigma=\partial\Omega$ is $\nu$. 
 Fix any $p\in\HH^n$ and set $V(x)\coloneqq \cosh d(p,x)$. The function $V$ is smooth also at $p$, and
\begin{eqnarray}\label{eq:static}
\nabla^2V=Vg,\qquad \Delta V=nV,\qquad |\nabla V|^2=V^2-1.
\end{eqnarray}
When the domain $\Omega$ is strictly geodesically convex and $p\in\Omega$, one has $\frac{V_\nu}{V}>0$ on $\Sigma$: the radial geodesic from $p$ meets the supporting tangent hyperplane with strictly positive outward component.

We define the second fundamental form and the normalized
mean curvature by
\begin{eqnarray}\label{eqn:H-defn}
h(Y,Z)\coloneqq \langle \nabla_Y\nu, Z\rangle,\qquad
H\coloneqq \frac{\tr_{g_\Sigma}h}{n-1},
\end{eqnarray}
and the normalized $j$-th mean curvature \begin{eqnarray}\label{eqn:H-j-defn}
    H_j\coloneqq \binom{n-1}{j}^{-1}\sigma_j(h),
\end{eqnarray} with $H_0=1$
and $H_1=H$. 
Thus the principal curvatures of a geodesic sphere of radius $R$ are $\coth R$. 
A bounded smooth geodesically convex domain has a connected boundary. The volume and boundary measures are denoted by $d\mu_g $ and $d\sigma_g $, respectively; they are occasionally suppressed in integrals. All tensor inequalities are quadratic-form inequalities, and eigenvalues of a covariant tensor are taken using the indicated ambient metric.

The following lemma records direct consequences of the
spatial-projection formulas and convexity correspondence of
Hu--Li--Wan--Xu
\cite[Lemmas~2.1--2.2 and Corollary~2.1]{HuLiWanXu}. We include a proof for the reader's convenience.
\begin{lem}\label{lem:static-convex-domain}
Let $\Omega\Subset\HH^n$ be connected with $C^2$ boundary, fix
$p\in\HH^n$, and set $V=\cosh d(p,\cdot)$. If
$h-\frac{V_\nu}{V}g_\Sigma\ge0$, the spatial projection of
$\bar\Omega$ is a Euclidean convex body; it is strictly convex when
the tensor inequality is strict. If $p\in\Omega$, then $\Omega$ is
starshaped about $p$ and strictly geodesically convex,
$\Sigma$ is connected, and
\begin{eqnarray*}
 \min_\Sigma\frac{V_\nu}{V}\ge\tanh\dist(p,\Sigma)>0,
 \qquad h\ge\tanh\dist(p,\Sigma)g_\Sigma>0.
\end{eqnarray*}
In particular, the second fundamental form has a uniform positive
lower bound on the fixed domain, even when static convexity is
non-strict.
\end{lem}
 
\begin{proof}
Place $p$ at $(1,0)$ in the hyperboloid model and write
$\iota(x)=(\sqrt{1+|x|^2},x)$ and $K=\iota^{-1}(\bar\Omega)$.
Let $\eta$ be the outward Euclidean unit normal, set $s=x\cdot\eta$
and $\ell=\sqrt{1+s^2}$, and denote the Euclidean second fundamental
form by $h_E$. The projection identities in
\cite[Lemmas~2.1 and~2.2]{HuLiWanXu} give
\begin{eqnarray*}
 \nu=\ell^{-1}(sV,\eta+sx),\qquad
 \frac{V_\nu}{V}=\frac{s}{\ell},
\end{eqnarray*}
and
\begin{eqnarray*}
 h_E=\ell\left(h-\frac{V_\nu}{V}g_\Sigma\right).
\end{eqnarray*}
By \cite[Corollary~2.1]{HuLiWanXu}, $K$ is a Euclidean
convex body, with strict convexity when the tensor inequality
is strict.

If $p\in\Omega$, then $0\in\operatorname{int}K$ and
$\bar B_{\dist(p,\Sigma)}(p)\subset\bar\Omega$. Its spatial projection
is the Euclidean ball of radius $\sinh\dist(p,\Sigma)$ about the
origin. Every supporting hyperplane of $K$ therefore satisfies
$s\ge\sinh\dist(p,\Sigma)$, and hence
\begin{eqnarray*}
 \frac{V_\nu}{V}=\frac{s}{\sqrt{1+s^2}}
 \ge\tanh\dist(p,\Sigma)>0.
\end{eqnarray*}
The assumed static inequality now gives
$h\ge\tanh\dist(p,\Sigma)g_\Sigma>0$.
Convexity of $K$ and containment of the origin make every ray meet
$\partial K$ once. Since the projection sends the radial hyperbolic
geodesic to $x=\sinh r \vartheta$, the same statement holds in
$\Omega$, proving star-shapedness. Connectedness is preserved by the
projection.  
\end{proof}
 
Recall the notations in \eqref{eq:static-pair-intro} and \eqref{eqn:defn-X}. Then $A[u]$ is a Codazzi tensor and satisfies 
\begin{eqnarray}\label{eq:X-derivative}
\nabla_{ij}\log V=g_{ij}-(\log V)_i (\log V)_j,\qquad
\nabla_iX_j=A_{ij}-X_i (\log V)_j.
\end{eqnarray}
Indeed, the first identity in \eqref{eq:X-derivative} follows from \eqref{eq:static}, and differentiating $X_j=u_j-u (\log V)_j$ gives the second.
\begin{lem}
If $F(A[u])=f^{1/k}$, then
\begin{flalign}
& \text{(1)}\quad
  \mathscr{L}u = f^{1/k} + u\mathcal{F}.
  && \label{eqn:Lu-eq}\\
& \text{(2)}\quad
  \mathscr{L}(u_\ell)
  = \nabla_\ell(f^{1/k}) + F^{\ell i}u_i.
  && \label{eq:gradient-commutation}\\
& \text{(3)}\quad
  \mathscr{L}A_{\xi\xi}
  \ge (f^{1/k})_{;\xi\xi}
  - \mathcal{F}A_{\xi\xi} + f^{1/k},
  && \label{eq:tensor-commutation}
\end{flalign}
for any unit direction $\xi$, interpreted by radial parallel extension at the point of evaluation.
\end{lem}
\begin{proof}
The assertion \eqref{eqn:Lu-eq} follows from Lemma \ref{lem:spectral} (1). To show \eqref{eq:gradient-commutation}, by 
differentiating $F(A)=f^{1/k}$ and commuting one derivative of $u_\ell$ gives
\begin{eqnarray*}
F^{ij}\nabla_{ij}u_\ell
=F^{ij}\nabla_\ell u_{ij}-\mathcal{F}  u_\ell+F^{\ell i}u_i
=\nabla_\ell(f^{1/k})+F^{\ell i}u_i, 
\end{eqnarray*} it yields \eqref{eq:gradient-commutation}.
Using
the second-derivative commutation formula yields
\begin{eqnarray*}
F^{ij}\nabla_{ij}A_{pq}
=F^{ij}\nabla_{pq}A_{ij}-\mathcal{F}  A_{pq}+f^{1/k}g_{pq}.
\end{eqnarray*}
In an orthonormal eigenframe of $A$, the diagonal
commutation identity is
\begin{eqnarray*}
\nabla_{ii}A_{pp}
=\nabla_{pp}A_{ii}-A_{pp}+A_{ii}.
\end{eqnarray*} Twice differentiating the equation and using concavity proves \eqref{eq:tensor-commutation}. Radial parallel transport of $\xi$ has $\nabla\xi=0$ and $\operatorname{Sym}\nabla^2\xi=0$ at the base point, so no additional terms arise from differentiating the direction. This completes the proof.
\end{proof}

\subsection{Weighted Reilly formula}
We recall the weighted  Reilly formula using the boundary
quantity in the Neumann boundary condition.

\begin{lem}\label{lem:reilly}
For $u\in C^2(\bar\Omega)$,  
there hold 
\begin{eqnarray}\label{eq:flux}
\int_\Omega V\tr A d\mu_g =\int_\Sigma V (B_V u) d\sigma_g ,
\end{eqnarray}
and
\begin{eqnarray}
 2\int_\Omega V\sigma_2(A) d\mu_g 
&=&\int_\Sigma V\bigg[(n-1)  (B_V u)^2 H
+\left(h-\frac{V_\nu}{V}g_\Sigma \right)(X^T,X^T) 
 \notag \\&&  -2 \left\langle X^T, \nabla_\Sigma(B_V u) \right\rangle \bigg] d\sigma_g.  \label{eq:reilly}
\end{eqnarray}
\end{lem}

These are specializations of the generalized Reilly formulas of
Qiu--Xia \cite{QXReilly} and Li--Xia \cite{LiXia}, also used in Xia's derivation of the
weighted Minkowski inequality \cite{Xia}. For the reader's convenience, we include a  proof here.  
\begin{proof}
Taking the trace in \eqref{eq:X-derivative} gives $\dive(VX)=V\tr A$, and integration by parts proves \eqref{eq:flux}. 

Denote $N\coloneqq (\tr A)g-A$ for the first Newton tensor. The Codazzi property gives $\nabla_iN^{ij}=0$. Symmetry and \eqref{eq:X-derivative} then imply
\begin{eqnarray*}
\dive(VNX)=VN^{ij}A_{ij}=2V\sigma_2(A).
\end{eqnarray*}
Decompose $X=(B_V u) \nu+X^T$ at the boundary. Differentiating $B_V u=\langle X, \nu\rangle$ and using \eqref{eq:X-derivative}, we obtain
\begin{eqnarray}\label{eqn:mixed-bdry-A-nu-alpha}
A_{\nu\alpha}=(B_V u)_\alpha-\left(h_{\alpha\gamma}-\frac{V_\nu}{V} g_{\alpha\gamma} \right)X^T_\gamma.
\end{eqnarray}
The tangential trace of the same identity gives
\begin{eqnarray*}
\tr_{T\Sigma}A
=V^{-1}\dive_\Sigma(VX^T)+(n-1)H(B_V u).
\end{eqnarray*}
Therefore
\begin{eqnarray*}
VN(\nu,X) &=& (B_V u)\dive_\Sigma(VX^T) + V(n-1)H(B_V u)^2 \\& & +V\left(h-\frac{V_\nu}{V}g_\Sigma \right)(X^T,X^T)-V\langle X^T, \nabla_\Sigma(B_V u)\rangle.
\end{eqnarray*}
Integrating the tangential divergence term by parts on the closed hypersurface $\Sigma$ proves \eqref{eq:reilly}.
\end{proof}

\begin{rem}
In the convention \eqref{eqn:H-defn}, identity~\eqref{eq:reilly}
is equivalently
\begin{eqnarray*}
\begin{aligned}
&\int_\Omega V\bigl[(\Delta u-nu)^2-|\nabla^2u-ug|^2\bigr] d\mu_g = \int_\Sigma V_\nu
 \bigl[|\nabla_\Sigma u|^2+(n-1)u^2\bigr] d\sigma_g   \\
& \qquad + \int_\Sigma V\bigl[2u_\nu\Delta_\Sigma u
 +(n-1)Hu_\nu^2+h(\nabla_\Sigma u,\nabla_\Sigma u)
 -2(n-1)uu_\nu\bigr] d\sigma_g .
\end{aligned}
\end{eqnarray*}

More generally, if $B_Vu=\lambda+\varphi $ on $\Sigma$ and $h-\frac{V_\nu}{V}g_\Sigma>0$, completing
the square in~\eqref{eq:reilly} gives the exact identity
\begin{eqnarray*}
\begin{aligned}
2\int_\Omega V\sigma_2(A[u]) d\mu_g 
={}&(n-1)\int_\Sigma VH(\lambda+\varphi )^2 d\sigma_g \\
&+\int_\Sigma V\left|
 \left(h-\frac{V_\nu}{V}g_\Sigma\right)^{1/2}X^T
 -\left(h-\frac{V_\nu}{V}g_\Sigma\right)^{-1/2}\nabla_\Sigma\varphi
 \right|^2 d\sigma_g \\
&-\int_\Sigma V \left\langle \left(h-\frac{V_\nu}{V}g_\Sigma\right)^{-1}\nabla_\Sigma\varphi, \nabla_\Sigma\varphi \right\rangle d\sigma_g .
\end{aligned}
\end{eqnarray*}
Here inverses and positive square roots are taken after identifying
$h-\frac{V_\nu}{V}g_\Sigma$ with a positive self-adjoint endomorphism
using the induced metric. 
For constant $\varphi$, the lower bound remains valid under
the weaker hypothesis $\left(h-\frac{V_\nu}{V}g_\Sigma\right)\ge0$ directly from~\eqref{eq:reilly},
without taking an inverse.
\end{rem}

For every $k$-admissible solution of $\sigma_k(A[u])=f>0$ in $\Omega$ with
$B_Vu=\lambda+\varphi$ on $\Sigma$, the Newton--Maclaurin inequality gives
\begin{eqnarray*}
\lambda \int_\Sigma V d\sigma_g 
 \ge n\int_\Omega V\left(\frac{f}{\binom nk}\right)^{1/k} d\mu_g 
       -\int_\Sigma V\varphi d\sigma_g .
\end{eqnarray*}
For $k=1$ this is an equality and is the linear compatibility condition, used in \cite{Xia}.
For $k>1$ it is only a necessary bound; it is not a sufficiency criterion
for nonlinear solvability.
The estimates below essentially address the nonlinear range $2\le k\le n$.

\section{A priori estimates}\label{sec:existence}

Throughout this section, $C$ denotes a positive constant in the a priori estimates that may change from line to line. Unless otherwise stated, it depends only on the fixed geometry, $V$, $n$, $k$, a positive lower bound for $f$, $\|f\|_{C^2}$, and $\|\varphi\|_{C^3}$. Any additional dependence will be indicated explicitly. 

 We continue to use $$A\coloneqq A[u],  \qquad F(A)\coloneqq \sigma_k(A)^{1/k}, \qquad \mathcal{F} \coloneqq \tr(F^{ij}),$$ and the linearized operator $$\mathscr L\coloneqq F^{ij}\nabla_{ij}.$$
We first derive the estimates for smooth domains, data, and solutions.
Their constants depend only on the indicated finite-regularity norms;
the approximation argument in the proof of Theorem~\ref{thm:existence}
then gives the stated regularity. Under the principal hypotheses of
that theorem, Lemma~\ref{lem:static-convex-domain} supplies $h>0$,
while strict static convexity supplies the remaining boundary gap. 

\subsection{The gradient estimate}

We will prove the gradient estimate for the general boundary data in Eq. 
\eqref{eq:main-problem}. Throughout this subsection, the domain is
strictly geodesically convex; this follows from
Lemma~\ref{lem:static-convex-domain} in the main static-convex setting.
The proof of this subsection also permits arbitrary poles when
$h>0$ is assumed separately. 

We first give the bounds for $\lambda$ in Eq. \eqref{eq:main-problem}, which do not require convexity of the domain $\Omega$ or a sign for $V_\nu$.
\begin{lem} \label{lem:c-bound}
Let $(u,\lambda) $ be a $k$-admissible solution of Eq. 
\eqref{eq:main-problem}. Denote 
\begin{eqnarray*}
 C_1\coloneqq \left(\frac{\min_{\bar\Omega}f}{\binom nk}\right)^{1/k},
 \qquad
 C_2\coloneqq \left(\frac{\max_{\bar\Omega}f}{\binom nk}\right)^{1/k}.
\end{eqnarray*}
Then
\begin{eqnarray}\label{eq:c-bound}
 \min_\Sigma\left(C_1\frac{V_\nu}{V}-\varphi\right)
 \le \lambda \le
 \max_\Sigma\left(C_2\frac{V_\nu}{V}-\varphi\right).
\end{eqnarray}
In particular, $|\lambda|$ is bounded in terms of the fixed data. 
\end{lem}
\begin{proof}
The constant function $-C$ satisfies
$A[-C]=Cg$ and $B_V(-C)=C\frac{V_\nu}{V}$.
First choose any $0<C<C_1$. If $(u+C)/V$ attained its maximum in
$\mathring\Omega$, at that interior point would give
$$A[u]=\n^2 u- ug = V \n^2 \left( \frac{u+C}{V}\right) + \n \left( \frac{u+C}{V}\right) \otimes \n V+\n V\otimes \n \left( \frac{u+C}{V}\right)+Cg  \le Cg, $$ it follows   
$$f=\sigma_k(A[u]) \le\binom nk C^k,$$ a contradiction to $C<C_1$. At a boundary maximum point, 
the outward normal derivative of $(u+C)/V$ is nonnegative, whence
\begin{eqnarray*}
 0\le V\partial_\nu\left(\frac{u+C}{V}\right)
   =\lambda +\varphi-C\frac{V_\nu}{V}.
\end{eqnarray*}
This proves the lower bound of $\lambda$ with $C$ in place of $C_1$, and passing
to the limit $C\uparrow C_1$ proves the stated lower bound in \eqref{eq:c-bound}.

For $C>C_2$ an interior minimum of $(u+C)/V$ would give
$A[u]\ge Cg$ and the opposite contradiction. Its outward derivative
at a boundary minimum is nonpositive. Letting $C\downarrow C_2$
proves the upper bound in \eqref{eq:c-bound}. 
\end{proof}

We next establish a gradient estimate that is uniform under rescaling.
The modified boundary function in \eqref{eqn:defn-w-C1} and the logarithmic
test function in \eqref{eqn:C1-defn-P} are motivated by the
Neumann gradient estimates of Ma--Xu \cite[Section~3]{MX-2016} for the prescribed mean
curvature equation and of Ma--Qiu--Xu \cite{MaQiuXu2016}
for Hessian equations in Euclidean space; the latter argument is
also presented in \cite[Theorem~10 and Section~3.1]{MQ}.

\begin{lem}\label{lem:gradient}
Let $U$ be a $k$-admissible solution of the problem: 
\begin{eqnarray}\label{eq:main-problem-F(U)}
\begin{cases}
 F(A[U])=s(x)>0  \qquad &\text{ in }\Omega,\\
U_\nu=\frac{V_\nu}{V}U+b(x) \qquad &\text{ on }\Sigma,
\end{cases}
\end{eqnarray}
  then
  \begin{eqnarray}\label{eqn:C1-bound-U}
\|\nabla U\|_{C^0(\bar\Omega)}\le C,  \end{eqnarray} where  the constant $C>0$ depends only on $\Omega, V$,  $\|U\|_{C^0(\bar\Omega)}$, $\|s\|_{C^1(\bar\Omega)}$, and 
$\|b\|_{C^3(\Sigma)}$, but not on $\inf_{\bar\Omega}s$.
\end{lem}
\begin{proof}
Extend $\frac{V_\nu}{V}$ and $b$ from $\Sigma$ to $\bar\Omega$ with
controlled $C^3$ norms. These extensions are used only in this proof.
Choose a $C^3$ function $d(x) $ which agrees with the interior distance
to $\Sigma$ in a small neighborhood of $\Sigma$, vanishes outside a slightly larger neighborhood,
and satisfying
\begin{eqnarray*}
 \rho(x)\coloneqq 1+\frac{V_\nu}{V}d>\frac12 \qquad \text{ in } \bar \Omega. 
\end{eqnarray*}
Thus \begin{eqnarray}\label{eqn:d-function}
    d=0, \qquad d_\nu=-1\qquad \text{ on } \Sigma.
\end{eqnarray} Define the function  \begin{eqnarray}\label{eqn:defn-w-C1}
    w\coloneqq \left(1+\frac{V_\nu}{V}d\right) U+b(x)d(x).
\end{eqnarray} Direct differentiation on $\Sigma$ and using \eqref{eqn:d-function} imply
 \begin{eqnarray}\label{eqn:w-nu=0}
    w_\nu=0 \qquad \text{ on } \Sigma. 
\end{eqnarray} 
Denote the tensor
\begin{eqnarray*}
 K:=\rho A[U]
 =\nabla^2w+\beta\otimes\mathrm dw
       +\mathrm dw\otimes\beta+C_2w+D_2,
\end{eqnarray*}
where the fixed coefficient tensors are
\begin{align*}
 \beta&\coloneqq-\mathrm d\log\rho,\\
 C_2&\coloneqq \mathrm d\log\rho\otimes\mathrm d\log\rho
                  -\nabla^2\log\rho-g,\\
 D_2&\coloneqq -(b(x)d(x)) C_2
       +\mathrm d\log\rho\otimes\mathrm d(b(x)d(x))
       +\mathrm d(b(x)d(x))\otimes\mathrm d\log\rho-\nabla^2(b(x)d(x)).
\end{align*}
It follows from \eqref{eq:main-problem-F(U)} that $K\in\Gamma_k$ and  
\begin{eqnarray}\label{eq:F(K)}
    F(K)=\rho(x) s(x), 
\end{eqnarray} 

In the remainder of this proof, $F^{ij}=\partial F/\partial K_{ij}$
is evaluated at $K$, $\mathcal F=\tr(F^{ij})$, and
\begin{eqnarray*}
 \mathcal L\coloneqq F^{ij}\nabla_{ij}+2F^{ij}\beta_i\nabla_j. 
\end{eqnarray*} 
The third-order commutation formula on $\mathbb H^n$ gives
\begin{eqnarray}\label{eqn:commutation-C1}
 \nabla_iw_{j\ell}-\nabla_\ell w_{ij}
 =g_{j\ell}w_i-g_{ij}w_\ell.
\end{eqnarray}
Differentiating \eqref{eq:F(K)} in $e_\ell$ direction gives 
\begin{eqnarray*}
 \nabla_\ell(\rho s)
 &=&F^{ij}\nabla_\ell K_{ij}
 =F^{ij}\nabla_\ell w_{ij}
      +2F^{ij}\beta_iw_{j\ell}
      +2F^{ij}(\nabla_\ell\beta_i)w_j\\
 &&{}+F^{ij}(C_2)_{ij}w_\ell
      +F^{ij}\bigl(w\nabla_\ell(C_2)_{ij}
                        +\nabla_\ell(D_2)_{ij}\bigr).
\end{eqnarray*}
Together with \eqref{eqn:commutation-C1}, we compute 
\begin{eqnarray}
 \mathcal L\left(\frac12|\nabla w|^2\right)
 &=&F^{ij}w_{\ell i}w_{\ell j}
      +F^{ij}w_\ell\nabla_\ell w_{ij}
      +F^{ij}w_iw_j-\mathcal F|\nabla w|^2
+2F^{ij}\beta_iw_\ell w_{j\ell}
\notag \\&=& F^{ij}w_{\ell i}w_{\ell j}
      +w_\ell\nabla_\ell(\rho s)
      +F^{ij}w_iw_j-\mathcal F|\nabla w|^2 \notag\\
 &&-2F^{ij}(\nabla_\ell\beta_i)w_jw_\ell
      -F^{ij}(C_2)_{ij}|\nabla w|^2 \notag\\
 &&-wF^{ij}\nabla_\ell(C_2)_{ij} w_\ell
      -F^{ij}\nabla_\ell(D_2)_{ij} w_\ell
      \notag \\ 
 &\ge& F^{ij}w_{\ell i}w_{\ell j}
      -C\bigl(|\nabla w|^2+|\nabla w|\bigr)\mathcal F
      -|\nabla w| |\nabla(\rho s)|.
 \label{eq:gradient-drift-estimate}
\end{eqnarray}
By the degree-one homogeneity of $F$ in \eqref{eq:F(K)}, we have
\begin{eqnarray}
 \mathcal Lw
 &=&F^{ij}w_{ij}+2F^{ij}\beta_iw_j 
 =F^{ij}K_{ij}
     -F^{ij}\left((C_2)_{ij}w+(D_2)_{ij}\right) \notag
\\
 &\ge& \rho(x) s(x)-\bigl(|C_2| |w|+|D_2|\bigr)\mathcal F
 \ge-C\mathcal F.
 \label{eq:gradient-drift-w}
\end{eqnarray}

Define the function 
\begin{eqnarray}\label{eqn:C1-defn-P}
    P\coloneqq |\nabla w|^2 e^{\eta(w)}, \qquad \text{ in } \bar\Omega,
\end{eqnarray}  where
\begin{eqnarray*}
 \eta(w)\coloneqq -\frac12\log(S_0-w),
\end{eqnarray*} and  $S_0\coloneqq \|w\|_{C^0}+1$. 
Moreover,
\begin{eqnarray}\label{eqn:eta-ineq}
    0<c_0\le\eta'\le C ,\qquad \text{ and } \qquad  \eta''=2(\eta')^2.
\end{eqnarray}

If $ P$ has a maximum value of zero, there is nothing to prove. Otherwise its positive maximum cannot lie on $\Sigma$.
Indeed, strict geodesic convexity and \eqref{eqn:w-nu=0}
give, at every boundary point where $P>0$,
\begin{eqnarray*}
 P_\nu
 =-2e^{\eta(w)}h(\nabla_\Sigma w,\nabla_\Sigma w)
 \le-2\left(
 \min_{\substack{x\in\Sigma, 
 \tau\in T_x\Sigma,\ |\tau|=1}}
 h_x(\tau,\tau)\right)P<0.
\end{eqnarray*}
This contradicts the nonnegative outward derivative
at a boundary maximum.

  At an interior maximum point $x_0\in\Omega$ of $P$, choose
an orthonormal frame with $w_1=|\nabla w|>0$ and
$w_\alpha=0$ for $2\le\alpha\le n$. Rotating
$e_2,\ldots,e_n$ if necessary, we may also assume that
$(K_{\alpha\beta})_{2\le\alpha,\beta\le n}$ is diagonal.
At $x_0$, we have
\begin{eqnarray}\label{eqn:fermat-lemma-C1}
 0= \nabla \log P=   \nabla(\log |\nabla w|^2+\eta(w)),
\end{eqnarray}
which follows
\begin{eqnarray*}
 w_{11}=-\frac12\eta' w_1^2,\qquad w_{1\alpha}=0, \qquad 2\leq \alpha\leq n. 
\end{eqnarray*}
Consequently,
\begin{eqnarray}\label{eqn:K-11}
     K_{11}=-\frac12\eta' w_1^2+O(w_1)+O(1).
\end{eqnarray}
 If $w_1$ is large enough to exceed a fixed constant, then \eqref{eqn:K-11} implies $$K_{11}<0.$$ For $k=n$, this contradicts $K\in \Gamma_n$. For $k<n$,
Lemma~\ref{lem:spectral} (4) gives
$$F^{11}\ge\mathcal{F} /(n-k+1).$$ Together with \eqref{eq:gradient-drift-estimate}, 
\eqref{eq:gradient-drift-w}, \eqref{eqn:fermat-lemma-C1} and \eqref{eqn:eta-ineq}, we obtain, for $w_1\ge1$,
\begin{eqnarray*}
 0&\ge &\mathcal L(\log |\nabla w|^2+\eta(w)) \ge (\eta''-(\eta')^2) |\nabla w|^2 F^{11}
          -C\mathcal{F} -\frac{2|\nabla(\rho s)|}{|\nabla w|}\\
  &\ge & \frac{c_0^2}{n-k+1}|\nabla w|^2 \mathcal{F} -C\mathcal{F} -\frac{C}{|\nabla w|}.
\end{eqnarray*} 
Together with Lemma \ref{lem:spectral} (2), we therefore derive the upper bounds of $|\nabla w|^2(x_0)$, which in turn also bounds $\max_{\bar \Omega} P$. In view of
\begin{eqnarray*}
\nabla U=\rho^{-1}
\bigl(\nabla w-U\nabla\rho-\nabla(bd)\bigr), 
\end{eqnarray*}
 we obtain
the desired estimate \eqref{eqn:C1-bound-U}.
\end{proof}

We refine the gradient estimate to the normalized admissible solution of \eqref{eq:main-problem}.
\begin{prop}\label{prop:c1}
Every $k$-admissible solution $u$ of \eqref{eq:main-problem} satisfying 
\begin{eqnarray}\label{eqn:normalization}
\int_\Omega u/V d\mu_g =0 
\end{eqnarray}
satisfies 
\begin{eqnarray}\label{eq:c1-estimate}
|\lambda|+\|u\|_{C^1(\bar\Omega)}\le C,
\end{eqnarray}
where the constant $C>0$ depends only on $n, k$, $\Omega$, $V$, 
$\|f^{1/k}\|_{C^1}$, and $\|\varphi\|_{C^3}$.
\end{prop}

Oscillation control is essential for Neumann problems whose solutions
are determined only up to a kernel; compare \cite{QX}.
Here the invariant quantity is $u/V$. The estimate for its normalized
oscillation is uniform as the rescaled right hand side tends to zero. 
\begin{proof}
Lemma~\ref{lem:c-bound} gives the bounds for $\lambda$ in \eqref{eq:c1-estimate}. In order to show the upper bound of $\|u\|_{C^1(\bar \Omega)}$, we denote $$v\coloneqq  u/V, \qquad\text{ and } M\coloneqq \operatorname{osc}_{\bar\Omega}v.$$ The normalization \eqref{eqn:normalization} gives
\begin{eqnarray}\label{eqn:normalization-min-max}
    \min_{\bar \Omega} v\leq 0\leq \max_{\bar\Omega} v,
\end{eqnarray}

If $M\leq 1$, it follows
\begin{eqnarray*}
 \|u\|_{C^0(\bar\Omega)}
 \le\|V\|_{C^0(\bar\Omega)}
       \|v\|_{C^0(\bar\Omega)}
 \le\|V\|_{C^0(\bar\Omega)}M
 \le C.
\end{eqnarray*}
Applying Lemma~\ref{lem:gradient} with $U=u$ yields the desired $C^1$ estimate in \eqref{eq:c1-estimate}.

We may therefore assume that $M>1$ from now on, and define
\begin{eqnarray*}
 z\coloneqq \frac{v-\sup_{\bar\Omega}v}{M},\qquad U\coloneqq Vz= \frac u M- \frac V M \sup_{\bar\Omega}\left(\frac u V\right).
\end{eqnarray*}
Then $-1\le z\le0$, $\operatorname{osc}_{\bar \Omega}z=1$, and
\begin{eqnarray*} 
\begin{cases}
F(A[U])=\frac{f^{1/k}}{M}  \qquad &\text{ in }\Omega,\\
U_\nu-\frac{V_\nu}{V}U=\frac{\lambda+\varphi }{M} \qquad &\text{ on }\Sigma,
\end{cases}
\end{eqnarray*}
The hypotheses of Lemma~\ref{lem:gradient} hold uniformly in $M\in(1,+\infty)$, since 
\begin{eqnarray*}
    \left\|\frac{f^{1/k}}{M} \right\|_{C^1(\bar \Omega)}+\left\|\frac{\lambda+\varphi}{M} \right\|_{C^3(\Sigma)}\leq C, 
\end{eqnarray*}then  Lemma~\ref{lem:gradient} implies 
\begin{eqnarray*}
    \|\nabla U\|_{C^0}\leq C,
\end{eqnarray*}which in turn follows 
\begin{eqnarray}\label{eqn:c1-bound-z}
    \|\nabla z\|_{C^0}=\left\|\frac{\nabla U}{V}-z\nabla \log V \right\|_{C^0} \le C.
\end{eqnarray} From $\Gamma_k\subset\Gamma_1$, 
\begin{eqnarray}\label{eqn:div-V2-nabla-z}
 \dive(V^2\nabla z)=V\tr A[U]\ge0 \qquad \text{ in }\Omega, 
 \qquad
 V^2z_\nu=\frac{V(\lambda+\varphi )}{M}, \qquad \text{ on }\Sigma. 
\end{eqnarray}
Integration by parts therefore yields
\begin{eqnarray}\label{eqn:div-V2-2}
 0\le\int_\Omega \dive(V^2\nabla z) d\mu_g 
   =\frac1M\int_\Sigma V(\lambda+\varphi ) d\sigma_g 
   \le\frac CM.
\end{eqnarray}
Multiplication \eqref{eqn:div-V2-nabla-z} by $z$, integration by parts and using \eqref{eqn:div-V2-2} imply
\begin{eqnarray*}
 \int_\Omega V^2|\nabla z|^2 d\mu_g 
 =\frac1M\int_\Sigma Vz(\lambda+\varphi ) d\sigma_g 
       -\int_\Omega z \dive(V^2\nabla z) d\mu_g 
 \le\frac CM.
\end{eqnarray*}
Denote $\bar z\coloneqq |\Omega|^{-1}\int_\Omega z d\mu_g$.
The Poincar\'e inequality and the preceding estimate give
\begin{eqnarray*}
 \int_\Omega |z-\bar z|^2 d\mu_g
 \le C\int_\Omega |\nabla z|^2 d\mu_g\le C\int_\Omega V^2|\nabla z|^2 d\mu_g
 \le \frac{C}{M}.
\end{eqnarray*}
Choose $x_0\in\bar\Omega$ such that
$|z(x_0)-\bar z|=\|z-\bar z\|_{C^0(\bar\Omega)}$.
Since $\operatorname{osc}_{\bar\Omega}z=1$,
\begin{eqnarray*}
 1=\max_{\bar\Omega}z-\min_{\bar\Omega}z\le2\|z-\bar z\|_{C^0(\bar\Omega)}
 =2|z(x_0)-\bar z|.
\end{eqnarray*}
Uniform interior and boundary charts, together with \eqref{eqn:c1-bound-z}, provide constants $r_0,c_\Omega>0$ and $L\ge1$,
independent of $M$, such that
\begin{eqnarray*}
 \operatorname{Vol}_g(\Omega\cap B_r(x))
 &\ge&c_\Omega r^n,
 \qquad x\in\bar\Omega,\quad 0<r<r_0,\\
 |z(x)-z(y)|
 &\le&L d_g(x,y),
 \qquad x,y\in\bar\Omega,\quad d_g(x,y)<r_0.
\end{eqnarray*}
Here $B_r(x)$ denotes the ambient geodesic ball; the second
estimate follows by joining points within the corresponding
interior or flattened boundary chart.

Denote
\begin{eqnarray*}
 r_*\coloneqq\min\left\{\frac{r_0}{2},\frac{1}{4L}\right\}>0.
\end{eqnarray*}
Then, for every $y\in\bar\Omega\cap B_{r_*}(x_0)$,
\begin{eqnarray*}
 |z(y)-\bar z|
 \ge|z(x_0)-\bar z|-|z(y)-z(x_0)|\ge\frac12-Lr_*
 \ge\frac14.
\end{eqnarray*}
Consequently,
\begin{eqnarray*}
 \frac{C}{M}
 &\ge&\int_\Omega |z-\bar z|^2 d\mu_g
 \ge\int_{\Omega\cap B_{r_*}(x_0)}
                  |z-\bar z|^2 d\mu_g\\
 &\ge&\frac1{16}\operatorname{Vol}_g
                 \left  (\Omega\cap B_{r_*}(x_0) \right)
 \ge\frac{c_\Omega r_*^n}{16}.
\end{eqnarray*}
Since $r_*$ and $c_\Omega$ are independent of $M$, this yields
$$ \operatorname{osc}_{\bar\Omega} v=M\le C.$$  
Together with \eqref{eqn:normalization-min-max}, we derive 
\begin{eqnarray}\label{eqn:C0-bound}
\|u\|_{C^0(\bar\Omega)}\le\|V\|_{C^0(\bar\Omega)}\|v\|_{C^0(\bar\Omega)}\le\|V\|_{C^0(\bar\Omega)}\operatorname{osc}_{\bar\Omega}v\le C. \end{eqnarray}
Moreover, the prescribed bounds on the data and
Lemma~\ref{lem:c-bound} give
\begin{eqnarray*}
 \|f^{1/k}\|_{C^1(\bar\Omega)}
 \le \|f\|_{C^0(\bar\Omega)}^{1/k}
 +\frac1 k
   \left(\inf_{\bar\Omega}f\right)^{1/k-1}
   \|\nabla f\|_{C^0(\bar\Omega)}
 \le C,
 \end{eqnarray*}
 and 
 \begin{eqnarray*}
 \|\lambda+\varphi\|_{C^3(\Sigma)}
 \le |\lambda|+\|\varphi\|_{C^3(\Sigma)}
 \le C.
\end{eqnarray*}
We now apply Lemma~\ref{lem:gradient} with
$U=u$, $s=f^{1/k}$, and $b=\lambda+\varphi$, obtaining
\begin{eqnarray*}
 \|\nabla u\|_{C^0(\bar\Omega)}\le C.
\end{eqnarray*}
Together with the $C^0$ bound \eqref{eqn:C0-bound}, we prove 
\eqref{eq:c1-estimate}.
\end{proof}

\subsection{Global \texorpdfstring{$C^2$}{} estimate}\ 

Before estimating second-order derivatives for \eqref{eq:main-problem}, we recall a spectral inequality, which will be used later.
\begin{lem} \label{lem:top-cluster}
Let $\lambda\in\Gamma_k$, let $\lambda_1=\max_i\lambda_i$, and suppose
$\lambda_j\ge c_0\lambda_1$ for some $c_0\in(0,1]$. Then
\begin{eqnarray}\label{eqn:lemma-Fj-leq-F1}
 F_j\le C(n,k,c_0)F_1.
\end{eqnarray}
\end{lem}
\begin{proof}
The assertion is immediate if $j=1$ or $k=1$. If $k=n$, then all the
eigenvalues are positive and $F_j/F_1=\lambda_1/\lambda_j\le c_0^{-1}$. From now on, we suppose that 
$2\le k<n$. By denoting 
$\lambda|1j$ be the $\lambda$ after deleting the first and the $j$-th entries. It follows that  $\sigma_{k-2}(\lambda|1j)=1$ when $k=2$. Recall that 
\begin{eqnarray}\label{eqn:two-identities-lemma3.4}
\frac{\partial\sigma_k }{\partial \lambda_1}=\sigma_{k-1}(\lambda|1j)+\lambda_j \sigma_{k-2}(\lambda|1j),\qquad \frac{\partial\sigma_k }{\partial \lambda_j}=\sigma_{k-1}(\lambda|1j)+\lambda_1 \sigma_{k-2}(\lambda|1j).
\end{eqnarray}
The Newton--MacLaurin inequality for $\lambda|1j\in \R^{n-2}$ gives 
\begin{eqnarray}\label{eqn:newton-mac-ineq}
 \sigma_{k-2}(\lambda|1j) \sigma_k(\lambda|1j) \le\frac{(k-1)(n-k-1)}{k(n-k)}\sigma_{k-1}(\lambda|1j)^2
 \quad (k\le n-2).
\end{eqnarray}
 When $k=n-1$,
$\sigma_k(\lambda|1j)=0$, and we use the same inequality with coefficient zero. Denote the
coefficient in either case by $\rho\coloneqq \frac{(k-1)(n-k-1)}{k(n-k)}\in[0,1)$. Expansion of $\sigma_k(\lambda)$ and using \eqref{eqn:newton-mac-ineq} give
\begin{eqnarray}
 0<\sigma_{k-2}(\lambda|1j) \sigma_k(\lambda)
& =& \sigma_{k-2}(\lambda|1j) \Big( \sigma_k(\lambda|1j) +(\lambda_1+\lambda_j)\sigma_{k-1}(\lambda|1j)+\lambda_1\lambda_j\sigma_{k-2}(\lambda|1j)\Big)   \notag
 \\&\le & \frac{\partial\sigma_k }{\partial \lambda_1}\frac{\partial\sigma_k }{\partial \lambda_j}-(1-\rho)\sigma_{k-1}(\lambda|1j)^2. \label{eqn:lem3.4-expansion}
\end{eqnarray}
We split the proof into the following three cases to finish the proof. 

\begin{enumerate}[label=(\arabic*), leftmargin=*, labelindent=0pt]
 \item 
If $\sigma_{k-1}(\lambda|1j)\ge0$, it follows from \eqref{eqn:two-identities-lemma3.4} that  $$\frac{\partial\sigma_k }{\partial \lambda_1}/\frac{\partial\sigma_k }{\partial \lambda_j}\geq \frac{\lambda_j}{\lambda_1}\geq  c_0.$$ 

\item If
$-\lambda_j \sigma_{k-2}(\lambda|1j)/2\le \sigma_{k-1}(\lambda|1j)<0$, it follows from \eqref{eqn:two-identities-lemma3.4} that
$$\frac{\partial\sigma_k }{\partial \lambda_1}/\frac{\partial\sigma_k }{\partial \lambda_j}\geq  \frac{\lambda_j/2}{\lambda_1} \geq  c_0/2.$$

\item In the remaining case, $\frac{\partial\sigma_k }{\partial \lambda_j}\le\lambda_1 \sigma_{k-2}(\lambda|1j)$ and \eqref{eqn:lem3.4-expansion} yield
\begin{eqnarray*}
  {\frac{\partial\sigma_k }{\partial \lambda_1}} /{\frac{\partial\sigma_k }{\partial \lambda_j}}>
 (1-\rho)\frac{\sigma_{k-1}(\lambda|1j)^2}{\left(\frac{\partial\sigma_k }{\partial \lambda_j}\right)^2}
 \ge{(1-\rho)c_0^2}/{4}.
\end{eqnarray*}
\end{enumerate}
In view of $F_j/F_1=\frac{\partial\sigma_k }{\partial \lambda_j}/\frac{\partial\sigma_k }{\partial \lambda_1}$, the conclusion follows. 
\end{proof}

 We now reduce the global $C^2$ estimate to a bound for the
double normal derivative on $\Sigma$, following the strategy in 
\cite{LTU,MQ} and \cite{CZ}.
For the correction of the mixed normal--tangential
components on a Riemannian manifold, see also
\cite[Section~1.3]{Hanani}.
From now on, we denote
\begin{eqnarray}\label{eqn:defn-M}
 M\coloneqq\max_\Sigma|A[u](\nu,\nu)|,
\end{eqnarray}
and for some small $\delta>0$, 
\begin{eqnarray}\label{eqn:defn-collar}
 \Omega_{\delta}\coloneqq\{x\in \Omega:  0<\text{dist}(x,\Sigma)<\delta\}.   
\end{eqnarray}
Choose $d\in C^{4,\alpha_0}(\bar\Omega)$ which agrees
with the interior distance to $\Sigma$ in
$\bar\Omega_{2\delta_0}$ and vanishes outside a slightly
larger neighborhood of $\Sigma$. The normal field and the projection below
are used only in $\bar\Omega_{2\delta_0}$, where
$|\nabla d|=1$. Define
\begin{eqnarray}\label{eqn:defn-distance}
 N\coloneqq -\nabla d,\qquad
 PX\coloneqq X-\langle X,N\rangle N,\qquad
 \bar h\coloneqq -\nabla^2d,
\end{eqnarray}where $X\in T\Omega$ is any ambient vector field. 
In this thin region, $|\nabla d|=1$, and hence
\begin{eqnarray}\label{eqn:distance-property}
 N|_\Sigma=\nu,\qquad \nabla_NN=0,\qquad
 \bar h(N,\cdot)=0,\qquad
 \bar h|_{T\Sigma\times T\Sigma}=h.
\end{eqnarray}
Extend $\frac{V_\nu}{V}$ and $\varphi$ constantly along normal
geodesics from $\Sigma$ to $\Omega_{\delta}$.

\begin{lem}\label{lem:reduction}
Let $(u,\lambda)$ be a $k$-admissible solution of
\eqref{eq:main-problem} with a uniform bound for
$|\lambda|+\|u\|_{C^1(\bar\Omega)}$. Suppose that
$h>0$ and $h-\frac{V_\nu}{V}g_\Sigma>0$ on $\Sigma$. Then
\begin{eqnarray}\label{eq:global-normal-reduction}
 \|A[u]\|_{C^0(\bar\Omega)}\le C(1+M),
\end{eqnarray}
where $C>0$ depends only on $n,k,\Omega,V$,
$(\inf_{\bar\Omega}f)^{-1}$, $\|f\|_{C^2}$,
$\|\varphi\|_{C^3(\Sigma)}$, and the given $C^1$ and $\lambda$ bounds.
\end{lem}
\begin{proof}
Denote
\begin{eqnarray*}
\kappa_0\coloneqq\min_{\substack{x\in\Sigma\\|\tau|=1}}
 \left(h-\frac{V_\nu}{V}g_\Sigma\right)(\tau,\tau)>0.
\end{eqnarray*}
Choose $\gamma\coloneqq 1+\varepsilon$ with
$0<\varepsilon<\kappa_0/ \left(1+\max_\Sigma|\frac{V_\nu}{V}|\right)$, then 
\begin{eqnarray}\label{eq:reduction-boundary-gap}
 2h-(\gamma+1)\frac{V_\nu}{V}g_\Sigma
 =2\left(h-\frac{V_\nu}{V}g_\Sigma\right)
   -\varepsilon\frac{V_\nu}{V}g_\Sigma
 \ge\kappa_0g_\Sigma.
\end{eqnarray}

Adding a fixed multiple of $V$ to $u$ makes $u\le-K<0$ while
preserving $A[u]$, the Neumann boundary condition in \eqref{eq:main-problem}, and a fixed $C^1$ bound.
Choose $K\gg 1$ large enough that
\begin{eqnarray}\label{eq:reduction-static-translation}
 \left|\frac{\gamma f^{1/k}}u\right| \leq  \frac{\gamma f^{1/k}} K \leq \frac{\varepsilon F(g)}8, \qquad \left|\frac{\gamma(\lambda+\varphi)}u\right| \leq  \frac{\left|\gamma(\lambda+\varphi)\right|}  K   \le\frac{\kappa_0}4.
\end{eqnarray}
Let $\{e_\alpha\}_{\alpha=1}^{n-1}$ be a tangential orthonormal frame on $\Sigma$, differentiating the boundary condition in \eqref{eq:main-problem} gives
\begin{eqnarray}
 A_{\alpha\nu}
 &=&\nabla^2u(e_\alpha,\nu)
   =e_\alpha(u_\nu)-h_{\alpha\beta}u_\beta \notag\\
 &=&\frac{V_\nu}{V}u_\alpha
  +u\nabla_\alpha\left(\frac{V_\nu}{V}\right)
  +\nabla_\alpha\varphi-h_{\alpha\beta}u_\beta  .\label{eq:reduction-mixed-correction}
\end{eqnarray}
For every ambient vector $X$, in view of \eqref{eqn:defn-distance}, define the differential one-form $\mu$ as 
\begin{eqnarray}\label{eqn:defn-mu}
 \mu(X)\coloneqq \frac{V_\nu}{V}\mathrm{d}u(PX)+u \mathrm{d}\left(\frac{V_\nu}{V}\right)(PX)
       +\mathrm{d}\varphi(PX)-\bar h(P\nabla u,PX).
\end{eqnarray}
This is a well-defined one-form, with $\mu(N)=0$ because $PN=0$.
On $\Sigma$, the tangential components of $\mu$ equal the right hand side
of \eqref{eq:reduction-mixed-correction}, namely
\begin{eqnarray}\label{eqn:mu-alpha-A-mixed}
    \mu(e_\alpha) =A_{\alpha \nu},~~~~\text{ on }\Sigma,~~~1\leq \alpha\leq n-1. 
\end{eqnarray}

Choose a fixed smooth cutoff $\chi$ equal to one near
$\Sigma$ and supported in $\bar\Omega_{\delta_0}$.
Define the following tensor in the thin region and extend it
by zero to the rest of $\bar\Omega$:
\begin{eqnarray}\label{eqn:defn-B}
 B\coloneqq \chi(x) \bigl(N^\flat\otimes\mu+\mu\otimes N^\flat\bigr),
 \qquad \widehat A\coloneqq A-B.
\end{eqnarray}
 The $C^1$ bound for
$u$ implies $|B|\le C$, and $\tr B=2\chi(x) \mu(N)=0$  for all $x\in \bar\Omega$. On $\Sigma$, there holds 
\begin{eqnarray*}
B_{\alpha\nu}=\mu(e_\alpha), \qquad   B_{\alpha\beta}=0, \qquad B_{\nu\nu}=0,
\end{eqnarray*}together with \eqref{eqn:defn-mu} and \eqref{eqn:mu-alpha-A-mixed} yield
\begin{eqnarray*}
 \widehat A_{\alpha\nu}=0,\qquad
 \widehat A_{\alpha\beta}=A_{\alpha\beta},\qquad
 \widehat A_{\nu\nu}=A_{\nu\nu}, ~~\forall 1\leq \alpha,\beta\leq n-1. 
\end{eqnarray*}
We \textbf{claim} that, in a geodesic orthonormal frame $\{e_i\}_{i=1}^n$ diagonalizing
$A$ such that $u_{ij}-u \delta_{ij}=\lambda_i\delta_{ij}$, then 
\begin{eqnarray}\label{eq:reduction-correction-bounds}
 |\nabla_iB|\le C(1+|\lambda_i|),\qquad
 |\mathscr LB|\coloneqq  |F^{ij}\nabla_{ij} B|\le
 C\left(1+\mathcal F+\sum_iF_i|\lambda_i|\right).
\end{eqnarray}

Indeed, we rewrite the explicit construction of \eqref{eqn:defn-B} as
\begin{eqnarray*}
 B=\mathcal P^\ell u_\ell+\mathcal C u+\mathcal D,
\end{eqnarray*}
where the coefficient tensors are independent of $u$ and have
uniformly bounded covariant derivatives through order two.
They depend only on $\Omega$, $V$ and $\varphi$.
Hence the first estimate in \eqref{eq:reduction-correction-bounds} follows directly. We only need to show the second inequality in \eqref{eq:reduction-correction-bounds}. Using \eqref{eqn:Lu-eq} and
\eqref{eq:gradient-commutation}, we compute
\begin{eqnarray*}
 \mathscr LB
 &= &(\mathscr L\mathcal P^\ell)u_\ell
     +2F^{ij}(\nabla_i\mathcal P^\ell)u_{\ell j}
     +\mathcal P^\ell
        \bigl[\nabla_\ell(f^{1/k})+F^{\ell j}u_j\bigr]\\
 &&+(\mathscr L\mathcal C)u
     +2F^{ij}(\nabla_i\mathcal C)u_j
     +\mathcal C(f^{1/k}+u\mathcal F)
     +\mathscr L\mathcal D.
\end{eqnarray*}
The bounded coefficient derivatives give
\begin{eqnarray*}
 |\mathscr L\mathcal P|+|\mathscr L\mathcal C|
   +|\mathscr L\mathcal D|\le C\mathcal F.
\end{eqnarray*}
Moreover, in the chosen frame,
\begin{eqnarray*}
 \left|2F^{ij}(\nabla_i\mathcal P^\ell)u_{\ell j}\right|
 \le C\sum_iF_i|\lambda_i+u|
 \le C\left(\mathcal F+\sum_iF_i|\lambda_i|\right).
\end{eqnarray*}
The remaining terms are bounded by $C(1+\mathcal F)$, using
the $C^1$ bound for $u$ and the $C^1$ bound for $f^{1/k}$.
This proves \eqref{eq:reduction-correction-bounds}.

Denote $\Lambda\coloneqq \lambda_{\max}(\widehat A)>0$ and consider the function
\begin{eqnarray}\label{eq:reduction-test}
 Q\coloneqq \log\Lambda+\gamma\log(-u)+\frac{\theta_0}{2(1+M)} |\nabla u|^2,
\end{eqnarray}
where $\theta_0>0$ is a constant to be determined below.  
It suffices to estimate $\Lambda$ at a maximum point of $Q$. Suppose $Q$ attains its maximum value at some point $x_0\in \bar \Omega$. 

\ 

\noindent{\it Case 1: $x_0\in \Sigma$.}
Since $\widehat A_{\alpha\nu}=0$ on $\Sigma$, its largest
eigenvalue is attained in either the normal or a tangential direction.
If $\Lambda(x_0)=A_{\nu\nu}(x_0)$, then
\begin{eqnarray*}
 \Lambda(x_0)=A_{\nu\nu}(x_0)\le M,
\end{eqnarray*} this proves \eqref{eq:global-normal-reduction}.
Otherwise, choose a unit tangential eigenvector $\tau$ corresponding
to $\Lambda(x_0)$. Extend $\tau$ smoothly as a unit tangential field
along $\Sigma$, and then parallel along the normal geodesics.
Thus $\nabla_\nu\tau=0$ and $\tau$ remains orthogonal to $N$.
The definition of $B$ in \eqref{eqn:defn-B} consequently gives
\begin{eqnarray}\label{eqn:hopf-B-tau-tau}
 B(\tau,\tau)=0,\qquad
 (\nabla_\nu B)(\tau,\tau)=0.
\end{eqnarray}
After shrinking the neighborhood of $x_0$, we have
$\widehat A(\tau,\tau)>0$ and may define
\begin{eqnarray*}
\wt Q\coloneqq\log\bigl(\widehat A(\tau,\tau)\bigr)
       +\gamma\log(-u)+\frac{\theta_0}{2(1+M)} |\nabla u|^2.
\end{eqnarray*}
It follows 
\begin{eqnarray*}
\wt Q\le Q,\qquad \wt Q(x_0)=Q(x_0).
\end{eqnarray*}
Hence $\wt Q$ also has a local maximum at $x_0$ relative to
$\bar\Omega$, and 
\begin{eqnarray}\label{eqn:hopf-lemma}
\wt Q_\nu(x_0)\ge0.
\end{eqnarray}
By the Codazzi identity for $A$ and \eqref{eqn:mu-alpha-A-mixed},
\begin{eqnarray*}
 A_{\tau\tau;\nu}
=A_{\nu\tau;\tau}=(\nabla^\Sigma_\tau\mu)_\tau
       -h_{\tau\beta}A_{\beta\tau}
       +h_{\tau\tau}A_{\nu\nu}.
\end{eqnarray*}
Tangential differentiation of \eqref{eqn:defn-mu} for $X=\tau$ yields
\begin{eqnarray*}
 (\nabla^\Sigma_\tau\mu)_\tau
 &=&\frac{V_\nu}{V}(\nabla^2_\Sigma u)_{\tau\tau}
       -h_{\tau\beta}(\nabla^2_\Sigma u)_{\beta\tau} +2u_\tau\nabla^\Sigma_\tau\!\left(\frac{V_\nu}{V}\right)
       +u\left(\nabla^2_\Sigma\frac{V_\nu}{V}\right)_{\tau\tau}\\
 &&+(\nabla^2_\Sigma\varphi)_{\tau\tau}
       -(\nabla^\Sigma_\tau h)_{\tau\beta}u_\beta.
\end{eqnarray*}
Using the $C^1$ bound for $u$ and
\begin{eqnarray*}
 \nabla^2_\Sigma u&=&A|_{T\Sigma}+ug_\Sigma-hu_\nu,\\
 A(\tau,e_\beta)&=&\Lambda g(\tau,e_\beta),
\end{eqnarray*}
we obtain
\begin{eqnarray}
 A_{\tau\tau;\nu}
 &=&\frac{V_\nu}{V}A_{\tau\tau}
       -2h_{\tau\beta}A_{\beta\tau}
       +h_{\tau\tau}A_{\nu\nu}+O(1)
 \nonumber\\
 &\le&\left(\frac{V_\nu}{V}-2h_{\tau\tau}\right)\Lambda
       +C(1+M).
 \label{eq:reduction-tangential-normal}
\end{eqnarray}
Moreover, 
\begin{eqnarray*}
  \left|\partial_\nu\left(\frac 12|\nabla u|^2\right)\right|
 = \left| u_\nu(A_{\nu\nu}+u)+\sum_{\alpha=1}^{n-1} u_\alpha\mu_\alpha \right| \leq C(1+M).
\end{eqnarray*} 
Together with the boundary condition in \eqref{eq:main-problem}, \eqref{eq:reduction-tangential-normal} and  \eqref{eqn:hopf-B-tau-tau}, it follows from \eqref{eqn:hopf-lemma} that
\begin{eqnarray*}
 0\le \wt Q_\nu
 &=&\frac{A_{\tau\tau;\nu}}{\Lambda}
       +\gamma\frac{u_\nu}{u}
       + \frac{\theta_0}{1+M} \partial_\nu\left(\frac12|\nabla u|^2\right)\\
 &\le&(\gamma+1)\frac{V_\nu}{V}-2h_{\tau\tau}
       +\frac{\gamma(\lambda+\varphi)}{u}
       +\frac{C(1+M)}{\Lambda}+C\theta_0.
\end{eqnarray*}
 By choosing $\theta_0>0$ sufficiently small such that
\begin{eqnarray}\label{eqn:first-choice-theta0}
    C\theta_0\le\kappa_0/4,
\end{eqnarray} in view of \eqref{eq:reduction-boundary-gap} and  \eqref{eq:reduction-static-translation}, we derive that 
\begin{eqnarray*}
 0 \leq \wt Q_\nu  &\leq &  -\kappa_0 +\frac{\kappa_0}{4} + \frac{C(1+M)}{\Lambda}+ C \theta_0 
 \\&\le&-\frac{\kappa_0}{2}+\frac{C(1+M)}{\Lambda},
\end{eqnarray*}it follows
\begin{eqnarray}\label{eqn:case1-global-C2}
     \Lambda(x_0) \le C(1+M).
\end{eqnarray}
Thus the desired bound \eqref{eq:global-normal-reduction}  holds in this case.

\ 

\noindent{\it Case 2: $x_0\in \Omega$.} For the moment, assume that $\Lambda(x_0)$ is simple; the general case
will be justified below by a smooth lower support. All local calculations below are evaluated at $x_0$. Choose an orthonormal eigenframe $\{e_i\}_{i=1}^n$ of $A$ at $x_0$,
with the eigenvalues ordered so that
\begin{eqnarray*}
 Ae_i=\lambda_i e_i,\qquad
 \lambda_1\ge\lambda_2\ge\cdots\ge\lambda_n.
\end{eqnarray*}
Complete the unit largest eigenvector $\xi$ of $\widehat A$ to an
orthonormal eigenbasis
$\widehat e_1,\widehat e_2,\ldots,\widehat e_n$ of $\widehat A$, where
\begin{eqnarray*}
 \widehat e_1=\xi,\qquad
 \widehat A\xi=\Lambda\xi,\qquad
 \widehat A\widehat e_\alpha=m_\alpha\widehat e_\alpha,
 \qquad  \text{ with } \quad m_\alpha<\Lambda\quad(2\le\alpha\le n).
\end{eqnarray*}
The first derivative to $Q$ at $x_0$ gives 
\begin{eqnarray}\label{eq:reduction-first-derivative}
\frac{\Lambda_i}{\Lambda}=-\gamma \frac{u_i}u-\frac{\theta_0}{1+M} (\lambda_i+u)u_i.
\end{eqnarray}
The derivative formulas for a simple eigenvalue give
\begin{eqnarray*}
 \Lambda_i=(\nabla_i\widehat A)(\xi,\xi),\qquad 
 \Lambda_{ii}
 =(\nabla_{ii}\widehat A)(\xi,\xi)
       +2\sum_{\alpha=2}^n
       \frac{|(\nabla_i\widehat A)(\xi,\wh e_\alpha)|^2}
            {\Lambda-m_\alpha}.
\end{eqnarray*}
It follows that
\begin{eqnarray}\label{eq:reduction-log-eigenvalue}
 \mathscr L\log\Lambda=\frac{(\mathscr LA)(\xi,\xi)
       -(\mathscr LB)(\xi,\xi)}\Lambda
       +\frac2\Lambda\sum_iF_i\sum_{\alpha=2}^n
       \frac{|(\nabla_i\widehat A)(\xi,\wh e_\alpha)|^2}
            {\Lambda-m_\alpha}
       -\sum_iF_i\left(\frac{\Lambda_i}\Lambda\right)^2. 
\end{eqnarray}
In view of \eqref{eqn:Lu-eq} and \eqref{eq:gradient-commutation}, we further compute to get
\begin{align*}
 \mathscr L\log(-u)
 &=\mathcal F+\frac{f^{1/k}}u-\sum_iF_i\left( \frac{u_i}u\right)^2,\\
 \mathscr L\left(\frac{|\nabla u|^2}{2}\right)
 &=\sum_iF_i\lambda_i^2 +2uf^{1/k}+u^2\mathcal F+\sum_iF_iu_i^2
   +\left\langle\nabla u,\nabla(f^{1/k}) \right\rangle, 
   \end{align*} and
   \begin{align*}
    (\mathscr LA)(\xi,\xi)
 = (f^{1/k})_{;\xi\xi}+f^{1/k}
    -\mathcal F A(\xi,\xi)-F^{ij,rs}A_{ij;\xi}A_{rs;\xi}.
\end{align*}
Since $A(\xi,\xi)=\Lambda+B(\xi,\xi)$,  combining the above identities with 
\eqref{eq:reduction-correction-bounds} and
\eqref{eq:reduction-log-eigenvalue} give
\begin{eqnarray}\label{eq:reduction-interior-raw}
 0&\ge &\mathscr LQ
 =\mathscr L\log\Lambda+\gamma\mathscr L\log(-u)
       +\frac{\theta_0}{1+M} \mathscr L\left(\frac{|\nabla u|^2}{2}\right)
 \nonumber\\
 &=&\varepsilon\mathcal F+\frac{\gamma f^{1/k}}u
       -\gamma\sum_iF_i\left(\frac{u_i}u\right)^2 +\frac{(f^{1/k})_{;\xi\xi}+f^{1/k}
       -\mathcal F B(\xi,\xi)-(\mathscr LB)(\xi,\xi)}{\Lambda}
 \nonumber\\
 && +\frac{\theta_0}{1+M}\left(\sum_iF_i\lambda_i^2+2uf^{1/k}+u^2\mathcal F
       +\sum_iF_i u_i^2
       +\langle\nabla u,\nabla(f^{1/k})\rangle\right)
 \nonumber\\
 &&-\frac{F^{ij,rs}A_{ij;\xi}A_{rs;\xi}}{\Lambda}
       +\frac2\Lambda\sum_iF_i\sum_{\alpha=2}^n
       \frac{|(\nabla_i\widehat A)(\xi,\wh e_\alpha)|^2}
            {\Lambda-m_\alpha}
-\sum_iF_i\left(\frac{\Lambda_i}{\Lambda}\right)^2  
 \nonumber\\
 &\ge&\varepsilon\mathcal F+\frac{\gamma f^{1/k}}u
       +\frac{\theta_0}{1+M} \left(\sum_iF_i\lambda_i^2-C \right)
-\frac C{\Lambda}
       \left(1+\mathcal F+\sum_iF_i|\lambda_i|\right)
 \nonumber\\
 &&-\frac{F^{ij,rs}A_{ij;\xi}A_{rs;\xi}}{\Lambda}
       +\frac2\Lambda\sum_iF_i\sum_{\alpha=2}^n
       \frac{|(\nabla_i\widehat A)(\xi,\wh e_\alpha)|^2}
            {\Lambda-m_\alpha}-\sum_iF_i\left(\frac{\Lambda_i}{\Lambda}\right)^2
       -\gamma\sum_iF_i\left(\frac{u_i}u\right)^2.~~~\qquad
\end{eqnarray}

We next estimate the last line (involving four terms) directly. Let $C_{\rm top}$ be the
constant in \eqref{eqn:lemma-Fj-leq-F1} with $c_0=1/2$, enlarged so that
$C_{\rm top}\ge1$. For $\tau,\delta\in(0,1/10)$, to be fixed below, denote 
\begin{eqnarray*}
 \mathcal I\coloneqq \left\{j:\lambda_j\ge \frac{\lambda_1} 2\right\},\qquad 
 J\coloneqq \left\{i:\lambda_i>-\tau\lambda_1,\quad
       F_i>\frac{C_{\rm top}}\delta F_1\right\}.
\end{eqnarray*}
Lemma~\ref{lem:top-cluster} gives $F_j\le C_{\rm top}F_1$ for
$j\in\mathcal I$. Hence $\mathcal I\cap J=\varnothing$, and, for
$i\in J$ and $j\in\mathcal I$,
\begin{eqnarray}\label{eq:reduction-separated-weights}
 F_i-F_j\ge(1-\delta)F_i,\qquad
 0<\lambda_j-\lambda_i\le(1+\tau)\lambda_1.
\end{eqnarray}

Write $\xi$ in the eigenframe $\{e_i\}_{i=1}^n$ of $A$ as
\begin{eqnarray*}
 \xi=\sum_{j=1}^n\xi_j e_j,\qquad
 \xi_j\coloneqq\langle\xi,e_j\rangle\in\mathbb R,\qquad
 \sum_{j=1}^n\xi_j^2=1.
\end{eqnarray*}
Since $\widehat A=A-B$ and $|B|\le C$, then
\begin{eqnarray}\label{eqn:Lambda-lambda-1}
 |\Lambda-\lambda_1|\le C.
\end{eqnarray}
In particular, for $\Lambda\geq 3C$, we have
\begin{eqnarray*} 
    \lambda_1\geq \Lambda-C \geq \frac{2\Lambda} 3>0,
\end{eqnarray*}
which follows   
\begin{eqnarray}\label{eqn:lambda1-to-Lambda-2}
    \left|\frac{\Lambda}{\lambda_1}-1\right| \leq \frac C {\lambda_1} \leq \frac{C}{\Lambda}.
\end{eqnarray}

Moreover, $\widehat A\xi=\Lambda\xi$ implies
\begin{eqnarray*}
 (A-\Lambda g)\xi=B\xi,
\end{eqnarray*}
and hence, taking the inner product with $e_j$,
\begin{eqnarray}\label{eqn:B-xi-e_j}
 (\lambda_j-\Lambda)\xi_j=\langle B\xi,e_j\rangle.
\end{eqnarray}

Denote $\mathcal I^c\coloneqq\{1,\ldots,n\}\setminus\mathcal I$, and 
define the orthogonal projections of $\xi$ onto
$\operatorname{span}\{e_j:j\in\mathcal I\}$ and its orthogonal
complement by
\begin{eqnarray*}
 \xi_{\mathcal I}\coloneqq\sum_{j\in\mathcal I}\xi_j e_j,
 \qquad
 \xi_{\mathcal I^c}\coloneqq\sum_{j\notin\mathcal I}\xi_j e_j.
\end{eqnarray*}
Thus
\begin{eqnarray*}
&& \xi=\xi_{\mathcal I}+\xi_{\mathcal I^c},\qquad
 \langle\xi_{\mathcal I},\xi_{\mathcal I^c}\rangle=0,\\
&& |\xi_{\mathcal I}|^2=\sum_{j\in\mathcal I}\xi_j^2,\qquad
 |\xi_{\mathcal I^c}|^2=\sum_{j\notin\mathcal I}\xi_j^2.
\end{eqnarray*}
For $j\notin\mathcal I$, using \eqref{eqn:Lambda-lambda-1}, we have $\lambda_j<\frac{\lambda_1}{2}
 \le\frac{\Lambda+C}{2}.$
Consequently, for $\Lambda\ge3C$,
\begin{eqnarray*}
 \Lambda-\lambda_j
 >\frac{\Lambda-C}{2}
 \ge\frac{\Lambda}{3}.
\end{eqnarray*}
Together with \eqref{eqn:B-xi-e_j}, \eqref{eqn:lambda1-to-Lambda-2} and $|\xi|=1$, we obtain
\begin{eqnarray}\label{eq:reduction-eigenvector-localization}
 |\xi_{\mathcal I^c}|^2=\sum_{j\notin\mathcal I}
       \frac{|\langle B\xi,e_j\rangle|^2}
            {(\Lambda-\lambda_j)^2}\le\frac9{\Lambda^2}|B\xi|^2 \le\frac C{\Lambda^2},\qquad \frac{\Lambda}{\lambda_1}=1+O(\Lambda^{-1}).
\end{eqnarray}

In the following calculation, symmetric tensors are identified with
their associated self-adjoint endomorphisms. Denote 
\begin{eqnarray*}
 \eta&\coloneqq&\nabla_\xi A,\qquad
       \zeta\coloneqq\xi_{\mathcal I^c}-|\xi_{\mathcal I^c}|^2\xi,\\
 W_i&\coloneqq&\Pi_{\xi^\perp}\bigl((\nabla_i\widehat A)\xi\bigr),
       \qquad\Pi_{\xi^\perp}X=X-\langle X,\xi\rangle\xi.
\end{eqnarray*}
By \eqref{eq:reduction-eigenvector-localization},
\begin{eqnarray*}
&& \langle\zeta,\xi\rangle=0,\qquad
       |\zeta|\le C\Lambda^{-1},\\
&& \xi_{\mathcal I}
 =(1-|\xi_{\mathcal I^c}|^2)\xi-\zeta.
\end{eqnarray*}
The full symmetry of $\nabla A$ gives
\begin{eqnarray*}
 \eta e_i=(\nabla_\xi A)(e_i)=(\nabla_{e_i}A)(\xi) =(\nabla_iA)\xi
       =(\nabla_i\widehat A)\xi+(\nabla_iB)\xi=\Lambda_i\xi+W_i+(\nabla_iB)\xi.
\end{eqnarray*}
Taking the inner product with $\xi_{\mathcal I}$ therefore yields  
\begin{eqnarray*}
 \langle\xi_{\mathcal I},\eta e_i\rangle
 =(1-|\xi_{\mathcal I^c}|^2)\Lambda_i
       +\langle\xi_{\mathcal I},(\nabla_iB)\xi\rangle
       -\langle\zeta,W_i\rangle.
\end{eqnarray*}
Applying Young's inequality twice and using
\eqref{eq:reduction-correction-bounds}, we obtain
\begin{eqnarray}\label{eq:reduction-projection-bound}
 \langle\xi_{\mathcal I},\eta e_i\rangle^2
 \ge(1-\delta)^2(1-|\xi_{\mathcal I^c}|^2)^2\Lambda_i^2-C_\delta(1+\lambda_i^2)
       -\frac{C_\delta}{\Lambda^2}|W_i|^2.
\end{eqnarray}  
The Cauchy--Schwarz inequality implies
\begin{eqnarray}\label{eq:reduction-projection-bound-2}
\langle\xi_{\mathcal I},\eta e_i\rangle^2
=\left(\sum_{j\in\mathcal I}
\langle\xi,e_j\rangle\eta_{ij}\right)^2\leq
\left(\sum_{j\in\mathcal I}\langle\xi,e_j\rangle^2\right)
\left(\sum_{j\in\mathcal I}\eta_{ij}^2\right)=|\xi_{\mathcal I}|^2
\sum_{j\in\mathcal I}\eta_{ij}^2\leq\sum_{j\in\mathcal I}\eta_{ij}^2.
\end{eqnarray}
The spectral second derivative formula reads
\begin{eqnarray*}
 -F^{ij,rs}\eta_{ij}\eta_{rs}= -\sum_{p,q}\frac{\partial^2F}{\partial\lambda_p\partial\lambda_q}
       \eta_{pp}\eta_{qq}
       +2\sum_{p<q}\frac{F_q-F_p}{\lambda_p-\lambda_q}\eta_{pq}^2.
\end{eqnarray*}
The divided differences at coincident eigenvalues are interpreted by
continuity. Keeping the pairs
from \eqref{eq:reduction-separated-weights}, whose gaps are strictly
positive, using \eqref{eq:reduction-projection-bound-2} gives
\begin{eqnarray*}
 -\frac{F^{ij,rs}\eta_{ij}\eta_{rs}}\Lambda
 &\ge&\frac2\Lambda\sum_{i\in J}\sum_{j\in\mathcal I}
       \frac{F_i-F_j}{\lambda_j-\lambda_i}\eta_{ij}^2 \ge\frac{2(1-\delta)}{(1+\tau)\Lambda\lambda_1}
       \sum_{i\in J}F_i\sum_{j\in\mathcal I}\eta_{ij}^2\\
 &\ge&\frac{2(1-\delta)}{(1+\tau)\Lambda\lambda_1}
       \sum_{i\in J}F_i
       \langle\xi_{\mathcal I},\eta e_i\rangle^2.
\end{eqnarray*}
Here the last inequality uses $|\xi_{\mathcal I}|\le1$.
Combining \eqref{eq:reduction-projection-bound} with
\eqref{eq:reduction-eigenvector-localization}, we conclude that
\begin{eqnarray}\label{eq:reduction-off-diagonal-control}
 -\frac{F^{ij,rs}\eta_{ij}\eta_{rs}}\Lambda
 &\ge&\frac{2(1-\delta)^3(1-|\xi_{\mathcal I^c}|^2)^2}{1+\tau}
       \frac\Lambda{\lambda_1}
       \sum_{i\in J}F_i\left(\frac{\Lambda_i}\Lambda\right)^2
\notag \\&& -\frac{C_\delta}{\Lambda^2}
       \left(\mathcal F+\sum_iF_i\lambda_i^2\right)
       -\frac{C_\delta}{\Lambda^4}\sum_iF_i|W_i|^2
 \nonumber\\
 &\ge&\left(\frac2{1+\tau}-C\delta-\frac C\Lambda\right)
       \sum_{i\in J}F_i\left(\frac{\Lambda_i}\Lambda\right)^2
\notag \\&& -\frac{C_\delta}{\Lambda^2}
       \left(\mathcal F+\sum_iF_i\lambda_i^2\right)
       -\frac{C_\delta}{\Lambda^4}\sum_iF_i|W_i|^2.~\qquad~
\end{eqnarray}
The constant $C$ in the leading coefficient is independent of
$\tau$ and $\delta$.

To absorb the last error term in \eqref{eq:reduction-off-diagonal-control}, use $\operatorname{tr}A>0$ and $|B|\le C$ to obtain
\begin{eqnarray*}
 \lambda_n>-(n-1)\lambda_1,\qquad m_\alpha\ge\lambda_n-C,
 \end{eqnarray*}it yields
 \begin{eqnarray}\label{eqn:difference-Lambda-m-alpha}
 0<\Lambda-m_\alpha\le\lambda_1-\lambda_n+C\le C\Lambda.
\end{eqnarray}
Since $\{\widehat e_2,\ldots,\widehat e_n\}$ form an orthonormal basis of $\xi^\perp$,
\begin{eqnarray*}
 \sum_{\alpha=2}^n
       |(\nabla_i\widehat A)(\xi,\wh e_\alpha)|^2=|W_i|^2.
\end{eqnarray*}
Together with \eqref{eqn:difference-Lambda-m-alpha}, it follows that
\begin{eqnarray}\label{eq:reduction-eigenvector-energy}
 \frac2\Lambda\sum_iF_i\sum_{\alpha=2}^n
       \frac{|(\nabla_i\widehat A)(\xi,\wh e_\alpha)|^2}
            {\Lambda-m_\alpha} 
 \ge\frac c{\Lambda^2}\sum_iF_i|W_i|^2.
\end{eqnarray}
Thus, for $\Lambda$ sufficiently large depending on $\delta$,
\begin{eqnarray*}
 \frac c{\Lambda^2}-\frac{C_\delta}{\Lambda^4}\ge0.
\end{eqnarray*}
Adding \eqref{eq:reduction-eigenvector-energy} to
\eqref{eq:reduction-off-diagonal-control} proves
\begin{eqnarray}
 -\frac{F^{ij,rs}A_{ij;\xi}A_{rs;\xi}}\Lambda
       +\frac2\Lambda\sum_iF_i\sum_{\alpha=2}^n
       \frac{|(\nabla_i\widehat A)(\xi,\wh e_\alpha)|^2}
            {\Lambda-m_\alpha}
& \ge & \left(\frac2{1+\tau}-C\delta-\frac C\Lambda\right)
       \sum_{i\in J}F_i\left(\frac{\Lambda_i}\Lambda\right)^2  \notag \\&&-\frac{C_\delta}{\Lambda^2}
       \left(\mathcal F+\sum_iF_i\lambda_i^2\right). \label{eq:reduction-integrated-spectral}
\end{eqnarray}
For $k=1$, $F_i=F_1$ and $J=\varnothing$; the same calculation also holds.

We now fix the parameters and return to the last line in 
\eqref{eq:reduction-interior-raw}. Choose
\begin{eqnarray*}
 c_*\coloneqq\frac{\gamma+1}{2\gamma},\qquad
 \eta_0\coloneqq\frac{\gamma-1}{2(\gamma+1)}.
\end{eqnarray*}
Then $0<c_*<1$, $0<\eta_0<1$, and
\begin{eqnarray*}
 c_*(1-\eta_0)\gamma^2-\gamma
       =\frac{\gamma(\gamma-1)}4>0.
\end{eqnarray*}
First choose $\tau,\delta\in(0,1/10)$ sufficiently small, and then a
fixed threshold $\Lambda_0$ sufficiently large, so that all preceding
estimates hold for $\Lambda\ge\Lambda_0$ and
\begin{eqnarray*}
 2\tau+C\delta\le\frac{1-c_*}{2},\qquad
 \frac C{\Lambda_0}\le\frac{1-c_*}{2}.
\end{eqnarray*}
Consequently,
\begin{eqnarray*}
 \frac2{1+\tau}-1-C\delta-\frac C\Lambda
       \ge c_*\qquad(\Lambda\ge\Lambda_0).
\end{eqnarray*}
All constants and thresholds are chosen to cover also the bounded
perturbation of $B$ used below. The case $\Lambda<\Lambda_0$ already
has the required bound.

For $i\in J$, \eqref{eq:reduction-first-derivative} and Young's inequality imply
\begin{eqnarray*}
 \left(\frac{\Lambda_i}\Lambda\right)^2
 &=&\left(\gamma\frac{u_i}u+\frac{\theta_0}{1+M}(\lambda_i+u)u_i\right)^2
 \\& \ge &(1-\eta_0)\gamma^2\left(\frac{u_i}u\right)^2
       -(\eta_0^{-1}-1) \left(\frac{\theta_0}{1+M} \right)^2(\lambda_i+u)^2u_i^2,
\end{eqnarray*}it follows
\begin{eqnarray}\label{eqn:global-C2-young-ineq-2}
     c_*\left(\frac{\Lambda_i}\Lambda\right)^2
       -\gamma\left(\frac{u_i}u\right)^2
& \ge &\frac{\gamma(\gamma-1)}4\left(\frac{u_i}u\right)^2
       -C\left(\frac{\theta_0}{1+M} \right)^2(\lambda_i+u)^2u_i^2
\notag  \\& \ge &-C\left(\frac{\theta_0}{1+M} \right)^2(1+\lambda_i^2).
\end{eqnarray}
For $i\notin 
J$, either $\lambda_i\le-\tau\lambda_1$ or
$F_i\le C_{\rm top}F_1/\delta$. In the first case,
\begin{eqnarray*}
 F_i\le\frac{F_i\lambda_i^2}{\tau^2\lambda_1^2}
       \le\frac1{\tau^2\lambda_1^2}\sum_{p=1}^n F_p\lambda_p^2, 
\end{eqnarray*}
in the second case,
\begin{eqnarray*}
 F_i\le\frac{C_{\rm top}}\delta F_1
       \le\frac{C_{\rm top}}{\delta\lambda_1^2}
       \sum_{p=1}^n F_p\lambda_p^2.
\end{eqnarray*}
Since $\lambda_1\ge\Lambda/2$ after increasing $\Lambda_0$, both cases give
\begin{eqnarray}\label{eq:reduction-complement-weights}
 \sum_{i\notin J}F_i
 \le \frac C{\Lambda^2}\sum_iF_i\lambda_i^2.
\end{eqnarray}
Using \eqref{eq:reduction-first-derivative}, the fixed $C^1$ bound, and
$u\le-K<0$, we also have
\begin{eqnarray*}
 \left(\frac{\Lambda_i}\Lambda\right)^2
       +\gamma\left(\frac{u_i}u\right)^2
 &\le&(2\gamma^2+\gamma)\left(\frac{u_i}u\right)^2
       +2\left(\frac{\theta_0}{1+M} \right)^2(\lambda_i+u)^2u_i^2\\
 &\le&C+C\left(\frac{\theta_0}{1+M} \right)^2(1+\lambda_i^2).
\end{eqnarray*}
Thus \eqref{eq:reduction-complement-weights} and \eqref{eqn:global-C2-young-ineq-2} yield
\begin{eqnarray*}
 &&\sum_{i\notin J}F_i\left [
       \left(\frac{\Lambda_i}\Lambda\right)^2
       +\gamma\left(\frac{u_i}u\right)^2\right] \le\frac C{\Lambda^2}\sum_iF_i\lambda_i^2
       +C\left(\frac{\theta_0}{1+M} \right)^2\left(\mathcal F+\sum_iF_i\lambda_i^2\right),\\
 &&\sum_{i\in J}F_i\left[
       c_*\left(\frac{\Lambda_i}\Lambda\right)^2
       -\gamma\left(\frac{u_i}u\right)^2\right] \ge-C\left(\frac{\theta_0}{1+M} \right)^2
       \left(\mathcal F+\sum_iF_i\lambda_i^2\right).
\end{eqnarray*}
Subtracting the gradient terms from
\eqref{eq:reduction-integrated-spectral} and applying these bounds, we obtain
\begin{eqnarray}\label{eq:reduction-third-order-control}
 &&-\frac{F^{ij,rs}A_{ij;\xi}A_{rs;\xi}}\Lambda
       + \frac2\Lambda\sum_iF_i\sum_{\alpha=2}^n
       \frac{|(\nabla_i\widehat A)(\xi,\wh e_\alpha)|^2}
            {\Lambda-m_\alpha}  
       -\sum_iF_i\left(\frac{\Lambda_i}\Lambda\right)^2
       -\gamma\sum_iF_i\left(\frac{u_i}u\right)^2
 \nonumber\\
 &\geq &  c_* \sum_{i\in J} F_i \left(\frac{\Lambda_i}\Lambda\right)^2-\sum_{i\notin J}  F_i \left(\frac{\Lambda_i}\Lambda\right)^2 -\gamma \sum_{i\in J} F_i \left(\frac{u_i}{u}\right)^2 -\gamma \sum_{i\notin J} F_i \left(\frac{u_i}{u}\right)^2
\notag \\ &&{}-\frac{C_\delta}{\Lambda^2}
       \left(\mathcal F+\sum_iF_i\lambda_i^2\right)
\notag \\ & \ge &-C\left[ \Lambda^{-2}+\left(\frac{\theta_0}{1+M} \right)^2 \right]
       \left(\mathcal F+\sum_iF_i\lambda_i^2\right).
\end{eqnarray}

It remains to absorb the lower-order terms in
\eqref{eq:reduction-interior-raw}.  
By Cauchy's inequality and Young's inequality,
\begin{eqnarray*}
 \frac C\Lambda\sum_iF_i|\lambda_i|
 \le \frac C\Lambda\mathcal F^{1/2}
       \left(\sum_iF_i\lambda_i^2\right)^{1/2} \le \frac\varepsilon8\mathcal F
       +\frac{C_\varepsilon}{\Lambda^2}\sum_iF_i\lambda_i^2.
\end{eqnarray*}
Moreover, \eqref{eq:reduction-static-translation} and
Lemma \ref{lem:spectral} imply
\begin{eqnarray*}
 \frac{\gamma f^{1/k}}u
       \ge-\frac\varepsilon8F(g)
       \ge-\frac\varepsilon8\mathcal F.
\end{eqnarray*}
Substituting these estimates and \eqref{eq:reduction-third-order-control}
into \eqref{eq:reduction-interior-raw}, we find
\begin{eqnarray*}
 \mathscr LQ
 &\ge&\frac{3\varepsilon}{4}\mathcal F
       +\frac{\theta_0}{1+M}\sum_iF_i\lambda_i^2
       -\frac{C\theta_0}{1+M}-\frac C\Lambda(1+\mathcal F)\\
 &&-C\left[\Lambda^{-2}
       +\left(\frac{\theta_0}{1+M}\right)^2\right]
       \left(\mathcal F+\sum_iF_i\lambda_i^2\right).
\end{eqnarray*}
Using Lemma \ref{lem:spectral} once more, this becomes,
for a fixed $C_0\ge1$,
\begin{eqnarray}\label{eq:reduction-interior-collected}
 \mathscr LQ
 &\ge&\left\{\frac{3\varepsilon}{4}
       -C_0\left[\frac{\theta_0}{1+M}+\left(\frac{\theta_0}{1+M} \right)^2+\Lambda^{-1}+\Lambda^{-2}\right]\right\}
       \mathcal F
 \nonumber\\
 &&{}+\left\{\frac{\theta_0}{1+M}-C_0\left[\left(\frac{\theta_0}{1+M} \right)^2+\Lambda^{-2}\right]\right\}
       \sum_iF_i\lambda_i^2.
\end{eqnarray}
Choose $\theta_0$ sufficiently small to satisfy
\eqref{eqn:first-choice-theta0} and 
\begin{eqnarray}\label{eq:reduction-interior-parameters}
&&  0<\theta_0
 \le\min\left\{1,\frac\varepsilon{16C_0},\frac1{4C_0}\right\},
 \nonumber\\ && 
 C_*\ge\max\left\{1,\Lambda_0,
       \frac{16C_0}\varepsilon,2\sqrt{C_0}\right\}.
\end{eqnarray}
If $\Lambda\ge C_*/\sqrt{\theta_0/(1+M)}$, then
\eqref{eq:reduction-interior-parameters} gives
\begin{eqnarray*}
 &&\Lambda^{-1}\le\frac{\sqrt{\theta_0/(1+M)}}{C_*}\le\frac1{C_*},\qquad
       \Lambda^{-2}\le\frac{\theta_0/(1+M)}{C_*^2},\\
&&C_0 \left[\frac{\theta_0}{1+M}+\left(\frac{\theta_0}{1+M} \right)^2\right]\le2C_0\theta_0\le\frac\varepsilon8,\\&&
 C_0(\Lambda^{-1}+\Lambda^{-2})
 \le\frac{2C_0}{C_*}\le\frac\varepsilon8,\\&&
 C_0\left(\frac{\theta_0}{1+M} \right)^2 \le \frac{\theta_0}{4(1+M)} ,\qquad
       C_0\Lambda^{-2}\le  \frac{\theta_0}{4(1+M)}.
\end{eqnarray*}
Together with \eqref{eq:reduction-interior-collected} and
Lemma \ref{lem:spectral}, these inequalities imply
\begin{eqnarray*}
 0\ge\mathscr LQ \ge\frac\varepsilon2\mathcal F
       +\frac{\theta_0}{2(1+M)}   \sum_iF_i\lambda_i^2
       \ge\frac\varepsilon2F(g)>0,
\end{eqnarray*}
a contradiction. Hence
\begin{eqnarray}\label{eq:reduction-interior-eigenvalue-bound}
 \Lambda(x_0)< \frac{C_*}{\sqrt{\theta_0}}\sqrt{1+M}
       \le C\sqrt{1+M}.
\end{eqnarray}

We finally remove the simplicity assumption. Choose a unit largest eigenvector $\xi$ at $x_0$ and extend it by radial parallel transport
near $x_0$. Then
\begin{eqnarray*}
 \Pi_{\xi^\perp}=g-\xi^\flat\otimes\xi^\flat,\qquad
       \nabla\Pi_{\xi^\perp}(x_0)=0,\qquad
       |\nabla^2\Pi_{\xi^\perp}(x_0)|\le C.
\end{eqnarray*}
Define locally
\begin{eqnarray*}
 \widetilde B\coloneqq B+\Pi_{\xi^\perp},\qquad
       \widetilde A\coloneqq A-\widetilde B,\qquad
       \widetilde\Lambda\coloneqq \lambda_{\max}(\widetilde A).
\end{eqnarray*}
Then 
\begin{eqnarray*}
 \widetilde\Lambda(x)\le\Lambda(x),\qquad
       \widetilde\Lambda(x_0)=\Lambda(x_0)>0.
\end{eqnarray*}
If $m_2,\ldots,m_n$ are the other eigenvalues of $\widehat A(x_0)$,
the eigenvalues of $\widetilde A(x_0)$ are
\begin{eqnarray*}
 \Lambda(x_0),\quad m_2-1,\ldots,m_n-1,
       \qquad \text{ with } m_\alpha-1<\Lambda(x_0).
\end{eqnarray*}
Thus $\widetilde\Lambda$ is simple and smooth near $x_0$. Moreover,
\begin{eqnarray*}
 |\widetilde B|\le C,\qquad
       |\nabla_i\widetilde B(x_0)|\le C(1+|\lambda_i|),
       \end{eqnarray*} and 
       \begin{eqnarray*}
 |\mathscr L\widetilde B(x_0)|
 \le |\mathscr LB(x_0)|+C\mathcal F\le C\left(1+\mathcal F+\sum_iF_i|\lambda_i|\right).
\end{eqnarray*}
On a neighborhood $U\Subset\Omega$ of $x_0$ where
$\widetilde\Lambda>0$, consider the function
\begin{eqnarray*}
 \widetilde Q\coloneqq\log\widetilde\Lambda
       +\gamma\log(-u)+\frac{\theta_0}{2(1+M)} |\nabla u|^2.
\end{eqnarray*}
Then
\begin{eqnarray*}
 \widetilde Q\le Q  \quad \text{ in }U,\qquad \text{ and } \qquad
       \widetilde Q(x_0)=Q(x_0),
\end{eqnarray*}
       consequently, 
       \begin{eqnarray*}
 \nabla\widetilde Q(x_0)&=&0,\qquad \text{ and }\qquad 
       \mathscr L\widetilde Q(x_0)\le0.
\end{eqnarray*}
The preceding calculations apply with $B,\widehat A,\Lambda,Q$ replaced
by $\widetilde B,\widetilde A,\widetilde\Lambda,\widetilde Q$.
In particular, the projection calculation and the absorption leading
to \eqref{eq:reduction-integrated-spectral} remain valid. Therefore
\eqref{eq:reduction-interior-collected} and
\eqref{eq:reduction-interior-parameters} give the same contradiction
when $\widetilde\Lambda(x_0)\ge C_*/\sqrt{\theta_0/(1+M)}$.
Since $\widetilde\Lambda(x_0)=\Lambda(x_0)$,
\eqref{eq:reduction-interior-eigenvalue-bound} follows without any
simplicity assumption.

\ 

Combining \eqref{eq:reduction-interior-eigenvalue-bound} with \eqref{eqn:case1-global-C2},
we have $$\Lambda(x_0)\le C(1+M).$$ For every $x\in\bar\Omega$,
\eqref{eq:reduction-test} and $Q(x)\le Q(x_0)$ imply
\begin{eqnarray*}
 \Lambda(x)
 &\le&\Lambda(x_0)
       \left(\frac{-u(x_0)}{-u(x)}\right)^\gamma
       \exp\left\{\frac{\theta_0}{2(1+M)} 
       \left[ |\nabla u(x_0)|^2-|\nabla u(x)|^2\right] \right\}\\
 &\le&C\Lambda(x_0)\le C(1+M),
\end{eqnarray*}
where the second inequality uses the fixed $C^1$ bound and $u\le-K<0$.
Thus we conclude that
\begin{eqnarray}\label{eqn:upper-bound-lambda-max}
 \lambda_{\max}(A(x))
\le \Lambda(x)+|B(x)|\le C(1+M).
\end{eqnarray}
Since $A\in\Gamma_k\subset\Gamma_1$,  
\begin{eqnarray*}
 \lambda_{\min}(A(x))\geq -(n-1)\lambda_{\max}(A(x)) \geq -(n-1) C(1+M),
       \end{eqnarray*} together with \eqref{eqn:upper-bound-lambda-max} imply
       \begin{eqnarray*}
 \|A[u]\|_{C^0(\bar\Omega)}\le C(1+M).
\end{eqnarray*}
This proves \eqref{eq:global-normal-reduction}.
Hence we complete the proof. 
\end{proof}

\subsection{Double normal derivative estimate}\ 

We now estimate \eqref{eqn:defn-M} by using the barrier argument; the test functions (cf. \eqref{eq:positive-barrier-definition} and \eqref{eqn:defn-P-barrier-function}) chosen here are motivated by the idea of Ma--Qiu in 
\cite[Lemmas~15 and~16]{MQ} for the Hessian equation of the Neumann problem in a strictly convex domain in Euclidean space; see also \cite[Section~4]{CZ} for an analogous setting of Hessian quotient equations.

\begin{prop}\label{prop:c2}
Under the hypotheses of Theorem~\ref{thm:existence}, a normalized
admissible solution of \eqref{eq:main-problem} satisfies
\begin{eqnarray*}
    M\leq C,
\end{eqnarray*}and consequently, 
\begin{eqnarray}\label{eqn:c2-global-estimate}
 |\lambda|+\|u\|_{C^2(\bar\Omega)}\le C.
\end{eqnarray}
More generally, the same conclusion holds when $h>0$, a uniform
$C^1$ bound for the normalized solution and a bound for $|\lambda|$
are given, and
\begin{eqnarray*}
 \|A[u]\|_{C^0(\bar\Omega)}\le C_0(1+M)
\end{eqnarray*}
holds with $C_0$ independent of $M$. In this case, the constant also
depends on $C_0$ and the given bounds.
\end{prop}
\begin{proof}
Under the hypotheses of Theorem~\ref{thm:existence},
Proposition~\ref{prop:c1} and Lemma~\ref{lem:reduction} give the
bounds in the general statement. We prove the estimate under
these bounds, using a representative with a fixed $C^1$ bound.
All constants below are independent of $M$;
unless a dependence on $\beta$ is indicated, they are also independent
of $\beta$, $\mu$, and $A_0$.
 
Since $h>0$, we may shrink the fixed neighborhood region in
\eqref{eqn:defn-collar} so that, for some fixed $c_d>0$,
\begin{eqnarray*}
 |\nabla d|=1,\qquad
 -\nabla^2d\ge c_d\bigl(g-\mathrm d d\otimes\mathrm d d\bigr),
 \qquad 0\le d\le\delta_0\le\frac14.
\end{eqnarray*}
For $h_0\coloneqq-d+d^2$, it follows that
\begin{eqnarray*}
 \nabla^2h_0
 &=&(-1+2d)\nabla^2d+2 \mathrm d d\otimes\mathrm d d\\
 &\ge&\frac{c_d}{2}\bigl(g-\mathrm d d\otimes\mathrm d d\bigr)
       +2 \mathrm d d\otimes\mathrm d d
 \ge c_0g,
\end{eqnarray*}
where $c_0\coloneqq \min\{c_d/2,2\}>0$. Therefore  
\begin{eqnarray}\label{eq:collar-convex-barrier}
 \mathscr Lh_0
 \ge c_0\mathcal F
 \ge c_h(1+\mathcal F),\qquad
 c_h\coloneqq\frac{c_0F(g)}{1+F(g)}>0.
\end{eqnarray}
Define
\begin{eqnarray*}
 w\coloneqq-\langle\nabla u,\nabla d\rangle
       -\frac{V_\nu}{V}u-\lambda-\varphi
       \qquad\text{in }\bar\Omega_{\delta_0}.
\end{eqnarray*}
In view of \eqref{eqn:defn-distance} and \eqref{eqn:distance-property},
\begin{eqnarray}\label{eqn:w=0-Sigma}
 w|_\Sigma=0,\qquad |w|\le C
       \quad\text{in }\bar\Omega_{\delta_0}.
\end{eqnarray}
At a point in the region $\Omega_{\delta_0}$, choose a normal orthonormal frame diagonalizing
$A$, so that $u_{ij}=(\lambda_i+u)\delta_{ij}$. Differentiation yields
\begin{eqnarray}\label{eq:normal-error-first}
     w_i= -\lambda_i d_i-ud_i-u_jd_{ij}-\left( \frac{V_\nu}{V} u+\varphi\right)_i
  \coloneqq -\lambda_i d_i+E_i
\end{eqnarray}where $E_i$ satisfy $|E_i|\le C$.

 Using \eqref{eqn:Lu-eq} and
\eqref{eq:gradient-commutation}, we compute
\begin{eqnarray*}
 \mathscr Lw
 &=&-2\sum_iF_i\lambda_i d_{ii}
       -\langle\nabla d,\nabla(f^{1/k})\rangle
       -\sum_iF_i d_i u_i\\
 &&-\sum_{i,j}F_i u_jd_{jii}
       -2u\sum_iF_i d_{ii}
       -\frac{V_\nu}{V}(f^{1/k}+u\mathcal F)\\
 &&-u\mathscr L\left(\frac{V_\nu}{V}\right)
       -2\sum_iF_i\left(\frac{V_\nu}{V}\right)_iu_i
       -\mathscr L\varphi.
\end{eqnarray*}
Note that 
\begin{eqnarray*}
 |\mathscr L\mathrm d d|
       +\left|\mathscr L\left(\frac{V_\nu}{V}\right)\right|
       +|\mathscr L\varphi|
 \le C\mathcal F.
\end{eqnarray*}
Together with the $C^1$ bounds for $u$ and $f^{1/k}$, this gives
\begin{eqnarray}\label{eq:normal-error-second}
 \left|\mathscr Lw+2\sum_iF_i\lambda_i d_{ii}\right|
 \le C(1+\mathcal F).
\end{eqnarray}
In particular, the errors in \eqref{eq:normal-error-first} and
\eqref{eq:normal-error-second} are independent of $M$.

For parameters $\beta\ge1$, $\mu>0$, and $A_0\ge1$, we define
\begin{eqnarray}\label{eq:positive-barrier-definition}
 P_+\coloneqq(1+\beta d)w
       +\left(A_0+\frac M2\right)h_0
       \qquad\text{in }\bar\Omega_\mu.
\end{eqnarray}
We shall choose the constants $\beta$, then $\mu$, and finally $A_0$ so that
\begin{eqnarray}\label{eq:positive-barrier-claim}
 P_+\le0\qquad\text{in }\bar\Omega_\mu.
\end{eqnarray}
To begin with, we impose
\begin{eqnarray*}
 0<\mu\le\min\left\{\delta_0,\frac1{10\beta}\right\}.
\end{eqnarray*}

On $\Sigma$, by \eqref{eqn:w=0-Sigma}, we have 
\begin{eqnarray}\label{eqn:P-bdry-1}
P_+=0.
\end{eqnarray}  
On $\partial\Omega_\mu\setminus \Sigma$,
we have  $d=\mu$, together with \eqref{eqn:w=0-Sigma} imply 
\begin{eqnarray}\label{eqn:P-bdry-2}
 P_+\le(1+\beta\mu)C
       -\left(A_0+\frac M2\right)\mu(1-\mu)
 \le\frac{11}{10}C-\frac9{10}A_0\mu\le0,
\end{eqnarray}
provided $A_0\ge2C/\mu$. Under this choice of $A_0$, we conclude that any positive maximum
of $P_+$ must lie in the interior of $\Omega_\mu$.

In view of \eqref{eq:normal-error-first}, we derive 
\begin{eqnarray*}
 \mathscr LP_+
 &=&-2\beta\sum_iF_i\lambda_i d_i^2
       -2(1+\beta d)\sum_iF_i\lambda_i d_{ii}
       +\left(A_0+\frac M2\right)\mathscr Lh_0\\
 &&{}+(1+\beta d)
       \left(\mathscr Lw+2\sum_iF_i\lambda_i d_{ii}\right)
       +\beta w\mathscr Ld+2\beta\sum_iF_i d_iE_i.
\end{eqnarray*}
Since $1+\beta d\le11/10$, $|\mathscr Ld|\le C\mathcal F$, and
$|d_i|\le1$, \eqref{eq:normal-error-second} implies
\begin{eqnarray*}
 &&\left|(1+\beta d)
       \left(\mathscr Lw+2\sum_iF_i\lambda_i d_{ii}\right)
       +\beta w\mathscr Ld+2\beta\sum_iF_i d_iE_i\right|\\
 &&\qquad\le C(1+\beta)(1+\mathcal F)
       \le C_\beta(1+\mathcal F).
\end{eqnarray*}
Consequently,
\begin{eqnarray}\label{eq:paired-barrier-operator}
 \mathscr LP_+=-2\beta\sum_iF_i\lambda_i d_i^2
       -2(1+\beta d)\sum_iF_i\lambda_i d_{ii}+\left(A_0+\frac M2\right)\mathscr Lh_0
       +O_\beta(1+\mathcal F),
\end{eqnarray}
where the absolute value of the last term is bounded by
$C_\beta(1+\mathcal F)$, independently of $\mu$, $A_0$, and $M$.

Suppose that $P_+$ has a positive interior maximum, and evaluate the
following calculations at that point in an eigenframe of $A$.
Differentiating \eqref{eq:positive-barrier-definition}, we obtain
\begin{eqnarray*}
 0=(P_+)_i=(1+\beta d)(-\lambda_i d_i+E_i)
       +\beta w d_i
       +\left(A_0+\frac M2\right)(-1+2d)d_i.
\end{eqnarray*}
Thus, whenever $d_i\ne0$,
\begin{eqnarray}\label{eq:normal-eigenvalue-equation}
 \lambda_i
=\frac{\beta w}{1+\beta d}
       +\frac{-1+2d}{1+\beta d}\left(A_0+\frac M2\right)
       +\frac{E_i}{d_i}.
\end{eqnarray}
Denote 
\begin{eqnarray*}
 I\coloneqq\left\{1\leq i \leq n :\beta d_i^2\ge\frac1n\right\}.
\end{eqnarray*}
Since $\sum_i d_i^2=1$, there exists an index $i_0$ such that
\begin{eqnarray*}
 d_{i_0}^2\ge\frac1n,\qquad
 \beta d_{i_0}^2\ge\frac\beta n\ge\frac1n.
\end{eqnarray*}
In particular, $i_0\in I$. For every $i\in I$,
\begin{eqnarray*}
 \left|\frac{E_i}{d_i}\right|\le C\sqrt\beta,
       \qquad
 \frac{1-2d}{1+\beta d}\ge\frac8{11}\ge\frac23.
\end{eqnarray*}
It follows from \eqref{eq:normal-eigenvalue-equation} that
\begin{eqnarray*}
 \lambda_i
 \le C(\beta+\sqrt\beta)
       -\frac23\left(A_0+\frac M2\right).
\end{eqnarray*}
Choosing $A_0\ge6C(\beta+\sqrt\beta+1)$ therefore gives
\begin{eqnarray}\label{eq:normal-eigenvalue-signs}
 \lambda_i\le-\frac12\left(A_0+\frac M2\right)
, \qquad  \forall i\in I.
\end{eqnarray}
On the other hand, the assumed global reduction and $A_0\ge1$ imply
\begin{eqnarray}\label{eq:normal-eigenvalue-global}
 |\lambda_i|
 \le C_0(1+M)
       \le2C_0\left(A_0+\frac M2\right)
       =C_*\left(A_0+\frac M2\right), \qquad  \forall i,
\end{eqnarray}
where $C_*=2C_0$ is independent of all the barrier parameters.

Assume first that $1\le k<n$. For $i_0\in I$, 
\eqref{eq:normal-eigenvalue-signs} implies $\lambda_{i_0}<0$.
Lemma~\ref{lem:spectral} gives 
\begin{eqnarray*} 
 F_{i_0}\ge\frac{\mathcal F}{n-k+1}.
\end{eqnarray*}
Together with $d_{i_0}^2\ge1/n$ and
\eqref{eq:normal-eigenvalue-signs}, we find
\begin{eqnarray*}
 -2\beta F_{i_0}\lambda_{i_0}d_{i_0}^2
 \ge\beta F_{i_0}\left(A_0+\frac M2\right)d_{i_0}^2
 \ge\frac{\beta}{n(n-k+1)}
       \left(A_0+\frac M2\right)\mathcal F.
\end{eqnarray*}
The remaining terms with indices in $I$ are nonnegative, by using \eqref{eq:normal-eigenvalue-signs}. For $i\notin I$,  we have $\beta d_i^2<1/n$ and
\eqref{eq:normal-eigenvalue-global} give
\begin{eqnarray*}
 -2\beta\sum_{i\notin I}F_i\lambda_i d_i^2
 \ge-\frac2n\sum_{i\notin I}F_i|\lambda_i|\ge-\frac{2C_*}{n}
       \left(A_0+\frac M2\right)\mathcal F.
\end{eqnarray*}
Hence the first term on the right hand side of \eqref{eq:paired-barrier-operator} satisfies
\begin{eqnarray}\label{eq:positive-barrier-first-sum}
 -2\beta\sum_iF_i\lambda_i d_i^2
 \ge\left(\frac{\beta}{n(n-k+1)}-\frac{2C_*}{n}\right)
       \left(A_0+\frac M2\right)\mathcal F.
\end{eqnarray}
For the second term, \eqref{eq:normal-eigenvalue-global} yields 
\begin{eqnarray}\label{eq:positive-barrier-second-sum}
 -2(1+\beta d)\sum_iF_i\lambda_i d_{ii}
 \ge -\frac{11}{5}\sum_iF_i|\lambda_i| |d_{ii}|
 \ge-\frac{11}{5}C_*\|\nabla^2d\|_{C^0(\bar\Omega_{\delta_0})}
       \left(A_0+\frac M2\right)\mathcal F.
\end{eqnarray}
Inserting \eqref{eq:positive-barrier-first-sum},
\eqref{eq:positive-barrier-second-sum}, and
\eqref{eq:collar-convex-barrier} into
\eqref{eq:paired-barrier-operator}, we obtain
\begin{eqnarray}\label{eq:positive-barrier-coercive}
 \mathscr LP_+
 \ge\left(\frac\beta{n(n-k+1)}-C_2\right)
       \left(A_0+\frac M2\right)\mathcal F+\left[c_h\left(A_0+\frac M2\right)-C_\beta\right]
       (1+\mathcal F),
\end{eqnarray}
where $C_2$ is a constant, independently of $\beta$, $\mu$, $A_0$, and $M$.

We now fix the parameters by choosing $\beta$,
then $\mu$, and finally $A_0$ so that
\begin{eqnarray*}
 \beta&\ge&\max\{1,n^2(C_2+1)\},
 \nonumber\\
 0<\mu&\le&\min\left\{\delta_0,\frac1{10\beta}\right\},
 \nonumber\\
 A_0&\ge&\max\left\{1,\frac{2C}{\mu},
       6C(\beta+\sqrt\beta+1),\frac{2C_\beta}{c_h}\right\}.
\end{eqnarray*}
These choices are independent of $M$ and satisfy all the preceding
restrictions. Since $n(n-k+1)\le n^2$, it follow from
\eqref{eq:positive-barrier-coercive} that, at the positive maximum interior point of $P_+$
\begin{eqnarray*}
 0\ge\mathscr LP_+
 \ge\left(A_0+\frac M2\right)\mathcal F
       +\frac{c_h}{2}\left(A_0+\frac M2\right)(1+\mathcal F)
 >0,
\end{eqnarray*}
which is a contradiction. If $k=n$, a positive interior maximum is
already impossible by \eqref{eq:normal-eigenvalue-signs}, because
$A\in\Gamma_n$ requires all its eigenvalues to be positive.
 Together with \eqref{eqn:P-bdry-1} and \eqref{eqn:P-bdry-2}, this proves
\eqref{eq:positive-barrier-claim}.

Suppose $A_{\nu\nu}(x_0)=\inf_{\Sigma} A[u](\nu,\nu)$. 
On $\Sigma$, we have $(h_0)_\nu=1$ and $w_\nu=A_{\nu\nu}+O(1)$.
Thus \eqref{eqn:P-bdry-1} and \eqref{eq:positive-barrier-claim} imply
\begin{eqnarray*}
 0\le(P_+)_\nu (x_0)=A_{\nu\nu}+\left(A_0+\frac M2\right)+O(1), 
\end{eqnarray*}it yields 
\begin{eqnarray}\label{eqn:-C2-double-normal-barrier}
    \inf_{\Sigma} A[u](\nu,\nu) \geq -C-\frac M 2.
\end{eqnarray}

Similarly, we define the function
\begin{eqnarray}\label{eqn:defn-P-barrier-function}
  P_-\coloneqq (1+\beta d)w - \left(A_0+\frac M2\right) h_0
 \quad\text{in }\Omega_\mu.
\end{eqnarray}
On $\Sigma$, \eqref{eqn:w=0-Sigma} gives
\begin{eqnarray}\label{eqn:P-bdry-3}
 P_-=0.
\end{eqnarray}
We choose the parameters for $P_-$ independently of those for
$P_+$, still in the order $\beta$, $\mu$, and $A_0$.
Assume
\begin{eqnarray*}
 \beta\ge1,\qquad
 0<\mu\le\min\left\{\delta_0,\frac1{10\beta}\right\},
 \qquad A_0\ge\frac{2C}{\mu}.
\end{eqnarray*}
Then $P_-\ge0$ on $\partial\Omega_\mu\setminus \Sigma$, since
\begin{eqnarray*}
 P_-\ge-\frac{11}{10}C+
       \frac9{10}\left(A_0+\frac M2\right)\mu\ge0.
\end{eqnarray*}
Suppose that $P_-$ has a negative interior minimum.
At that point choose an orthonormal eigenframe of $A$.
The first derivative equation gives, whenever $d_i\ne0$,
\begin{eqnarray*}
 \lambda_i=
 \frac{\beta w}{1+\beta d}
 +\frac{1-2d}{1+\beta d}\left(A_0+\frac M2\right)
 +\frac{E_i}{d_i}.
\end{eqnarray*}
With $I=\{i:\beta d_i^2\ge1/n\}$ as above, increasing
$A_0$ so that $A_0\ge6C(\beta+\sqrt\beta+1)$ gives
\begin{eqnarray*}
 \lambda_i\ge\frac12\left(A_0+\frac M2\right)
 \qquad(i\in I).
\end{eqnarray*}
Choose $i_0$ with $d_{i_0}^2\ge1/n$; then $i_0\in I$.
The global bound \eqref{eq:normal-eigenvalue-global} still gives
\begin{eqnarray*}
 |\lambda_i|\le C_*\left(A_0+\frac M2\right).
\end{eqnarray*}
The same differentiation as in
\eqref{eq:paired-barrier-operator} yields
\begin{eqnarray*}
 \mathscr LP_-=-2\beta\sum_iF_i\lambda_i d_i^2
 -2(1+\beta d)\sum_iF_i\lambda_i d_{ii}
-\left(A_0+\frac M2\right)\mathscr Lh_0
 +O_\beta(1+\mathcal F).
\end{eqnarray*}

Fix $\varepsilon>0$ sufficiently small, independently of
$\beta$, $\mu$, and $A_0$.
First suppose
$\lambda_{\min}\ge-\varepsilon(A_0+M/2)$.
Every negative eigenvalue has index outside $I$, and hence
\begin{eqnarray*}
 -2\beta\sum_iF_i\lambda_i d_i^2
 \le\frac2n\sum_{\lambda_i<0}F_i|\lambda_i|
 \le\frac{2\varepsilon}{n}
       \left(A_0+\frac M2\right)\mathcal F.
\end{eqnarray*}
Euler's identity gives
\begin{eqnarray*}
 \sum_{\lambda_i>0}F_i\lambda_i
 =f^{1/k}+\sum_{\lambda_i<0}F_i|\lambda_i|
 \le C+\varepsilon\left(A_0+\frac M2\right)\mathcal F.
\end{eqnarray*}
Since $d_{ii}\le0$ in the fixed collar, it follows that
\begin{eqnarray*}
 -2(1+\beta d)\sum_iF_i\lambda_i d_{ii}
 \le C+C\varepsilon
       \left(A_0+\frac M2\right)\mathcal F,
\end{eqnarray*}
where $C$ is independent of $\beta$.
Choose $\varepsilon$ so that the sum of the two
$\varepsilon$-coefficients is at most $c_h/2$.
Using \eqref{eq:collar-convex-barrier}, and increasing
$A_0$ so that $A_0\ge4C_\beta/c_h$ after $\beta$ is fixed, we obtain
\begin{eqnarray*}
 \mathscr LP_-
 \le-\frac{c_h}{2}\left(A_0+\frac M2\right)
       (1+\mathcal F)+C_\beta(1+\mathcal F)<0.
\end{eqnarray*}
This contradicts the minimum property.

It remains to consider
$\lambda_{\min}<-\varepsilon(A_0+M/2)$.
The preceding bounds imply
\begin{eqnarray*}
 \lambda_{i_0}\ge\frac{\lambda_{\max}}{2C_*},
 \qquad
 -\lambda_{\min}\ge\frac{\varepsilon}{C_*}\lambda_{i_0}.
\end{eqnarray*}
Lemma~\ref{lem:spectral} (5), with smaller fixed parameters
if necessary, gives $F_{i_0}\ge c_*\mathcal F$, where
$c_*>0$ is independent of the barrier parameters.
Consequently,
\begin{eqnarray*}
 -2\beta\sum_iF_i\lambda_i d_i^2
 \le -\frac{\beta c_*}{n}
       \left(A_0+\frac M2\right)\mathcal F
       +\frac{2C_*}{n}\left(A_0+\frac M2\right)\mathcal F,\end{eqnarray*}
       and 
       \begin{eqnarray*}
 -2(1+\beta d)\sum_iF_i\lambda_i d_{ii}
 \le C\left(A_0+\frac M2\right)\mathcal F.
\end{eqnarray*}
Here the first inequality uses the nonpositive terms with
indices in $I$ and $\beta d_i^2<1/n$ for $i\notin I$.
Thus
\begin{eqnarray*}
 \mathscr LP_-
 \le\left(-\frac{\beta c_*}{n}+C\right)
       \left(A_0+\frac M2\right)\mathcal F
 -c_h\left(A_0+\frac M2\right)(1+\mathcal F)
 +C_\beta(1+\mathcal F).
\end{eqnarray*}
Choose $\beta$ sufficiently large that
$\beta c_*/n\ge C+1$, then choose $\mu$ as above, and
finally choose $A_0$ to meet the preceding restrictions and
$A_0\ge4C_\beta/c_h$.
The right hand side is strictly negative, again a
contradiction. When $k=n$, all eigenvalues are positive,
so only the first case is needed.
We have proved
\begin{eqnarray}\label{eq:negative-barrier-claim}
 P_-\ge0\qquad\text{in }\bar\Omega_\mu.
\end{eqnarray}
Suppose $A_{\nu\nu}(y_0)=\sup_{\Sigma} A[u](\nu,\nu)$. 
On $\Sigma$, we have $(h_0)_\nu=1$ and $w_\nu=A_{\nu\nu}+O(1)$.
Thus \eqref{eqn:P-bdry-3} and \eqref{eq:negative-barrier-claim} imply
\begin{eqnarray*}
 0\ge(P_-)_\nu (y_0)=A_{\nu\nu}-\left(A_0+\frac M2\right)+O(1), 
\end{eqnarray*}which follows
\begin{eqnarray*}
    \sup_{\Sigma} A[u](\nu,\nu) \leq C+\frac M 2.
\end{eqnarray*}
Together with \eqref{eqn:-C2-double-normal-barrier}, we derive 
\begin{eqnarray*}
    M\leq C. 
\end{eqnarray*}
The given $C^1$ bound in Proposition~\ref{prop:c1}  and $\nabla^2u=A[u]+ug$ now give
\eqref{eqn:c2-global-estimate}, after restoring the normalization.
This completes the proof.
\end{proof}

\subsection{Proof of Theorem~\ref{thm:existence}}

The preceding estimates allow us to solve \eqref{eq:main-problem}
by using the continuity method.
\begin{proof}[Proof of Theorem~\ref{thm:existence}]
We first assume that the domain and the data are smooth. For
$0\leq t\leq 1$, consider the family of equations
\begin{eqnarray}\label{eq:homotopy}
 \begin{cases}
  \sigma_k(A[u_t])=(1-t)\binom nk+tf \eqqcolon f_t &\text{in }\Omega,\\
  B_V(u_t)=\lambda_t+ (1-t)\dfrac{V_\nu}{V}+t\varphi
       &\text{on }\Sigma,
 \end{cases}
\end{eqnarray}
with the normalization
\begin{eqnarray}\label{eq:gauge}
 \int_\Omega\frac {u_t} V d\mu_g=0.
\end{eqnarray}
Denote  
\begin{eqnarray*}
\mathcal I\coloneqq
\left\{t\in[0,1]:
\begin{array}{l}
\text{there exists }(u_t,\lambda_t)
 \in C^2(\bar\Omega)\times\mathbb R\\
\text{satisfying \eqref{eq:homotopy} and \eqref{eq:gauge},
 with }A[u_t]\in\Gamma_k\text{ on }\bar\Omega 
\end{array}
\right\}.
\end{eqnarray*}
At $t=0$, a simple solution is
\begin{eqnarray*}
 u_0=-1+ \frac V {|\Omega|}\int_\Omega V^{-1} d\mu_g,  
 \qquad \lambda_0=0.
\end{eqnarray*} Hence $0\in \mathcal I$. 

From 
Propositions~\ref{prop:c1} and \ref{prop:c2}, we have 
\begin{eqnarray*}
    |\lambda_t|+\|u_t\|_{C^2(\bar \Omega)}\leq C.
\end{eqnarray*} 
Thus Eq. \eqref{eq:homotopy} is uniformly elliptic. The boundary operator is linear and
uniformly oblique, so the boundary H\"older estimate of
Lieberman--Trudinger
\cite[Theorem~6.2 and the supplementary remarks on p.~544]{LT} gives
\begin{eqnarray}\label{eq:c2-holder}
|\lambda_t|+ \|u_t\|_{C^{2,\alpha_1}(\bar\Omega)}\le C, 
\end{eqnarray}
for some fixed $\alpha_1\in(0,1)$.  For $0<\alpha<\alpha_1$, the preceding bounds \eqref{eq:c2-holder} 
give compactness in $C^{2,\alpha}(\bar\Omega)\times\RR$. 
Since $f_t$ is uniformly positive, limits remain admissible.
For smooth data, elliptic and oblique regularity then give smooth
solutions, to which the preceding differential estimates apply.
This proves closedness of $\mathcal I$.

For openness of $\mathcal I $, write $u=Vv$ and linearize
$F(A[Vv])=f_t^{1/k}$ together with the boundary condition
and its normalization \eqref{eq:gauge}. The resulting operator has Fredholm
index zero, and its kernel consists of pairs $(z,\mu)$ satisfying
\begin{eqnarray*}
 VF^{ij}\nabla_{ij}z+2F^{ij}V_i z_j&=&0
 \qquad\text{in }\Omega,\\
 Vz_\nu&=&\mu
 \qquad\text{on }\Sigma,\qquad
 \int_\Omega z d\mu_g=0.
\end{eqnarray*}
The maximum principle and Hopf lemma force $z$ to be constant:
otherwise its maximum and minimum would give opposite signs for $\mu$. Thus $z=\mu=0$. The Fredholm alternative and implicit function theorem give openness. Hence $\mathcal I =[0,1]$. Smooth approximation and oblique Schauder estimates yield $u\in C^{3,\alpha}(\bar\Omega)$ for the stated domain and data, after decreasing $\alpha$ if necessary.  

For the arbitrary-pole extension, the gradient estimate
uses only $h>0$, and Lemma~\ref{lem:reduction} uses
$h>0$ and $h-(V_\nu/V)g_\Sigma>0$ separately.
The double normal derivative estimate and the continuity argument
therefore apply without change.

Finally, let $(u_1,\lambda_1)$ and $(u_2,\lambda_2)$ be solutions
for the same data, and put $z=(u_1-u_2)/V$. Linearizing $F$ along
the segment joining the two admissible matrices gives
\begin{eqnarray}\label{eqn:two-solution-unique}
 \overline F^{ij}\nabla_{ij}z
 +2\overline F^{ij}(\log V)_i z_j=0,
 \qquad Vz_\nu |_{\Sigma}=\lambda_1-\lambda_2,
\end{eqnarray}
with uniformly positive definite coefficients. The maximum principle
and the boundary point lemma, as above, give $z\equiv c$ and
$\lambda_1=\lambda_2$. This completes the proof.
\end{proof}

We next show that assuming only strict geodesic convexity of $\Omega$ is insufficient for the existence of Theorem \ref{thm:existence}, even for $k=n$.
\begin{prop}
\label{prop:strict-convexity-counterexample}
Let $n\ge3$ and fix $p\in\HH^n$. There exists a connected bounded
domain $\Omega\Subset\HH^n$ containing $p$, with $C^\infty$
boundary $\Sigma=\partial\Omega$ and strictly geodesic convex, and 
\begin{eqnarray*}
0< f\in C^\infty(\bar\Omega),\qquad
 \varphi\in C^\infty(\Sigma),
\end{eqnarray*}
such that  
\begin{eqnarray}\label{eqn:monge-ampere-counterexample}
 \begin{cases}
  \det(\nabla^2u-ug)=f
       &  \qquad \text{in }\Omega,\\
  u_\nu=\dfrac{V_\nu}{V}u+\lambda+\varphi
       & \qquad \text{on }\Sigma, 
 \end{cases}
\end{eqnarray}
has no $n$-admissible solution $u\in C^2(\bar\Omega)$ for any
$\lambda\in\RR$. Moreover, $\Omega$ is not statically convex
relative to $V=\cosh d(p,\cdot)$.
\end{prop}

\begin{proof}
The construction has three ingredients: a convex function
singular along an interior axis, a boundary chosen to make its
oblique derivative smooth, and a comparison argument excluding
every classical solution.

\smallskip
\noindent\emph{Step 1. Reduction in Klein coordinates.}
The Klein model identifies $\HH^n$ with the Euclidean unit ball
$B_1\subset\RR^n$, equipped with the metric
\begin{eqnarray*}
 g_{ij}=\frac{\delta_{ij}}{s}+\frac{y_i y_j}{s^2},
 \qquad s=1-|y|^2.
\end{eqnarray*}
The direct computation gives 
\begin{eqnarray*}
    g^{ij}=s(\delta_{ij}-y_iy_j), ~\qquad ~\det (g_{ij})=s^{-(n+1)}, 
\end{eqnarray*} and the Christoffel symbol
\begin{eqnarray*}
\Gamma^k_{ij}
=\frac12 g^{k l}
\bigl(\partial_i g_{j l}
+\partial_jg_{i l}-\partial_l g_{ij}\bigr)
=\frac{y_i\delta_j^k+y_j\delta_i^k}{s}.
\end{eqnarray*}
Its unparametrized geodesics are Euclidean straight segments.
Choose the coordinates so that the pole $p$ corresponds to
$y=0$. Then
\begin{eqnarray*}
 V=\cosh d(p,\cdot)=s^{-1/2}.
\end{eqnarray*}
It follows
\begin{eqnarray*}
    V_i=\frac{Vy_i}{s},
\qquad
V_{ij}
=V\left(\frac{\delta_{ij}}s+\frac{3y_iy_j}{s^2}\right),
\end{eqnarray*}
and 
\begin{eqnarray*}
    \nabla_{ij}V
=V_{ij}-\Gamma^k_{ij}V_k
=V\left(\frac{\delta_{ij}}s+\frac{y_iy_j}{s^2}\right)
=Vg_{ij}.
\end{eqnarray*}
For $u=Vv$, a direct computation gives
\begin{eqnarray*}
    \nabla_{ij}(Vv)-Vvg_{ij}
=Vv_{ij}
+V_iv_j+V_jv_i
-V\Gamma^k_{ij}v_k+v\bigl(\nabla_{ij}V-Vg_{ij}\bigr).
\end{eqnarray*}
 It follows 
 \begin{eqnarray*}
A[u]_{ij}=Vv_{ij},\qquad
\sigma_n(A[u])
=\det\bigl(g^{-1}(VD^2v)\bigr)
=s^{(n+2)/2}\det D^2v,
\end{eqnarray*}
where $D$ denotes Euclidean differentiation and the first
equality is an equality of covariant components in these
coordinates; $\sigma_n$ is computed with respect to $g$.
If $\eta$ is the Euclidean outward unit normal, then the hyperbolic outward unit normal is 
\begin{eqnarray}\label{eqn:defn-nu-hyperbolic}
    \nu =\frac{g^{-1}\eta}
{\sqrt{g^{-1}(\eta,\eta)}}=\frac{s(I-y\otimes y)\eta}
{\sqrt{s(1-(y\cdot \eta)^2)}}=\frac{\sqrt s (I-y\otimes y)\eta}
{\sqrt{1-(y\cdot \eta)^2}}.
\end{eqnarray} Then
\begin{eqnarray*}
u_\nu-\frac{V_\nu}{V}u
=(Vv)_\nu-V_\nu v=Vv_\nu=\frac{(I-y\otimes y)\eta}
{\sqrt{1-(y\cdot \eta)^2}}\cdot Dv\coloneqq \beta\cdot Dv, 
\end{eqnarray*}where 

\begin{eqnarray*}
 \beta\coloneqq \frac{(I-y\otimes y)\eta}
 {\sqrt{1-(y\cdot\eta)^2}}.
\end{eqnarray*}
In particular, $n$-admissibility is equivalent to $D^2v>0$, and
\begin{eqnarray*}
 \beta\cdot\eta=\sqrt{1-(y\cdot\eta)^2}>0.
\end{eqnarray*}
Thus the problem \eqref{eqn:monge-ampere-counterexample} becomes a Euclidean Monge--Amp\`ere equation
with an oblique derivative boundary condition.

\smallskip
\noindent\emph{Step 2. A singular comparison function with a smooth
positive determinant.} We use the axis-singular ansatz underlying Pogorelov's
counterexample \cite{Pogorelov1971}, with the profile chosen below
to permit smooth oblique boundary data. 
Write $y=(x,z)\in\RR^{n-1}\times\RR$ and $r=|x|$. Define 
\begin{eqnarray}\label{eqn:v-gamma}
v(x,z)\coloneqq |x|^\gamma H(z), \qquad  \gamma\coloneqq 2-\frac2n, 
\end{eqnarray}
where the scalar function
\begin{eqnarray*}
 H(z)\coloneqq (1-z^2)^{-\gamma/2}e^{-\gamma z^2/8}.
\end{eqnarray*}
Since $1<\gamma<2$, the function $v$ is $C^1$ but not $C^2$ on
the axis $\{x=0\}$; Indeed, $v(0,z)=0, D v(0,z)=0$, and for any horizontal unit vector $e\in \SS^{n-2}$, there holds
\begin{eqnarray*}
    \frac{v(te,z)+v(-te,z)-2v(0,z)}{t^2}
=2H(z) |t|^{\gamma-2} \to +\infty, \quad \text{ as } t \to 0. 
\end{eqnarray*}
Denote
\begin{eqnarray*}
 Y\coloneqq \frac{H'}{H}
   =\frac{\gamma z}{1-z^2}-\frac{\gamma z}{4}.
\end{eqnarray*}
For $|z|\le a\coloneqq \frac 1 {10}$,
\begin{eqnarray*}
 Y'\ge\frac{3\gamma}{4},\qquad
 |Y|\le\gamma|z|,
\end{eqnarray*}
and hence
\begin{eqnarray*}
 (\gamma-1)Y'-Y^2
 \ge\frac{11}{16}\gamma(\gamma-1)>0.
\end{eqnarray*}
The transverse tangential eigenvalues of $D^2v$ are
$\gamma H r^{\gamma-2}$, with multiplicity $n-2$, and its
radial--axial block is
\begin{eqnarray*}
 \begin{pmatrix}
 \gamma(\gamma-1)H r^{\gamma-2}
   & \gamma H'r^{\gamma-1}\\
 \gamma H'r^{\gamma-1}
   & H''r^\gamma
 \end{pmatrix}.
\end{eqnarray*}
Consequently, $D^2v>0$ for $r>0$, and
\begin{eqnarray*}
 \det D^2v
 =\gamma^{n-1} r^{n\gamma-2n+2} H^n
   \bigl((\gamma-1)Y'-Y^2\bigr)=\gamma^{n-1}   H^n
   \bigl((\gamma-1)Y'-Y^2\bigr)
 =:\rho(z)>0.
\end{eqnarray*}
The choice $\gamma=2-2/n$ in \eqref{eqn:v-gamma} is precisely what cancels the radial
factor $r^{n\gamma-2n+2}$. Hence the determinant extends to a
smooth positive function of $z$, although $v$ itself is not
$C^2$ on the axis.

\smallskip
\noindent\emph{Step 3. Choosing the boundary to obtain smooth data $f$ and $\varphi$.}
Defining $\varphi=\beta\cdot Dv$ on an arbitrary boundary
would not ensure smoothness where the boundary meets the
axis. We instead arrange that $\beta\cdot Dv$ vanishes on
whole boundary neighborhoods of the two axial poles.

Let $Z(t)$ be the analytic solution near $t=0$ of
\begin{eqnarray*}
 2Z'(t)=
 \frac{(1-Z^2)Y(Z)-\gamma Z}
      {\gamma(1-t)-tZY(Z)},
 \qquad Z(0)=a.
\end{eqnarray*}
This equation is obtained by imposing $\beta\cdot Dv=0$ on
the upper cap $z=Z(r^2)$: an outward normal is proportional
to $(-2Z'(r^2)x,1)$. Since $H$ is even, the reflected lower
cap also has $\beta\cdot Dv=0$. Moreover,
\begin{eqnarray*}
 Z'(0)=-\frac{a(1-a^2)}8<0.
\end{eqnarray*}

These caps can be completed to a smooth, uniformly Euclidean
convex closed boundary, as follows. Define 
\begin{eqnarray*}
 G_{\rm loc}(t)\coloneqq a^2-Z(t)^2,\qquad
 g_0\coloneqq G_{\rm loc}'(0)=\frac{a^2(1-a^2)}4>0.
\end{eqnarray*}
Choose a smooth nonincreasing cutoff
$\chi:[0,\infty)\to[0,1]$, equal to one on $[0,1]$ and zero
on $[2,\infty)$. Choose $\delta>0$ sufficiently small that
$2\delta<1/4$ and $G_{\rm loc}$ is defined on $[0,2\delta]$.
Define, for $t\ge0$,
\begin{eqnarray*}
 \bar g(t)\coloneqq 
 \begin{cases}
  g_0+\chi(t/\delta)
       \bigl(G_{\rm loc}'(t)-g_0\bigr),
       & 0\le t<2\delta,\\
  g_0, & t\ge2\delta.
 \end{cases}
\end{eqnarray*}
Since $\bar g=g_0+O(\delta)$ and
$t\bar g'=O(\delta)$, we may ensure
\begin{eqnarray*}
 \frac{g_0}2\le\bar g\le2g_0,\qquad
 \bar g+2t\bar g'>0.
\end{eqnarray*}
Denote $\psi(t)\coloneqq 1-\chi(t/\delta)$ and
\begin{eqnarray*}
 C\coloneqq 
 \frac{a^2-\int_0^{1/4}\bar g(t) dt}
      {\int_0^{1/4}\psi(t) dt}>0,
 \qquad
 G(t)\coloneqq \int_0^t(\bar g+C\psi)(\tau) d\tau.
\end{eqnarray*}
Then $G=G_{\rm loc}$ near zero, $G(1/4)=a^2$, and
\begin{eqnarray*}
 G'>0,\qquad G'+2tG''>0.
\end{eqnarray*}
It follows that
\begin{eqnarray*}
 \Omega_K\coloneqq 
 \left\{(x,z)\in \RR^{n-1}\times \RR: ~G(|x|^2)+z^2<a^2 \right\}, \end{eqnarray*}
is smooth and uniformly Euclidean convex: the Hessian of
its defining function has eigenvalues
\begin{eqnarray*}
 2,\qquad 2G',\qquad 2(G'+2tG''),
\end{eqnarray*}
with the appropriate multiplicities. Also $|z|<a$ and
$|x|<1/2$, so $\Omega_K\Subset B_1$. The corresponding
hyperbolic domain is strictly geodesically convex and
contains the pole.

Now define
\begin{eqnarray*}
 f(y)=s^{(n+2)/2}\rho(z),\qquad
 \varphi=\beta\cdot Dv
 \quad\text{on }\partial\Omega_K.
\end{eqnarray*}
The function $f$ is smooth and strictly positive. The
boundary datum $\varphi$ is also smooth: it is identically
zero near the axial poles, and every other boundary point
lies off the singular axis.

The domain fails static convexity. Indeed, at the upper
axial pole $p_+=(0,a)$, using \eqref{eqn:defn-nu-hyperbolic},  direct computation gives
\begin{eqnarray*}
 \kappa_i=-2Z'(0)=\frac{a(1-a^2)}4>0,
 \qquad \frac{V_\nu}{V}=\nu(\log V)=\frac{\sqrt{s}(y\cdot \eta)}{\sqrt{1-(y\cdot \eta)^2}} \bigg|_{p_+} =a.
\end{eqnarray*}
Thus every hyperbolic principal curvature there is strictly
less than $V_\nu/V$, and $h-(V_\nu/V)g_\Sigma=-\frac{a(3+a^2)}4g_\Sigma$ is negative
definite at $p_+$.

\smallskip
\noindent\emph{Step 4. Excluding every classical admissible solution.} Suppose that \eqref{eqn:monge-ampere-counterexample} admits an
$n$-admissible solution $u\in C^2(\bar\Omega)$ for some
$\lambda\in\RR$. Let $\pi:\HH^n\to B_1$ denote the Klein
coordinate map used in Step~1. Thus
\begin{eqnarray*}
 \pi(\Omega)=\Omega_K,\qquad
 \pi(\Sigma)=\partial\Omega_K.
\end{eqnarray*}
Define $w\in C^2(\bar\Omega_K)$ by
\begin{eqnarray*}
 w(y)\coloneqq
 \frac{u(\pi^{-1}(y))}{V(\pi^{-1}(y))}
 =\sqrt{1-|y|^2} u(\pi^{-1}(y)),
 \qquad y\in\bar\Omega_K.
\end{eqnarray*}
As in Step~3, we use the same symbols $f$ and $\varphi$
for their expressions in Klein coordinates. The transformation
identities in Step~1 and the choice of data in Step~3 give
\begin{eqnarray}\label{eqn:w-equation-bdry}
 \begin{cases}
  D^2w>0
       & \qquad \text{in }\bar\Omega_K,\\
  \det D^2w=\rho(z)
       & \qquad \text{in }\Omega_K,\\
  \beta\cdot Dw=\lambda+\varphi
       & \qquad \text{on }\partial\Omega_K.
 \end{cases}
\end{eqnarray} 
The equation for $w$ holds throughout $\Omega_K$, including
the interior axis $\{(0,z):-a<z<a\}$. By contrast, the
comparison function $v=|x|^\gamma H(z)$ is only $C^1$ across
the axis and satisfies its classical Hessian equation off
the axis. Its boundary derivative is nevertheless defined
everywhere on $\partial\Omega_K$ and equals $\varphi$.

At $p_\pm=(0,\pm a)$,
\begin{eqnarray*}
 \beta(p_\pm)=\pm\sqrt{1-a^2} e_z,\qquad
 \varphi(p_\pm)=0.
\end{eqnarray*}
The boundary condition and strict convexity along the axis
therefore give
\begin{eqnarray*}
 \frac{2\lambda}{\sqrt{1-a^2}}
 =w_z(p_+)-w_z(p_-)
 =\int_{-a}^{a}w_{zz}(0,z) dz>0.
\end{eqnarray*}
Thus 
\begin{eqnarray}\label{eqn:lambda>0}
    \lambda>0.
\end{eqnarray}

On the other hand, fix $\varepsilon>0$
and consider
\begin{eqnarray*}
 d_\varepsilon\coloneqq v-w+\varepsilon|y|^2.
\end{eqnarray*}
Its maximum value on $\bar\Omega_K$ cannot occur in the interior.
In fact, at an interior maximum off the axis we would have
\begin{eqnarray*}
 D^2w\ge D^2v+2\varepsilon I,
\end{eqnarray*}
and hence
\begin{eqnarray*}
 \det D^2w>\det D^2v=\rho(z),
\end{eqnarray*}
a contradiction. At an interior maximum on the axis,
$D_xw=0$; for any unit transverse direction $e$,
\begin{eqnarray*}
 d_\varepsilon(te,z)-d_\varepsilon(0,z)
 =H(z)|t|^\gamma+O(t^2)>0
\end{eqnarray*}
for small nonzero $t$, again a contradiction.

Therefore, $d_\varepsilon$ attains its maximum on the boundary. Since $\beta\cdot\eta>0$, using \eqref{eqn:w-equation-bdry}, at this boundary maximum point
\begin{eqnarray*}
 0\le\beta\cdot Dd_\varepsilon
 =-\lambda+2\varepsilon\beta\cdot y.
\end{eqnarray*}
Consequently,
\begin{eqnarray*}
 \lambda\le
 2\varepsilon
 \sup_{\partial\Omega_K}|\beta\cdot y|.
\end{eqnarray*}
Letting $\varepsilon\downarrow0$ yields $\lambda\le0$,
contradicting \eqref{eqn:lambda>0}. Hence we complete the proof.
\end{proof}

We record the comparison principle needed below for solutions of different
Hessian orders in Eq. \eqref{eq:main-problem-lambda-j}.
\begin{lem}\label{prop:constant-dependence}
Let $\Omega$ be smooth and connected, and let
$u,\widetilde u\in C^2(\bar\Omega)$ satisfy
$A[u],A[\widetilde u]\in\Gamma_k$ on $\bar\Omega$. Suppose that
\begin{eqnarray*}
\begin{cases}
 F(A[\widetilde u])\ge F(A[u]) \qquad &\text{ in }\Omega,\\
B_V(\widetilde u-u)=\widetilde\lambda-\lambda \qquad &\text{ on }\Sigma, 
\end{cases}
\end{eqnarray*}
where $\lambda,\widetilde\lambda$ are constants. Then
$\widetilde\lambda\ge\lambda$. Moreover, equality holds if and only if
$\widetilde u-u$ is a multiple of $V$.
\end{lem}
\begin{proof}
For $z\coloneqq (\widetilde u-u)/V$, the linearization as in the uniqueness argument (cf. \eqref{eqn:two-solution-unique}) gives
\begin{eqnarray*}
 \overline F^{ij}\nabla_{ij}z
 +2\overline F^{ij}(\log V)_i z_j\ge0,
\text{ in } \Omega,  \qquad Vz_\nu |_{\Sigma}=\widetilde\lambda-\lambda.
\end{eqnarray*}
If $z$ is nonconstant, the strong maximum principle and Hopf's lemma
give $z_\nu>0$ at a boundary maximum, so
$\widetilde\lambda>\lambda$. If $z$ is constant, both the boundary
constant and the tensors $A[\widetilde u]$, $A[u]$ agree.
\end{proof}

\section{Geometric applications}
\label{sec:geometry}

We first give a rigidity lemma for the equality cases, which will be used later.
\begin{lem}\label{lem:geometric-rigidity}
Let $\Omega\Subset\HH^n$ have $C^2$ connected boundary $\Sigma$,
and assume $h-\frac{V_\nu}{V}g_\Sigma\ge0$. If $u\in C^2(\bar \Omega)$ satisfies 
\begin{eqnarray} \label{eqn:rigidity-geodesic ball}
\begin{cases}
A[u]=g \qquad &\text{ in }\Omega,\\
B_Vu=\lambda \qquad &\text{ on }\Sigma
\end{cases}
\end{eqnarray}for some constant $\lambda$, 
then  $\Omega$ is a
geodesic ball of radius $R$, with $\lambda=\tanh R$ and the vector field
$X$ defined in \eqref{eqn:defn-X} satisfies $X^T=0$ on $\Sigma$.  
\end{lem}

\begin{proof}
The argument below is motivated by the
rigidity argument in \cite[Section~6]{LiXia}. 
The integral identity  gives $\lambda>0$, and \eqref{eqn:mixed-bdry-A-nu-alpha}, \eqref{eqn:rigidity-geodesic ball} give
\begin{eqnarray}\label{eqn:mixed-bdry-A-nu-alpha-2}
\left(h-\frac{V_\nu}{V}g_\Sigma\right)(X^T,Y)
=-A[u](\nu,Y)=-g(\nu,Y)=0,
\end{eqnarray} for any $Y\in T\Sigma$. 
At a point where $\Sigma$ touches an enclosing geodesic ball of
radius $R_0$ from the inside, $h\ge\coth R_0 g_\Sigma$.
Since $V_\nu/V<1$, the tensor in parentheses \eqref{eqn:mixed-bdry-A-nu-alpha-2} is positive definite
near that point, and hence $X^T=0$ there.  Hence 
\begin{eqnarray}\label{eqn:X=lambda-nu}
    X=\< X,\nu\> \nu= (B_V u) \nu=\lambda \nu.
\end{eqnarray} Tangential restriction of
\eqref{eq:X-derivative} and \eqref{eqn:rigidity-geodesic ball} imply
\begin{eqnarray*}
    \n_i X_j = A_{ij}- X_i (\log V)_j=A_{ij}=g_{ij}, 
\end{eqnarray*}together with \eqref{eqn:X=lambda-nu}, we obtain $  h= \lambda^{-1} g_\Sigma$ and
$\lambda^{-1}\geq \coth R_0 >1$.

Let $U$ be the component of the interior of $\{X^T=0\}$ containing
this neighborhood. On $U$, $h=\lambda^{-1}g_\Sigma$, and the same
holds on $\overline U$ by continuity. The static-convexity tensor is
therefore positive definite near every point of $\overline U$.
The identity \eqref{eqn:mixed-bdry-A-nu-alpha-2} forces $X^T=0$ in these neighborhoods, so $U$
is closed as well as open in $\Sigma$. Connectedness gives
$U=\Sigma$. Thus $\Sigma$ is totally umbilical with constant
principal curvature $\lambda^{-1}>1$, and is a geodesic sphere.
Its radius satisfies $\lambda^{-1}=\coth R$.
\end{proof}

\subsection{Proof of Theorem~\ref{thm:hierarchy}}

We start by proving the last inequality involving the sharp upper bound of $\lambda_n$ in \eqref{eq:constant-hierarchy} and state it separately in the following conclusion.
 \begin{prop}\label{prop:MA-boundary-constant-upper}
Let $\Omega$ and $V$ satisfy the domain and pole assumptions of
Theorem~\ref{thm:existence}, and let $(u_n,\lambda_n)$ be the
$n$-admissible solution of
\begin{eqnarray*}
\begin{cases}
\det(\nabla^2u_n-u_ng)=1
    & \qquad \text{in }\Omega,\\
B_Vu_n=\lambda_n
    & \qquad \text{on }\Sigma.
\end{cases}
\end{eqnarray*}
Then
\begin{eqnarray}\label{eq:MA-boundary-constant-upper}
\lambda_n^n\int_\Sigma VH_{n-1} d\sigma_g
\le n\int_\Omega V d\mu_g.
\end{eqnarray}
Equality holds if and only if $\Omega$ is a geodesic ball,
whose center need not coincide with $p$.
\end{prop}

\begin{proof}
The identity \eqref{eq:flux} and the Newton--Maclaurin
inequality give
\begin{eqnarray*}
\lambda_n\int_\Sigma V d\sigma_g
=\int_\Omega V\tr A[u_n] d\mu_g
\ge n\int_\Omega V d\mu_g>0.
\end{eqnarray*}
In particular, $\lambda_n>0$.

We use the hyperboloid coordinates and the notation in the proof of
Lemma~\ref{lem:static-convex-domain}. Thus $p=(1,0)$,
\begin{eqnarray*}
\iota(x)=(\sqrt{1+|x|^2},x),\qquad
K=\iota^{-1}(\bar\Omega),\qquad
s=x\cdot\eta,\qquad
\ell=\sqrt{1+s^2},
\end{eqnarray*}
where $\eta$ is the outward Euclidean unit normal to $\partial K$.
Lemma~\ref{lem:static-convex-domain} implies that $K$ is a smooth
strictly convex body containing the origin in its interior,
$s>0$, and $h>0$ on $\Sigma$. Recall also that
\begin{eqnarray*}
V=\sqrt{1+|x|^2},\qquad
\nu=\frac{(sV,\eta+sx)}{\ell},\qquad
V_\nu=\frac{sV}{\ell}.
\end{eqnarray*}

Pass to the Klein coordinates
\begin{eqnarray*}
y=\frac{x}{V},\qquad
\Omega_K=
\left\{\frac{x}{\sqrt{1+|x|^2}}:
             x\in\operatorname{int}K\right\},
\qquad
v=\frac{u_n}{V}.
\end{eqnarray*}
Here $v$ is regarded as a function of $y$, and $D$ denotes
Euclidean differentiation. Since $\Omega$ is geodesically convex
and Klein geodesics are straight segments, $\Omega_K$ is convex.
The calculation in the proof of
Proposition~\ref{prop:strict-convexity-counterexample} gives
\begin{eqnarray*}
A[u_n]_{ij}=Vv_{ij},\qquad
\det(A[u_n])
=(1-|y|^2)^{(n+2)/2}\det D^2v.
\end{eqnarray*}
Consequently,
\begin{eqnarray*}
D^2v>0\quad\text{on }\bar\Omega_K,\qquad
\det D^2v=(1-|y|^2)^{-(n+2)/2}.
\end{eqnarray*}
Moreover,
\begin{eqnarray*}
V d\mu_g=(1-|y|^2)^{-(n+2)/2} dy.
\end{eqnarray*}
Strict convexity of $v$ on the convex domain $\Omega_K$ implies
that $Dv$ is injective on $\bar\Omega_K$. Its restriction to
$\Omega_K$ is therefore a diffeomorphism onto its image.
The change-of-variables formula yields
\begin{eqnarray*}
|Dv(\Omega_K)|
=\int_{\Omega_K}\det D^2v dy
=\int_\Omega V d\mu_g,
\end{eqnarray*}
where the absolute value of a subset of $\RR^n$ denotes its
Euclidean volume.

We next describe the convex body determined by the boundary
condition. Under the Euclidean Gauss parametrization of
$\partial K$, the function $s=x\cdot\eta$ is the support function
of $K$. The function $\ell=\sqrt{1+s^2}$ is the support function
of a smooth strictly convex body whose boundary is parametrized by
\begin{eqnarray*}
x\longmapsto\frac{\eta+sx}{\ell},
\qquad x\in\partial K.
\end{eqnarray*}
Its outward Euclidean unit normal at this point is $\eta$. Equivalently, this convex body is
\begin{eqnarray*}
\bigcap_{x\in\partial K}
\left\{z\in\RR^n:z\cdot\eta(x)\le\ell(x)\right\}.
\end{eqnarray*}
We compute its volume directly.

The map $x\mapsto(\eta+sx)/\ell$ is the spatial component of the
hyperbolic unit normal $\nu$. The ambient differential of $\nu$
in directions tangent to $\Sigma$ is given by the shape operator
of $\Sigma$. Since $h>0$, its determinant is
$\sigma_{n-1}(h)=H_{n-1}$.
The Gram matrix of the spatial components of a hyperbolic
orthonormal tangential frame is
$g_\Sigma+\mathrm d_\Sigma V\otimes\mathrm d_\Sigma V$.
Thus the Euclidean area Jacobian of the above parametrization is
\begin{eqnarray*}
H_{n-1}\sqrt{1+|\nabla_\Sigma V|^2}
=H_{n-1}\sqrt{V^2-V_\nu^2}
=\frac{V}{\ell}H_{n-1}.
\end{eqnarray*}
Here we used \eqref{eq:static} and $V_\nu=sV/\ell$.
Its support function is $\ell$, so the Euclidean divergence
theorem gives
\begin{eqnarray*}
&&n\left|
\bigcap_{x\in\partial K}
\left\{z\in\RR^n:z\cdot\eta(x)\le\ell(x)\right\}
\right|\\
&&\qquad=
\int_\Sigma \ell \frac{V}{\ell}H_{n-1} d\sigma_g
=\int_\Sigma VH_{n-1} d\sigma_g.
\end{eqnarray*}

In the Klein coordinates, the boundary condition takes the form
\begin{eqnarray*}
\beta\cdot Dv=\lambda_n
\quad\text{on }\partial\Omega_K,\qquad
\beta=\frac{\eta}{\ell}.
\end{eqnarray*}
Indeed, this follows by differentiating $y=x/V$ in the direction
$V\nu$. The vector $\eta$ here is still the outward normal to
$\partial K$. The outward Euclidean unit normal to
$\partial\Omega_K$ at $y=x/V$ is
\begin{eqnarray*}
\frac{\eta+sx}{|\eta+sx|}.
\end{eqnarray*}
In particular, the boundary vector is strictly oblique, since
\begin{eqnarray*}
\beta\cdot\frac{\eta+sx}{|\eta+sx|}
=\frac{1+s^2}{\ell|\eta+sx|}
=\frac{\ell}{|\eta+sx|}>0.
\end{eqnarray*}

We now apply the minimum-point argument underlying the ABP method; see \cite{Cabre2008}.
Fix any $z\in\RR^n$ such that
\begin{eqnarray*}
z\cdot\eta(x)<\lambda_n\ell(x)
\qquad\text{for every }x\in\partial K.
\end{eqnarray*}
The function $v(y)-z\cdot y$ attains its minimum on
$\bar\Omega_K$. This minimum cannot occur on
$\partial\Omega_K$. Indeed, at a boundary minimum its tangential
derivatives vanish and its outward normal derivative is
nonpositive. Strict obliqueness would then imply
\begin{eqnarray*}
\beta\cdot D(v-z\cdot y)\le0.
\end{eqnarray*}
On the other hand, the boundary condition gives
\begin{eqnarray*}
\beta\cdot D(v-z\cdot y)
=\lambda_n-\frac{z\cdot\eta}{\ell}>0,
\end{eqnarray*}
a contradiction. Hence the minimum is attained in $\Omega_K$,
where $Dv=z$. We have proved
\begin{eqnarray*}
\left\{z\in\RR^n:
z\cdot\eta(x)<\lambda_n\ell(x)
\text{ for every }x\in\partial K\right\}
\subset Dv(\Omega_K).
\end{eqnarray*}
The set on the left is the interior of the $\lambda_n$-dilate
of the convex body described above. Its boundary has zero
Euclidean volume. Consequently,
\begin{eqnarray*}
\frac{\lambda_n^n}{n}
\int_\Sigma VH_{n-1} d\sigma_g
&=&
\left|
\bigcap_{x\in\partial K}
\left\{z\in\RR^n:
z\cdot\eta(x)\le\lambda_n\ell(x)\right\}
\right|\\
&\le&|Dv(\Omega_K)|
=\int_\Omega V d\mu_g.
\end{eqnarray*}
This proves \eqref{eq:MA-boundary-constant-upper}.

Suppose that equality holds in
\eqref{eq:MA-boundary-constant-upper}.
The preceding inclusion and equality of volumes imply
\begin{eqnarray*}
Dv(\Omega_K)=
\left\{z\in\RR^n:
z\cdot\eta(x)<\lambda_n\ell(x)
\text{ for every }x\in\partial K\right\}.
\end{eqnarray*}
Indeed, $Dv(\Omega_K)$ is open, so any additional point would
give a subset of positive volume outside the convex body.
Since $Dv$ is continuous and injective on the compact set
$\bar\Omega_K$, it is a homeomorphism onto its image and maps
$\partial\Omega_K$ onto the boundary of the same dilated
convex body.

At corresponding boundary points $y=x/V$, the Neumann condition
gives
\begin{eqnarray*}
Dv(y)\cdot\eta(x)=\lambda_n\ell(x).
\end{eqnarray*}
Strict convexity makes the supporting point with outward
normal $\eta(x)$ unique. The boundary parametrization therefore
implies
\begin{eqnarray*}
Dv(y)=\lambda_n\frac{\eta+sx}{\ell}
\qquad\text{on }\partial\Omega_K.
\end{eqnarray*}
Along $\Sigma$, we have
\begin{eqnarray*}
\eta\cdot\mathrm dx=0,\qquad
x\cdot\mathrm dx=V \mathrm dV,\qquad
(\eta+sx)\cdot x=sV^2.
\end{eqnarray*}
It follows that
\begin{eqnarray*}
\mathrm dv
&=&Dv\cdot\mathrm dy\\
&=&\lambda_n\frac{\eta+sx}{\ell}\cdot
\left(\frac{\mathrm dx}{V}
      -\frac{x \mathrm dV}{V^2}\right)
=0
\qquad\text{on }\Sigma.
\end{eqnarray*}
Since $\Sigma$ is connected by
Lemma~\ref{lem:static-convex-domain}, $u_n/V$ is constant on
$\Sigma$. Subtracting this constant multiple of $V$ from $u_n$
preserves both $A[u_n]$ and $B_Vu_n$. With this representative, we have 
\begin{eqnarray*}
\begin{cases}
\det(\nabla^2u_n-u_ng)=1
    &\text{in }\Omega,\\
u_n=0,\qquad (u_n)_\nu=\lambda_n>0
    &\text{on }\Sigma.
\end{cases}
\end{eqnarray*}
By  applying  
\cite[Theorem~1.2]{GaoJiaYan2023}, it follows that $\Omega$ is a geodesic ball.

Conversely, let $\Omega$ be a geodesic ball of radius $R$. Then
\begin{eqnarray*}
A[u_n]=g,\qquad
\lambda_n=\tanh R,\qquad
H_{n-1}=\coth^{n-1}R.
\end{eqnarray*}
Using \eqref{eq:flux}, we obtain
\begin{eqnarray*}
\lambda_n^n\int_\Sigma VH_{n-1} d\sigma_g
=\tanh R\int_\Sigma V d\sigma_g
=n\int_\Omega V d\mu_g.
\end{eqnarray*}
This argument applies to every center of the ball permitted
by the hypotheses. Hence equality holds, and the proof is complete.
\end{proof}

Now we turn to prove Theorem \ref{thm:hierarchy}.
\begin{proof}[Proof of Theorem~\ref{thm:hierarchy}]
The compatibility identity gives
\begin{eqnarray*}
 \lambda_1=\frac{n\int_\Omega V d\mu_g}
                  {\int_\Sigma V d\sigma_g}>0.
\end{eqnarray*}
For $1\le k<\ell\le n$, Newton--Maclaurin inequality shows that $u_\ell$ is
$k$-admissible and
\begin{eqnarray*}
 \left(\frac{\sigma_k(A[u_\ell])}{\binom nk}\right)^{1/k}
 \ge\left(\frac{\sigma_\ell(A[u_\ell])}{\binom n\ell}\right)^{1/\ell}
 =1.
\end{eqnarray*}
The comparison principle in Lemma~\ref{prop:constant-dependence} therefore gives $$\lambda_\ell\ge\lambda_k.$$ 

The sharp upper bound of $\lambda_n$  in \eqref{eq:constant-hierarchy} follows from Proposition \ref{prop:MA-boundary-constant-upper}.  Next we show the last two inequalities in \eqref{eq:constant-hierarchy}. 
We first recall the Minkowski integral formula (see, e.g., \cite[Proposition 2.5]{GuanLi2015})
\begin{eqnarray}\label{eq:geometric-upper-minkowski-identities}
\int_\Sigma VH_{j-1} d\sigma_g
=
\int_\Sigma V_\nu H_j d\sigma_g,
\qquad 1\le j\le n-1.
\end{eqnarray} 
It follows from $H_j>0$ that 
\begin{eqnarray*}
\int_\Sigma VH_{j-1} d\sigma_g
=\int_\Sigma
\frac{V_\nu}{V}VH_j d\sigma_g \le
\left(\max_\Sigma\frac{V_\nu}{V}\right)
\int_\Sigma VH_j d\sigma_g,
\qquad 1\le j\le n-1.
\end{eqnarray*}
Eq. \eqref{eq:static} and the divergence theorem yield
\begin{eqnarray*}
n\int_\Omega V d\mu_g
=\int_\Sigma V_\nu d\sigma_g
\le
\left(\max_\Sigma\frac{V_\nu}{V}\right)
\int_\Sigma V d\sigma_g.
\end{eqnarray*}
Applying the preceding comparisons successively, we obtain
\begin{eqnarray*}
n\int_\Omega V d\mu_g
&\le &\left(\max_\Sigma\frac{V_\nu}{V}\right)
\int_\Sigma V d\sigma_g\le
\left(\max_\Sigma\frac{V_\nu}{V}\right)^2
\int_\Sigma VH_1 d\sigma_g \\& \le & \cdots\le
\left(\max_\Sigma\frac{V_\nu}{V}\right)^n
\int_\Sigma VH_{n-1} d\sigma_g.
\end{eqnarray*} The inequality $\frac{V_\nu}{V}< 1$ is obvious. Thus we prove the last two inequalities in \eqref{eq:constant-hierarchy}. 

If $\lambda_k=\lambda_\ell$, the equality case of the comparison
implies that $u_\ell-u_k$ is a multiple of $V$. Equality in
Newton--Maclaurin then gives $A[u_\ell]=g$, and
Lemma~\ref{lem:geometric-rigidity} shows that $\Omega$ is a ball.

Conversely, on a geodesic ball $B_R(q)$, define
\begin{eqnarray}\label{eq:ball-reference-solution}
 u_q(x)\coloneqq\frac{\cosh d(q,x)}{\cosh R}-1.
\end{eqnarray}
Then $A[u_q]=g$, $u_q=0$ on $\Sigma$, and
$B_Vu_q=(u_q)_\nu=\tanh R$ for every pole. Adding a multiple
of $V$ gives the required normalization. Uniqueness therefore
gives $\lambda_j=\tanh R$ for all $j$.
 
\end{proof}

\subsection{Proof of Theorem~\ref{thm:deficits}}

We prove Theorem~\ref{thm:deficits} by applying the Reilly formula to the normalized $2$-Hessian solution in \eqref{eqn:sigma-2-neumann}, which yields geometric deficit identities.

\begin{proof}[Proof of Theorem~\ref{thm:deficits}]
Newton--Maclaurin inequality in Lemma \ref{lem:spectral} gives $$\tr A[u]\ge  n\left({\sigma_2(A[u])} /{\binom{n}{2}} \right)^{1/2}= n.$$ Hence
\begin{eqnarray}\label{eq:geometry-D}
 \int_\Omega V(\tr A[u]-n) d\mu_g\ge0,\qquad
 \int_\Sigma V\left(h-\frac{V_\nu}{V}g_\Sigma\right)
     (X^T,X^T) d\sigma_g\ge0.
\end{eqnarray}
Since $X_\nu=\lambda$ on $\Sigma$, the identity gives
\begin{eqnarray}\label{eq:geometry-flux}
 \lambda\int_\Sigma V d\sigma_g
 =n\int_\Omega V d\mu_g+\int_\Omega V(\tr A[u]-n) d\mu_g.
\end{eqnarray}
Thus $\lambda>0$. Lemma~\ref{lem:reilly} gives
\begin{eqnarray}\label{eq:geometry-reilly}
 n(n-1)\int_\Omega V d\mu_g
 =(n-1)\lambda^2\int_\Sigma VH d\sigma_g
 +\int_\Sigma V\left(h-\frac{V_\nu}{V}g_\Sigma\right)
     (X^T,X^T) d\sigma_g.
\end{eqnarray}
Combining these two identities proves \eqref{eq:main-minkowski}.

Since $2\sigma_2(A[u])=(\tr A[u])^2-|A[u]|^2=n(n-1)$, then
$$|A[u]-g|^2=(\tr A[u]-n)(\tr A[u]+n-2),$$ therefore
\begin{eqnarray}\label{eq:geometry-bulk-energy}
 \int_\Omega V(\tr A[u]-n) d\mu_g
 =\int_\Omega\frac{V|A[u]-g|^2}{\tr A[u]+n-2} d\mu_g.
\end{eqnarray}
Together with
\begin{eqnarray}\label{eq:geometry-cauchy-deficit}
 \int_\Sigma VH\left(\frac1H-
     \frac{\int_\Sigma V d\sigma_g}{\int_\Sigma VH d\sigma_g}\right)^2 d\sigma_g
 =\int_\Sigma\frac VH d\sigma_g-
     \frac{\left(\int_\Sigma V d\sigma_g\right)^2}
          {\int_\Sigma VH d\sigma_g},
\end{eqnarray}
and substitute \eqref{eq:main-minkowski} to obtain
\eqref{eq:main-hk}.  

Equality in either inequality forces the integral term in
\eqref{eq:geometry-D} to vanish. Equation~\eqref{eq:geometry-bulk-energy}
gives $A=g$, and Lemma~\ref{lem:geometric-rigidity} implies that
$\Omega$ is a geodesic ball. On a ball, the reference solution has $A=g$ and $X^T=0$, so equality
holds. Finally, addition of a multiple of $V$ changes neither
$A$ nor $X$.

\end{proof}

For comparison with the linear potential argument of \cite{Xia,LiXia}, let
$v\in C^2(\bar \Omega)$ solve the linear Neumann problem: 
\begin{eqnarray*}
 \Delta v-nv=n\quad\hbox{in }\Omega,\qquad\qquad 
 B_Vv=\frac{n\int_\Omega V d\mu_g}{\int_\Sigma V d\sigma_g}
 \quad\hbox{on }\Sigma.
\end{eqnarray*}
This problem is solvable on every smooth connected bounded domain.
Since $\tr A[v]=n$, retaining both remainder terms in
\eqref{eq:reilly} gives
\begin{eqnarray}\label{eq:linear-deficit}
\begin{aligned}
 &\left(\int_\Sigma V d\sigma_g\right)^2
 -n\left(\int_\Omega V d\mu_g\right)
       \left(\int_\Sigma VH d\sigma_g\right)\\
 &=\frac{\left(\int_\Sigma V d\sigma_g\right)^2}
                  {n(n-1)\int_\Omega V d\mu_g}
 \biggl[\int_\Omega V|A[v]-g|^2 d\mu_g+\int_\Sigma V^3\left(h-\frac{V_\nu}{V}g_\Sigma\right)
\bigl(\nabla_\Sigma(v/V),\nabla_\Sigma(v/V)\bigr) d\sigma_g\biggr].
\end{aligned}
\end{eqnarray}
Static convexity makes both remainder terms on the right hand side nonnegative. Eq. 
\eqref{eq:geometry-cauchy-deficit} similarly gives the weighted
Heintze--Karcher inequality when $H>0$.

\subsection{Proof of Theorem~\ref{thm:af}}
\label{sec:af-neumann}\ 

By Lemma~\ref{lem:static-convex-domain}, $\Sigma$ is connected
and $h>0$. We first establish a boundary integral identity
which will be used later.

We denote by $\sigma_m(B_1,\ldots,B_m)$ the complete
polarization of $\sigma_m$ in dimension $n-1$, normalized by
$\sigma_m(B[m])=\sigma_m(B)$, where square brackets indicate
repeated arguments. The mixed Newton tensor $T_{m-1}$ is defined by
\begin{eqnarray*}
\sigma_m(B_1,\ldots,B_m)= \frac 1 m \< T_{m-1}(B_1,\ldots,B_{m-1}),B_m\> 
  .
\end{eqnarray*}
We write $T_{m-1}(h)=T_{m-1}(h[m-1])$.
Products of symmetric tensors below are understood as
compositions after raising one index with $g_\Sigma$.
In particular, $g_\Sigma$ represents the identity endomorphism.
Since $h$ is a Codazzi tensor, mixed Newton tensors with
arguments chosen from $h$ and $g_\Sigma$ are divergence-free.

\begin{lem}\label{lem:af-mixed-boundary}
For $2\le m\le n-1$ and $f\in C^3(\Sigma)$, there holds
\begin{eqnarray}\label{eq:af-mixed-boundary}
\begin{aligned}
 &m\int_\Sigma V
 \sigma_m\left(
 \nabla_\Sigma^2f-\frac fV\nabla_\Sigma^2V,
 \nabla_\Sigma^2f-\frac fV\nabla_\Sigma^2V,
 h[m-2]\right)d\sigma_g
 \\
 &\quad=
 \int_\Sigma V^3
 \left\langle
 T_{m-1}(h)\left(h-\frac{V_\nu}{V}g_\Sigma\right)
 \nabla_\Sigma\left(\frac fV\right),
 \nabla_\Sigma\left(\frac fV\right)
 \right\rangle d\sigma_g
 \ge0.
\end{aligned}
\end{eqnarray}
Equality holds if and only if $f/V$ is constant on $\Sigma$.
\end{lem}

\begin{proof}
For this calculation, denote
\begin{eqnarray*}
 D_{\alpha\beta}
 \coloneqq \nabla^\Sigma_{\alpha\beta}f
  -\frac fV\nabla^\Sigma_{\alpha\beta}V,
 \qquad
 Z_\alpha
 \coloneqq \nabla^\Sigma_\alpha f
  -\frac fV\nabla^\Sigma_\alpha V.
\end{eqnarray*}
Then $Z=V\nabla_\Sigma(f/V)$ and
\begin{eqnarray}\label{eq:af-boundary-derivative}
 D_{\alpha\beta}
 =\nabla^\Sigma_\alpha Z_\beta
  +Z_\alpha\nabla^\Sigma_\beta\log V.
\end{eqnarray}

Choose an extension $\widetilde f$ to a collar of $\Sigma$
such that $\widetilde f=f$ and $B_V\widetilde f=0$ on $\Sigma$.
The static equation and the tangential Hessian formula give
\begin{eqnarray*}
 \nabla_\Sigma^2V=Vg_\Sigma-V_\nu h,
 \qquad
 A[\widetilde f]_{\alpha\beta}=D_{\alpha\beta}.
\end{eqnarray*}
By the mixed boundary identity in the proof of
Lemma~\ref{lem:reilly},
\begin{eqnarray*}
 A[\widetilde f]_{\nu\alpha}
 =-\left(h_{\alpha\beta}
          -\frac{V_\nu}{V}g_{\alpha\beta}\right)Z^\beta.
\end{eqnarray*}
Consequently, the Codazzi property of $A[\widetilde f]$
and the Gauss formula imply
\begin{eqnarray}\label{eq:af-boundary-codazzi}
\begin{aligned}
 \nabla^\Sigma_\alpha D_{\beta\gamma}
 -\nabla^\Sigma_\beta D_{\alpha\gamma}
 =
 h_{\alpha\gamma}
 \left(h_{\beta\delta}
       -\frac{V_\nu}{V}g_{\beta\delta}\right)Z^\delta -
 h_{\beta\gamma}
 \left(h_{\alpha\delta}
       -\frac{V_\nu}{V}g_{\alpha\delta}\right)Z^\delta.
\end{aligned}
\end{eqnarray}

Set $\mathcal T=T_{m-1}(D,h[m-2])$.
Using the generalized Kronecker symbol in a local
orthonormal frame, we obtain
\begin{eqnarray}\label{eq:af-mixed-newton-divergence}
\begin{aligned}
 \nabla^\Sigma_\alpha\mathcal T^\alpha{}_\beta
 &=
 \frac{1}{2(m-1)!}
 \delta^{\alpha\gamma\alpha_2\cdots\alpha_{m-1}}
        _{\beta\delta\beta_2\cdots\beta_{m-1}}
 \left(
 \nabla^\Sigma_\alpha D_\gamma{}^\delta
 -\nabla^\Sigma_\gamma D_\alpha{}^\delta
 \right)
 h_{\alpha_2}{}^{\beta_2}\cdots
 h_{\alpha_{m-1}}{}^{\beta_{m-1}}
 \\
 &=
 -T_{m-1}(h)^\alpha{}_\beta
 \left(h_{\alpha\gamma}
       -\frac{V_\nu}{V}g_{\alpha\gamma}\right)Z^\gamma.
\end{aligned}
\end{eqnarray}
Indeed, the terms differentiating $h$ vanish by its
Codazzi property, and the remaining terms are evaluated
by \eqref{eq:af-boundary-codazzi}.

By \eqref{eq:af-boundary-derivative} and the symmetry
of $\mathcal T$,
\begin{eqnarray*}
 V\mathcal T^{\alpha\beta}D_{\alpha\beta}
 =\mathcal T^{\alpha\beta}
   \nabla^\Sigma_\alpha(VZ_\beta).
\end{eqnarray*}
Integrating by parts on the closed hypersurface $\Sigma$
and using \eqref{eq:af-mixed-newton-divergence}, we have
\begin{eqnarray*}
\begin{aligned}
 m\int_\Sigma V\sigma_m(D,D,h[m-2])d\sigma_g
 &=
 -\int_\Sigma VZ^\beta
   \nabla^\Sigma_\alpha\mathcal T^\alpha{}_\beta
   d\sigma_g
 \\
 &=
 \int_\Sigma V
 \left\langle
 T_{m-1}(h)\left(h-\frac{V_\nu}{V}g_\Sigma\right)Z,Z
 \right\rangle d\sigma_g.
\end{aligned}
\end{eqnarray*}
Substituting $Z=V\nabla_\Sigma(f/V)$ proves the identity
in \eqref{eq:af-mixed-boundary}.

The two endomorphisms in the last integrand commute.
If $\kappa_\alpha$ are the principal curvatures, the
eigenvalues of their product are
\begin{eqnarray*}
 \sigma_{m-1}(\kappa|\alpha)
 \left(\kappa_\alpha-\frac{V_\nu}{V}\right)>0.
\end{eqnarray*}
Thus equality holds precisely when
$\nabla_\Sigma(f/V)=0$. Since $\Sigma$ is connected,
this is equivalent to $f/V$ being constant.
\end{proof}

Now we are ready to prove Theorem~\ref{thm:af}.
\begin{proof}[Proof of Theorem~\ref{thm:af}]
Let $(u,\lambda)$ be the $2$-admissible solution of
\eqref{eqn:sigma-2-neumann} supplied by
Theorem~\ref{thm:existence}. The Newton--Maclaurin inequality gives $\tr A[u]\ge n$. Applying \eqref{eq:geometry-flux} and \eqref{eq:geometry-reilly}, and using
static convexity, we obtain  
\begin{eqnarray}\label{eq:af-initial-comparison}
 \lambda\int_\Sigma Vd\sigma_g
 \ge n\int_\Omega Vd\mu_g>0,
 \qquad
 \lambda^2\int_\Sigma VHd\sigma_g
 \le n\int_\Omega Vd\mu_g.
\end{eqnarray}
In particular, $\lambda>0$.

We begin by proving \eqref{eq:af-adjacent-analytic}.
For each $2\le m\le n-1$, consider the following
linear equation on $\Sigma$:
\begin{eqnarray}\label{eq:af-boundary-equation}
\begin{aligned}
 \dive_\Sigma\left(
 V^2T_{m-1}(h)\nabla_\Sigma\left(\frac fV\right)
 \right) =
 m\binom{n-1}{m}V
 \left(
 \frac{\displaystyle\int_\Sigma VH_{m-1}d\sigma_g}
      {\displaystyle\int_\Sigma VH_md\sigma_g}
 H_m-H_{m-1}
 \right).
\end{aligned}
\end{eqnarray}
The right hand side has zero integral on $\Sigma$.
Since $h>0$, the tensor $T_{m-1}(h)$ is positive definite,
and the operator in the variable $f/V$ is uniformly elliptic.
Its kernel consists of constants. Hence the Fredholm
alternative gives a solution $f$, unique up to a multiple of $V$.
The boundary regularity and Schauder estimates give
$f\in C^{3,\alpha_0}(\Sigma)$, which is sufficient for
Lemma~\ref{lem:af-mixed-boundary}.

For this auxiliary boundary solution, denote
\begin{eqnarray}\label{eq:af-primitive-tensor}
 B\coloneqq 
 \nabla_\Sigma^2f-\frac fV\nabla_\Sigma^2V
 +g_\Sigma
 -\frac{\displaystyle\int_\Sigma VH_{m-1}d\sigma_g}
       {\displaystyle\int_\Sigma VH_md\sigma_g} h.
\end{eqnarray}
The normalization of the polarized symmetric functions gives
\begin{eqnarray}\label{eq:af-polarization-identities}
 \sigma_m(h[m])
 =\binom{n-1}{m}H_m,\qquad 
 \sigma_m(g_\Sigma,h[m-1])
 =\binom{n-1}{m}H_{m-1},
\end{eqnarray}
and
\begin{eqnarray}\label{eq:af-polarization-identities-2}
     \sigma_m(g_\Sigma,g_\Sigma,h[m-2])
 =\binom{n-1}{m}H_{m-2}.
\end{eqnarray}
Since $T_{m-1}(h)$ is divergence-free,
\begin{eqnarray*}
\begin{aligned}
 \dive_\Sigma\left(
 V^2T_{m-1}(h)\nabla_\Sigma\left(\frac fV\right)
 \right) =
 mV\sigma_m\left(
 \nabla_\Sigma^2f-\frac fV\nabla_\Sigma^2V,h[m-1]
 \right).
\end{aligned}
\end{eqnarray*}
Therefore, \eqref{eq:af-boundary-equation} and
\eqref{eq:af-polarization-identities} imply
\begin{eqnarray}\label{eq:af-primitive-condition}
 \sigma_m(B,h[m-1])=0.
\end{eqnarray}
Applying the concavity of $\sigma_m^{1/m}$ to $h+tB$
at $t=0$, we obtain
\begin{eqnarray}\label{eq:af-polarized-concavity}
 \sigma_m(h)\sigma_m(B,B,h[m-2])
 \le \sigma_m(B,h[m-1])^2=0.
\end{eqnarray}
In particular, $\sigma_m(B,B,h[m-2])\le0$ pointwise.

The tensors $T_{m-1}(h)$ and
$T_{m-1}(g_\Sigma,h[m-2])$ are divergence-free.
The integration by parts used in the proof of
Lemma~\ref{lem:af-mixed-boundary} therefore also gives
\begin{eqnarray}\label{eq:af-boundary-orthogonality}
\begin{aligned}
 &\int_\Sigma V
 \sigma_m\left(
 \nabla_\Sigma^2f-\frac fV\nabla_\Sigma^2V,
 h[m-1]\right)d\sigma_g=0,
 \\
 &\int_\Sigma V
 \sigma_m\left(
 \nabla_\Sigma^2f-\frac fV\nabla_\Sigma^2V,
 g_\Sigma,h[m-2]\right)d\sigma_g=0.
\end{aligned}
\end{eqnarray}
Expanding \eqref{eq:af-primitive-tensor} and using
\eqref{eq:af-polarization-identities-2} and
\eqref{eq:af-boundary-orthogonality}, we have
\begin{eqnarray*}
\begin{aligned}
\int_\Sigma V\sigma_m(B,B,h[m-2])d\sigma_g
 &=
 \int_\Sigma V
 \sigma_m\left(
 \nabla_\Sigma^2f-\frac fV\nabla_\Sigma^2V,
 \nabla_\Sigma^2f-\frac fV\nabla_\Sigma^2V,
 h[m-2]\right)d\sigma_g
 \\
 &\qquad+
 \binom{n-1}{m}
 \left[
 \int_\Sigma VH_{m-2}d\sigma_g
 -
 \frac{\left(\displaystyle\int_\Sigma
                    VH_{m-1}d\sigma_g\right)^2}
      {\displaystyle\int_\Sigma VH_md\sigma_g}
 \right].
\end{aligned}
\end{eqnarray*}
Substituting \eqref{eq:af-mixed-boundary} and rearranging
the terms yields the exact identity
\begin{eqnarray}\label{eq:af-adjacent-deficit}
\begin{aligned}
 &\left(\int_\Sigma VH_{m-1}d\sigma_g\right)^2
 -\left(\int_\Sigma VH_{m-2}d\sigma_g\right)
  \left(\int_\Sigma VH_md\sigma_g\right)
 \\
 &\quad=
 \frac{\displaystyle\int_\Sigma VH_md\sigma_g}
      {m\binom{n-1}{m}}
 \int_\Sigma V^3
 \left\langle
 T_{m-1}(h)\left(h-\frac{V_\nu}{V}g_\Sigma\right)
 \nabla_\Sigma\left(\frac fV\right),
 \nabla_\Sigma\left(\frac fV\right)
 \right\rangle d\sigma_g
 \\
 &\qquad+
 \frac{\displaystyle\int_\Sigma VH_md\sigma_g}
      {\binom{n-1}{m}}
 \int_\Sigma V
 \left[-\sigma_m(B,B,h[m-2])\right]d\sigma_g.
\end{aligned}
\end{eqnarray}
Both terms on the right hand side are nonnegative by
Lemma~\ref{lem:af-mixed-boundary} and
\eqref{eq:af-polarized-concavity}.
Taking $m=j+1$ proves \eqref{eq:af-adjacent-analytic}.

Suppose equality holds in
\eqref{eq:af-adjacent-analytic} for some $j=m-1$.
The first term in \eqref{eq:af-adjacent-deficit} vanishes,
so Lemma~\ref{lem:af-mixed-boundary} gives
$f/V=\mathrm{constant}$. Hence
\begin{eqnarray*}
 \nabla_\Sigma^2f-\frac fV\nabla_\Sigma^2V=0.
\end{eqnarray*}
The second term in \eqref{eq:af-adjacent-deficit} also
vanishes, and \eqref{eq:af-polarized-concavity} gives
$\sigma_m(B,B,h[m-2])=0$ pointwise.
Combining this with \eqref{eq:af-primitive-condition}
and \eqref{eq:af-polarization-identities}, we obtain
\begin{eqnarray*}
 H_{m-1}^2=H_{m-2}H_m.
\end{eqnarray*}
Since $h>0$, equality in the Newton--MacLaurin inequality implies
that all principal curvatures coincide.
Moreover, \eqref{eq:af-primitive-condition} gives
\begin{eqnarray*}
 h=
 \frac{\displaystyle\int_\Sigma VH_md\sigma_g}
      {\displaystyle\int_\Sigma VH_{m-1}d\sigma_g}
 g_\Sigma.
\end{eqnarray*}
Thus $\Sigma$ is compact, connected, and totally umbilical
with constant principal curvature. It follows that
$\Omega$ is a geodesic ball.

We now combine the adjacent-order inequalities with
the estimates for the solution $u$ in Eq. \eqref{eqn:sigma-2-neumann}.
Since all weighted curvature integrals are positive,
\eqref{eq:af-adjacent-analytic} gives
\begin{eqnarray*}
 \frac{\displaystyle\int_\Sigma VH_{j+1}d\sigma_g}
      {\displaystyle\int_\Sigma VH_jd\sigma_g}
 \le
 \frac{\displaystyle\int_\Sigma VH_jd\sigma_g}
      {\displaystyle\int_\Sigma VH_{j-1}d\sigma_g},
 \qquad 1\le j\le n-2.
\end{eqnarray*}
Multiplying the resulting ratio estimates yields
\begin{eqnarray}\label{eq:af-curvature-ratios}
 \int_\Sigma VH_jd\sigma_g
 \le
 \frac{\left(\displaystyle\int_\Sigma VHd\sigma_g\right)^j}
      {\left(\displaystyle\int_\Sigma Vd\sigma_g\right)^{j-1}},
 \qquad 1\le j\le n-1.
\end{eqnarray}
For $j=1$ this is an identity, which also covers $n=2$.

Using \eqref{eq:af-initial-comparison} in
\eqref{eq:af-curvature-ratios}, we obtain
\begin{eqnarray*}
 \lambda^{j+1}\int_\Sigma VH_jd\sigma_g
 \le
 \frac{\lambda^{j+1}
       \left(\displaystyle\int_\Sigma VHd\sigma_g\right)^j}
      {\left(\displaystyle\int_\Sigma Vd\sigma_g\right)^{j-1}}
\le
 n\int_\Omega Vd\mu_g
 \left(
 \frac{n\displaystyle\int_\Omega Vd\mu_g}
      {\lambda\displaystyle\int_\Sigma Vd\sigma_g}
 \right)^{j-1}
\le n\int_\Omega Vd\mu_g.
\end{eqnarray*}
Together with the first inequality in
\eqref{eq:af-initial-comparison}, this proves
\begin{eqnarray}\label{eq:af-boundary-constant-bounds}
 \frac{n\displaystyle\int_\Omega Vd\mu_g}
      {\displaystyle\int_\Sigma Vd\sigma_g}
 \le\lambda
 \le
 \left(
 \frac{n\displaystyle\int_\Omega Vd\mu_g}
      {\displaystyle\int_\Sigma VH_jd\sigma_g}
 \right)^{1/(j+1)},
 \qquad 1\le j\le n-1.
\end{eqnarray}
Eliminating $\lambda$ gives \eqref{eq:af-two-hessian-bounds}.
Raising both sides to the power $j+1$ and rearranging
gives \eqref{eq:global-weighted-AF}.
Its unnormalized form \eqref{eq:global-weighted-AF-unnormalized} follows from
$H_j=\binom{n-1}{j}^{-1}\sigma_j(h)$.

If equality holds in \eqref{eq:global-weighted-AF}
for some $j$, the two endpoints in
\eqref{eq:af-boundary-constant-bounds} coincide. Hence
\begin{eqnarray*}
 \lambda\int_\Sigma Vd\sigma_g
 =n\int_\Omega Vd\mu_g.
\end{eqnarray*}
Equations~\eqref{eq:geometry-flux} and
\eqref{eq:geometry-bulk-energy} imply that $A[u]=g$.
Lemma~\ref{lem:geometric-rigidity} then shows that
$\Omega$ is a geodesic ball.

Conversely, on a geodesic ball of radius $R$, the reference solution \eqref{eq:ball-reference-solution} 
satisfies $A[u]=g$ and $\lambda=\tanh R$.
Since $H_j=\coth^jR$ and \eqref{eq:geometry-flux} gives
\begin{eqnarray*}
 n\int_\Omega Vd\mu_g
 =\lambda\int_\Sigma Vd\sigma_g,
\end{eqnarray*}
all the stated inequalities are equalities.
The center of the ball need not coincide with the pole
defining $V$. This completes the proof.
\end{proof}

\bigskip

\bigskip

\noindent\textbf{Acknowledgment.} 
The author expresses his gratitude to Chao Xia for many enlightening discussions on this project over the years. The author used AI tools to assist with language editing and to identify points requiring mathematical and bibliographic
verification. The author checked all mathematical arguments, references, and revisions and take full responsibility for the
final manuscript.

\printbibliography[heading=bibliography]
\end{document}